\documentclass[10pt,letterpaper]{amsart}

\usepackage[T1]{fontenc}
\usepackage[utf8]{inputenc}
\usepackage[english]{babel}
\usepackage[letterpaper,top=1.5cm,bottom=1.5cm,left=1.5cm,right=1.5cm,marginparwidth=2cm]{geometry}
\usepackage{microtype}
\usepackage{amsmath,amssymb,amsfonts,amsthm,mathtools,mathrsfs,bm}
\usepackage{graphicx}
\usepackage{tikz-cd}

\usepackage{enumitem}
\usepackage{cite}
\usepackage[colorlinks=true,linkcolor=blue,citecolor=blue,urlcolor=blue]{hyperref}

\newlist{steps}{enumerate}{1}
\setlist[steps,1]{label=\textbf{Step \arabic*.},leftmargin=*}
\numberwithin{equation}{section}

\makeatletter
\renewcommand{\section}{\@startsection{section}{1}{0pt}{-3.5ex plus -1ex minus -.2ex}{2.3ex plus .2ex}{\normalfont\Large\bfseries}}
\renewcommand{\subsection}{\@startsection{subsection}{2}{0pt}{-3.25ex plus -1ex minus -.2ex}{1.5ex plus .2ex}{\normalfont\large\bfseries}}
\renewcommand{\subsubsection}{\@startsection{subsubsection}{3}{0pt}{-3.25ex plus -1ex minus -.2ex}{1.5ex plus .2ex}{\normalfont\normalsize\bfseries}}
\providecommand{\xleftrightarrow}[2][]{\ext@arrow 3359\leftrightarrowfill@{#1}{#2}}
\providecommand{\xdashrightarrow}[2][]{\ext@arrow 0359\rightarrowfill@@{#1}{#2}}
\providecommand{\xdashleftarrow}[2][]{\ext@arrow 3095\leftarrowfill@@{#1}{#2}}
\providecommand{\xdashleftrightarrow}[2][]{\ext@arrow 3359\leftrightarrowfill@@{#1}{#2}}
\def\rightarrowfill@@{\arrowfill@@\relax\relbar\rightarrow}
\def\leftarrowfill@@{\arrowfill@@\leftarrow\relbar\relax}
\def\leftrightarrowfill@@{\arrowfill@@\leftarrow\relbar\rightarrow}
\def\arrowfill@@#1#2#3#4{  $\m@th\thickmuskip0mu\medmuskip\thickmuskip\thinmuskip\thickmuskip
   \relax#4#1
   \xleaders\hbox{$#4#2$}\hfill
   #3$}
\makeatother
\title{Surfaces with canonical hyperplane sections and their moduli spaces}
\author[TRAN BAO]{Phan Quoc Bao Tran}
\address{Mathematisches Institut \\University of Bern \\
Bern, Switzerland}
\email{phan.tran@unibe.ch}
\date{}

\DeclareMathOperator{\Hom}{Hom}

\DeclareMathOperator{\id}{id}

\DeclareMathOperator{\depth}{depth}
\DeclareMathOperator{\rank}{rank}

\DeclareMathOperator{\Aut}{Aut}
\DeclareMathOperator{\Jac}{Jac}
\DeclareMathOperator{\Supp}{Supp}

\DeclareMathOperator{\Codim}{codim}
\DeclareMathOperator{\Pic}{Pic}

\DeclareMathOperator{\Div}{Div}

\DeclareMathOperator{\mult}{mult}

\DeclareMathOperator{\Ab}{Ab}
\DeclareMathOperator{\Bir}{Bir}
\DeclareMathOperator{\Cl}{Cl}
\DeclareMathOperator{\Def}{Def}
\DeclareMathOperator{\Diff}{Diff}
\DeclareMathOperator{\For}{For}
\DeclareMathOperator{\Lin}{Lin}
\DeclareMathOperator{\SLCY}{SLCY}
\DeclareMathOperator{\Sslc}{S_{slc}}

\DeclareMathOperator{\can}{can}
\DeclareMathOperator{\pr}{pr}
\DeclareMathOperator{\prox}{prox}

\DeclareMathOperator{\ter}{ter}

\DeclareMathOperator{\Sym}{Sym}
\DeclareMathOperator{\Spec}{Spec}
\DeclareMathOperator{\Proj}{Proj}
\DeclareMathOperator{\Bl}{Bl}
\DeclareMathOperator{\Bs}{Bs}

\newcommand{\F}{\mathbb{F}}
\newcommand{\PP}{\mathbb{P}}
\newcommand{\OO}{\mathcal{O}}

\newcommand{\Q}{\mathbb{Q}}
\newcommand{\Z}{\mathbb{Z}}
\newcommand{\NN}{\mathbb Z_{\ge 0}}
\DeclareMathOperator{\vol}{vol}
\newcommand{\slc}{\mathrm{slc}}
\DeclareMathOperator{\lc}{lc}

\subjclass[2020]{Primary 14J10; Secondary 14D22, 14E08,
14E30, 14J17, 14J26, 14J27.}
\hypersetup{
pdftitle={Surfaces with canonical hyperplane sections and their moduli spaces},
pdfsubject={Algebraic geometry},
pdfauthor={Phan Quoc Bao Tran},
pdfkeywords={surfaces with canonical hyperplane sections, log Calabi-Yau surfaces, KSBA moduli, stable minimal models}
}

\theoremstyle{plain}
\newtheorem{theo}{Theorem}[section]
\newtheorem{coro}[theo]{Corollary}
\newtheorem{lemm}[theo]{Lemma}
\newtheorem{prop}[theo]{Proposition}

\theoremstyle{definition}
\newtheorem{defi}[theo]{Definition}
\newtheorem{exam}[theo]{Example}
\newtheorem{conj}[theo]{Conjecture}
\newtheorem{ques}{Question}
\newtheorem{construction}{Construction}
\theoremstyle{remark}
\newtheorem{rema}[theo]{Remark}

\thanks{The author was funded by Schweizerischer Nationalfonds  [10005432] and [229554].}

\begin{document}
\maketitle

\begin{abstract}
We construct stable log Calabi-Yau surfaces from anticanonical pairs and marked clusters on rational and elliptic ruled surfaces. First, we blow up the cluster. Then we contract by the polarization. Sometimes the cluster resolves the whole intersection of the polarization with the boundary. In this case, the boundary is contracted, and we get a log canonical surface with canonical hyperplane sections. The surfaces of a fixed construction type form locally closed strata in families. Each stratum has a quasi-projective coarse moduli space. With this description, we prove that these moduli spaces are unirational when the boundary is nodal. Under a transversality condition, we get the same result when the boundary consists of smooth elliptic curves. We also study the Hacking--Keel--Yu moduli space of smooth polarized anticanonical pairs. Next, we study Birkar's moduli spaces of these models with fibre volume $e$. For them, we give sharp lower bounds for the self-intersection $\alpha$ of the polarization. Every integer $\alpha\ge-3e$ occurs.
\end{abstract}

\tableofcontents

\section{Introduction}\label{sec1}

The moduli theory for algebraic objects started from Mumford, who constructed the coarse moduli space $M_g$ of smooth curves as a quasi-projective variety using geometric invariant theory \cite{MumfordGIT}. Later, Deligne and Mumford introduced stable curves and the compactification $\overline{M}_g$ \cite{DM69}. Then, Mumford and Gieseker proved that $\overline{M}_g$ is projective by GIT \cite{Mum77,Gie82}. Continue from that, Knudsen and Mumford gave another proof with determinant line bundles \cite{KnM76,Knu83}. For surfaces of general type, Gieseker used GIT to construct the coarse moduli space of canonical models \cite{Gie77}, which is a quasi-projective variety. To compactify this space, Koll\'ar and Shepherd-Barron introduced stable surfaces \cite{KSB} and built an appropriate moduli theory for this new class of surfaces. From that, Alexeev proved the boundedness of stable surfaces \cite{Alexeev94,Alexeev96}. Notice that this setting is close to the moduli of polarized manifolds with semiample canonical sheaf. For these manifolds, Viehweg constructed quasi-projective moduli spaces \cite{Vie95}. His method combines GIT with the positivity of direct images. From that, Koll\'ar proved an ampleness lemma, and with it he proved the projectivity of complete moduli spaces \cite{Kol90}. Today, the theory also covers stable pairs in higher dimensions \cite{HMX1,KP17,Kol2,FI21,ABIP21} and the stable minimal models of Birkar \cite{Bir1,Bir2}.

In this paper, we look at the objects of some of these moduli spaces. We want to get them by explicit birational constructions. For this reason, we start with the following question.
\begin{ques}\label{querecon}
We consider a moduli problem for projective pairs $((X,\Delta),C)$ with some given numerical information. Here $C$ is an effective semiample Cartier divisor, and $(X,\Delta+\epsilon C)$ is log canonical for all sufficiently small rational $\epsilon>0$. Can we get these objects from pairs on minimal models of an explicitly classified type? Here we only allow explicit blow-ups and contractions, and finitely many discrete data.
\end{ques}

Our starting point is a normal Gorenstein surface $S$ with trivial canonical bundle and an embedding $S\hookrightarrow\mathbb P^g$ given by the complete system $|C|$. Here a general member of $|C|$ is a smooth curve of genus $g$. By adjunction, we have $\mathcal O_C(C)\simeq\omega_C$, and this canonical bundle is very ample. So the hyperplane sections are canonically embedded curves. Epema studied these surfaces in \cite{Epe}. He took their minimal resolutions and contracted $(-1)$-curves down to rational or relatively minimal ruled surfaces. During this process, he kept track of the polarization and of the blow-up centres, and in this way he classified these surfaces.

We recall the definition of \cite[Definition I.1.1]{Epe}.

\begin{defi}\label{defchs}
A projective surface $S$ is called a surface with canonical hyperplane sections if one of the following holds:
\begin{enumerate}
    \item[a)] There exists an embedding
    \begin{center}
    $i:S\hookrightarrow \mathbb P^g,\ g\geq 3,$
    \end{center}
    such that a general hyperplane section $C\subset \mathbb P^{g-1}$ is a smooth non-hyperelliptic curve of genus $g$ embedded by its complete canonical linear system. Equivalently, $\mathcal O_C(1)\simeq\omega_C.$
    \item[b)] $S$ contains a linear system $\mathcal{L}$ with $\dim \mathcal{L} =g$, $g \geq 2$, a general element of which is a smooth hyperelliptic curve $C$ of genus $g$ such that $\mathcal{L}|_C =|K_C|$ and the map $\phi_L: S \longrightarrow \mathbb{P}^g$ is a finite morphism of degree 2 onto a surface $\overline{S}$ in $\mathbb{P}^g$.
\end{enumerate}

\end{defi}

Epema's classification separates the $K3$ case from the rational and ruled cases that we study.

\begin{theo}[{\cite[Theorem 8.1]{CDS20}}, after {\cite{Epe}}]\label{threpema}
For a non-degenerate surface $S\subseteq\mathbb P^g$ with a canonical hyperplane section of genus $g\geq3$, we have the following four possibilities in Epema's classification:
\begin{enumerate}
    \item[1)] a $K3$ surface with canonical singularities;
    \item[2)] a rational surface with one minimally elliptic singularity;
    \item[3)] a ruled surface over a genus-$q$ curve, $q\geq1$, with one singularity of genus $q+1$;
    \item[4)] a ruled surface over an elliptic curve with two simple elliptic singularities.
\end{enumerate}
In cases 2)--4), the surface $S$ may also have canonical singularities.
\end{theo}

In this paper, we consider the log canonical surfaces in cases 2) and 4). They belong to the \textit{fake $K3$ surfaces} of \cite{CDS20}. By Epema \cite{Epe}, we have $K_S\sim0$ for these surfaces. We recall the numerical identities on their relatively minimal ruled surfaces in Subsection \ref{sec26}.

Epema also studied moduli and period maps for several special families \cite{Epe}. Ciliberto and Lopez used his classification to bound the number of moduli of extendable canonical curves \cite{CilibertoLopez02}. One can find an account of the classification in \cite{CDS20}. Related questions on hyperplane sections of $K3$ surfaces are studied in \cite{ABS17}.

Let $(S,C)$ be a normal log canonical surface with canonical hyperplane sections which is not a $K3$ surface. Let $f:X\to S$ be its minimal resolution. We denote by $W$ the reduced exceptional divisor over the non-canonical points of $S$. We still write $C$ for the pullback of the polarization. With this notation, we have
\begin{center}
$K_X+W=f^*K_S\sim0.$
\end{center}
Here the exceptional curves over the Du Val points do not appear in $W$. Now we run the $K_X$-minimal model program. This means that we contract $(-1)$-curves one by one. At the end, we get a rational surface or a relatively minimal ruled surface $X_{\min}$. Let $E$ be a $(-1)$-curve which is not contained in $W$. By adjunction, we have that $W\cdot E=-K_X\cdot E=1$. So the blow-down centre of $E$ lies on the image of the boundary. The same holds at each step. Let $\phi:X\to X_{\min}$ be the morphism we get in this way. We put $W_1=\phi_*W$ and $\mathcal L_1=\phi_*|C|$. We denote by $C_1$ a general member of $\mathcal L_1$. In general, the system $\mathcal L_1$ has base points, proper or infinitely near. At these points, its members can be singular, or they can be tangent to $W_1$. The inverse blow-ups resolve these base points. For each base point, we record its multiplicity. We also record the proximity relations between these points. Let $p$ be an earlier point, and let $p'$ be a later point. We say that $p'$ is proximate to $p$ if $p'$ lies on the strict transform of the exceptional curve of $p$. In this way, we get a weighted cluster in the sense of \cite{Roe04,KPT11} and \cite[Chapters 3--4]{Cas}. We recall this notation in Definitions \ref{defprox} and \ref{defmult}. In this situation, we can see that the centres lie on the strict transform of $W_1$ at each stage. Hence, we can label each centre by the boundary component that contains it. We call these data \textit{marked cluster}, cf. Definition \ref{defmark}.

Our aim is to reverse the construction of Epema with these clusters. In this direction, we allow the non-emptiness of the boundary. So, we use the following notion. We call $((Y,D),C)$ a \textit{stable log Calabi-Yau surface} if three conditions hold: First, the pair $(Y,D)$ is a projective log canonical surface pair with $K_Y+D\sim0$. Second, the divisor $C$ is an effective ample Cartier divisor. Third, the pair $(Y,D+\epsilon C)$ is log canonical for all sufficiently small rational $\epsilon>0$. In this case, the divisor $K_Y+D+\epsilon C\sim_{\mathbb Q}\epsilon C$ is ample. Also, the pullback of $C$ to a resolution is nef and big.

Now we ask which divisors and which marked clusters on $(X_{\min},W_1)$ give back $((Y,D),C)$, and with which intersection numbers. We also ask how this works in families.

We start with a smooth surface $X_{\min}$ and a reduced anticanonical divisor $W_1\in|-K_{X_{\min}}|$. We choose $C_1\in|aC_\infty|$ with $a\geq2$. We call $X_{\min}$ the \textit{initial surface} and $(X_{\min},W_1)$ the \textit{initial pair}. In the ruled cases, the section $C_\infty$ is disjoint from the negative section $C_0$. We have $C_0^2=-e$ and $C_\infty^2=e>0$. In the plane case, $C_\infty$ is a line. We describe these section classes in Subsection \ref{sec21}. We denote by $q$ the genus of the base of the ruling. For the plane, we put $q=-\infty$. We also put $e=C_\infty^2$. Next, we choose a marked cluster $\mathcal K$ with multiplicities $m_\alpha$ in $C_1|_{W_1}=\Sigma+R$. Here $\Sigma$ is the part which the cluster imposes and which its blow-ups resolve. The part $R$ is the intersection which is left on the boundary. When the initial pair and the construction satisfy Definition \ref{deftuple}, we call this an \textit{admissible tuple}. By blowing up the cluster, we get $K_X+W\sim0$. When $\Sigma=C_1|_{W_1}$ and $R=0$, we also have $C\cap W=\varnothing$. Then the polarization gives the contraction of $W$. When we fix $\Sigma+R$ before we choose $C_1$, we require $\mathcal O_{W_1}(\Sigma+R)\simeq\mathcal O_{W_1}(C_1)$. We call this the \textit{Abel--Jacobi--Picard condition}. Together with the surjectivity of the restriction on sections, it allows a member of $|C_1|$ to cut out the given divisor $\Sigma+R$. For $R=0$, resolving this intersection gives $\mathcal O_X(C)|_W\simeq\mathcal O_W$. This condition is needed, as Hacon and Jiang show that a non-torsion restriction to an elliptic curve can prevent its contraction to a projective surface \cite[Remark 4.3 and Lemma 4.4]{HaconJiang25}. From this, we introduce our main theorem

\begin{theo}\label{thrconst}\label{thrchs}
Every admissible tuple gives a stable log Calabi-Yau surface $((Y,D),C)$ and a surface $S$ with CHS.
\end{theo}

In this construction, the initial surface and its anticanonical boundary are only extra choices. For example, we can change a ruled initial surface by an elementary transformation. This transformation blows up a point and then contracts the strict transform of the fibre through this point, see Subsection \ref{sec21}. When the boundary is rational and reducible, we choose a taut anticanonical pair. For such a pair, the ordered self-intersection numbers determine the labelled pair up to isomorphism, see Definition \ref{deftaut}. In the rational case, we explan the inverse direction in Remark \ref{remrecon}. Here, we assume that the resolution of the CHS surface has a sequence of blow-downs, which ends at an initial pair (cf. Definition \ref{deftuple}). For surfaces with CHS with two simple elliptic singularities, we will see that Proposition \ref{prop2ell} gives the inverse construction from elliptic ruled surfaces.

From the point of view of $K3$ surfaces, it is also natural to look at fibration structures on surfaces with CHS. This leads to the following question. Is there a non-$K3$ surface with CHS which carries a pencil of curves of genus one? We want the general member of this pencil to be smooth. Here a pencil means a linear system of dimension one. To construct such surfaces, we take a pencil $\mathcal P\subset|M|$, where $M=-K_{X_{\min}}$ or $M=C_\infty$. We require that all proper and infinitely near base centres of $\mathcal P$ occur in $\mathcal K$. The multiplicities of the pencil at these centres are given separately from those of the polarization. If a base divisor has repeated points, then its resolution also contains the infinitely near centres over these points, as in Subsection \ref{sec33}. We assume that every curve contracted by the polarization is vertical for the resolved pencil. In this case, the pencil descends to a morphism $f_S:S\to\mathbb P^1$. When $\mathcal P$ is an anticanonical pencil on a rational surface, or a positive-section pencil on an elliptic ruled surface, this morphism is an elliptic fibration. For the case that a pencil in $|C_\infty|$ on $\mathbb F_e$, it gives a rational fibration. Thus, we have the following theorem which gives a positive answer to the question above for every genus $g\geq3$.
\begin{theo}\label{thrfib}
For every integer $g\geq3$, there is a surface $S$ with CHS of genus $g$ which has exactly two simple elliptic singularities and an elliptic fibration $f:S\to\mathbb P^1$.
\end{theo}
All these data can be parametrized. So it is natural to ask the following question about families.
\begin{ques}\label{quefam}
Can we extend the reconstruction to families? Can we identify its image in a moduli space of stable minimal models or of KSBA stable pairs, and give this image a moduli interpretation?
\end{ques}
Here a KSBA stable pair is a projective slc pair $(X,\Delta)$ with $K_X+\Delta$ ample, see Definition \ref{defpair}. We refer to \cite{Jiang26} for normalized moduli in the polarized klt case.

Now we fix some notation. Let $\Upsilon$ be the tuple of some important information from the initial anticanonical pair. These include the polarization, the marked cluster and its numerical invariants, cf. Definition \ref{deftuple} and Subsection \ref{sec42}. For a fixed $\Upsilon$, the image of the construction is a finite union of locally closed pieces, with the indexes are $ \nu \in V$ with $V$ finite. In the decomposable elliptic ruled case, this image is locally closed. Here, we work on each chosen piece indexing by $\nu$. So We drop its index for simplicity, through this is still true when we work with $q=1$
\begin{theo}\label{thrcoarse}
The class of stable log Calabi-Yau surfaces represented by $\Upsilon$ has a quasi-projective coarse moduli space.
\end{theo}
We denote this coarse moduli space by $\mathrm{SLCY}_{\mathrm{lc}}\Upsilon$. By Theorem \ref{thrmoduli}, its geometric points are the isomorphism classes of these surfaces. Now, let $q=1$ and $R=0$. In the condition that the multiplicities of singular points on $C_1$ of $X_{\mathrm{min}}$ satisfy $0<m_\alpha<a$, the stratums with a chosen admissible pencil also have quasi-projective coarse moduli spaces. These spaces are finite over the coarse spaces of the surface strata below them. This follows from Theorem \ref{thrpenmod} and Lemma \ref{lempencil}. By Proposition \ref{propelmod}, the same holds for the stable minimal models constructing from these elliptic fibrations.

The coarse moduli spaces of Section \ref{sec4} connect two things. The first one is the parameter space of initial pairs and marked clusters. The second one is the set of isomorphism classes of the surfaces which we construct. More exactly, the stable construction locus $I^{\mathrm{st}}_\Upsilon$ of Definition \ref{defstable} parametrizes the admissible tuples which give stable pairs. The construction gives a dominant rational map
\begin{center}
$ I^{\mathrm{st}}_\Upsilon\dashrightarrow \mathrm{SLCY}_{\mathrm{lc}}\Upsilon$
\end{center}
to the coarse moduli stratum. Two choices of initial surface and extra markings can give isomorphic stable pairs. Let $\mathrm M(\Upsilon)$ be the rational quotient which identifies such choices. By Definition \ref{defquot}, its function field consists of the rational functions on the construction space which do not change under these choices. For the birational statements, we take the dense open piece of each irreducible construction image.

One of our motivations is the conjecture of Hacking, Keel and Yu on the unirationality of moduli spaces of anticanonical pairs \cite[Conjecture 1.1]{HKY}. Recall that for a normal crossing divisor with smooth components, a $0$-stratum is a point which is a connected component of an intersection of boundary components. On a surface, a $0$-stratum is a node of the boundary. We describe the surface setting in Subsection \ref{sec25}.
\begin{conj}\label{conjhky}
Any connected component of the moduli space of triples $(X,E,\Theta)$ with the following properties is unirational.
\begin{enumerate}
    \item $X$ is a connected smooth projective complex variety.
    \item $E=E_1+\cdots+E_n\in|-K_X|$ is a normal crossing divisor with a $0$-stratum, and every $E_i$ is smooth.
    \item $\Theta$ is ample and contains no $0$-stratum of $E$.
\end{enumerate}
\end{conj}

In \cite[Theorem A]{AAB}, Alexeev, Arg\"uz and Bousseau prove that every irreducible component of the coarse KSBA moduli space of stable log Calabi-Yau surfaces has a finite surjective morphism from a projective toric variety. This proves \cite[Conjecture 1.2]{HKY} for surfaces. In particular, these components are unirational, and Conjecture \ref{conjhky} holds for surfaces. Our construction gives another proof of the unirationality, see Proposition \ref{prophky}. We use \cite[Lemma 2.8]{AAB} to put a general polarized pair in a good toric model. Our proof also works without a $0$-stratum. Indeed, we also include smooth elliptic boundaries and disjoint boundaries. For them, we need a transversality condition and a cluster condition. Let us state them. Let $W_1$ be a disjoint union of smooth elliptic curves. We let $W_1\in|-K_{X_{\min}}|$ and smooth $C_1\in|aC_\infty|$ vary in their complete linear systems. We assume that
\begin{center}
$ C_1|_{W_1}=p_1+\cdots+p_N.$
\end{center}
Here the $p_i$ are distinct points, and $C_1$ meets $W_1$ transversely. The cluster $\mathcal K$ consists of these $N$ distinct points with multiplicity one. So we have $R=0$.

\begin{theo}\label{thrunirat}
Let $a\geq2$. For the nodal-boundary tuples, $\mathrm M(\Upsilon)$ and $\mathrm{SLCY}_{\mathrm{lc}}\Upsilon$ are unirational. The same holds for the smooth elliptic boundaries above.
\end{theo}
Notice that we must impose the condition on distinct transverse intersection points here. Indeed, Example \ref{extang} gives an elliptic construction stratum which dominates $M_{1,11}$, and $M_{1,11}$ is not unirational.

Beside the elliptic ruled case, the rational nodal case is Theorem \ref{thrnodal} and the rational smooth-boundary case follows from Proposition \ref{propsmbd} for rational initial surfaces and from Construction \ref{constr5}.

The construction also covers canonical modifications of surfaces with CHS. Their boundary has no $0$-strata, and their polarization is semiample. For their unirationality, we need the distinct transverse intersection points and the proper centres of multiplicity one, as in Theorem \ref{thrunirat}. Notice that our approach gives a criterion for the rationality of some of KSBA irreducible components by constructing a dominant rational map with the quotient description. We explain this now. For the smooth rational pairs in \cite[Conjecture 1.1]{HKY}, the ampleness gives a positive polarization degree on every boundary component. Let $D=D_1+\cdots+D_r$. Gross, Hacking and Keel define a boundary marking \cite[Definition 1.11]{GHK1}. It is a choice of a smooth point on each component of the boundary. The group $\Aut^0(D)\simeq(\mathbb G_m)^r$ acts on these choices. By \cite[Lemma 2.1]{GHK1}, it rescales the coordinates on the curves $D_i^{\mathrm{sm}}$ one by one. A general divisor $C|_D$ has trivial stabilizer under this action. We will use this fact to shorten the computation of the quotient. A Picard marking identifies a fixed lattice $P$ with $\Pic(X)$. This identification respects the intersection product, the class of $L$ and the classes of the $D_i$. Let $g$ be a change of marking. It replaces the identification $\mu$ by $\mu\circ g^{-1}$, and it permutes the boundary labels. We only keep the changes which come from presentations in the class $\Upsilon$. For these presentations, we identify the boundary and the polarizing systems on a common resolution. These actions on $P$ and on the boundary labels give us a finite group, which we denote by $F_\Upsilon$.

Now, let $\alpha\in P$ be a class with $\alpha^2=-2$. We assume that $\alpha\cdot L=0$ and $\alpha\cdot D_i=0$ for every $i$. Then the reflection in $\alpha$ is
\begin{center}
$ s_\alpha(\lambda)=\lambda+(\lambda\cdot\alpha)\alpha.$
\end{center}
Let $W_\Upsilon$ be the group generated by the reflections $s_\alpha$ in $F_\Upsilon$. We put $H_\Upsilon=F_\Upsilon/W_\Upsilon$. At first, we take the quotient by $\Aut^0(D)$. after that, we take the invariants under $W_\Upsilon$. After these two steps, the group $H_\Upsilon$ acts monomially on a torus. If $|H_\Upsilon|\leq2$, Then this action is trivial or a monomial involution. In both cases, Lemma \ref{lemmoninv} shows that the invariant field is rational. For the section bundle, we use Lemma \ref{lemaffine}. This is the reason for the bound on the order in the following theorem. Let $M_0$ be an irreducible component of the moduli space of smooth rational pairs $(X,D)$ with an ample polarization $C$. Here $D$ is a normal crossing anticanonical boundary with smooth components and a $0$-stratum.
\begin{theo}\label{thrrat}
If the general member of $M_0$ has $|H_\Upsilon|\leq2$, then $M_0$ is rational.
\end{theo}
Here $(X,D)$ varies in its deformation type, and the polarizing curve varies in the complete linear system. See Theorem \ref{thrsmooth}. Corollary \ref{coraab} gives the same conclusion for KSBA components with dense smooth locus. We will treat the cases with $|H_\Upsilon|>2$ in a separate paper.

A stable minimal model is a triple $(X,B),A$ with the following properties. The variety $X$ is projective, and $A$ is an effective $\mathbb Q$-divisor. The divisor $K_X+B$ is semiample. For every sufficiently small $t>0$, the pair $(X,B+tA)$ is slc with ample log canonical divisor, see Definition \ref{defmodel}. At a simple base point of the pencil $\mathcal P$, we first blow up once. Then we blow up the point where the first exceptional curve meets the boundary. The strict transform of the first exceptional curve is a $(-2)$-section. After one more blow-up, these sections give the polarizing curves of Examples \ref{exell}, \ref{exruled} and \ref{exsharp}. From an elliptic fibration, Construction \ref{constr3} gives stable minimal models of Kodaira dimension one. Construction \ref{constr4} treats rational fibrations with $K+B\sim\beta F$, where $F$ is a general fibre and $\beta>0$. Its boundary can be the union of the CHS boundary with some chosen fibres, or it can be a general member of the same system.

Birkar defines moduli spaces for these stable minimal models \cite{Bir2}. Here the fibration is the one defined by $K_X+B$. He asks which numerical tuples occur \cite[Question 11.13]{Bir2}. We denote by $\mathcal V_{\mathrm{fib}}$ the set of permitted volumes on general fibres of the Iitaka fibration. This set is called $\Gamma$ in \cite{Bir2}. We denote by $S_{\mathrm{slc}}(d,\Phi,\mathcal V_{\mathrm{fib}},\sigma)$ the set of stable minimal models of dimension $d$ with coefficients in $\Phi$, fibre volumes in $\mathcal V_{\mathrm{fib}}$ and volume polynomial $\sigma(t)$. We denote by $M_{\mathrm{slc}}$ their coarse moduli space, see Definition \ref{defmodel}.
\begin{ques}\label{queeff}
For which tuples $(d,\Phi,\mathcal V_{\mathrm{fib}},\sigma)$ is $S_{\mathrm{slc}}(d,\Phi,\mathcal V_{\mathrm{fib}},\sigma)$ non-empty?
\end{ques}

To address this question, we consider surfaces for which the base of the Iitaka fibration is a curve. Let $F$ be a general fibre. We put $e=A\cdot F$, $b=(K_X+B)\cdot A$ and $\alpha=A^2$. Recall the uniform bound on the multiple which defines this fibration in Proposition \ref{propiit} (see \cite{Todorov08,BirkarZhang16,CHL23}). In the next theorem, the coefficients are integral, $b>0$ is fixed, $\mathcal V_{\mathrm{fib}}=\{e\}$, and $A$ is $\mathbb Q$-Cartier.
\begin{theo}\label{thrlower}
The number $\alpha$ is bounded below. If $X$ is normal and $K_X+B\equiv F$, then we have $\alpha\geq-3e^2$. If moreover $A$ is reduced, then we have $\alpha\geq-3e$.
\end{theo}
The two bounds are different because $A$ can be non-reduced. Indeed, the multiplicities of the components of $A$ enter $A^2$ with a square. Both bounds are sharp. In Example \ref{exsharp}, we have a section $R$ with $R\cdot F=1$ and $R^2=-3$. For $A=eR$, we get $A\cdot F=e$ and $A^2=-3e^2$. On the other hand, the rational families of Construction \ref{constr4} give us a reduced sum of $e$ disjoint $(-3)$-sections, which means that $A^2 = -3e$.

Now, we put $\sigma(t)=2e\beta t+\alpha t^2$. Then we have the following theorem
\begin{theo}\label{threxist}
For every $e\in\mathbb N$, the space $M_{\mathrm{slc}}(2,\mathbb Z^{\geq0},\{e\},\sigma)$ is non-empty for every $\beta\in\mathbb N$ and every integer $\alpha\geq-3e$.
\end{theo}

Thus the moduli space is non-empty for every $\alpha\geq-3e$. If we replace $A$ by $mA$, the stability does not change, because here $K_X+B+t(mA)=K_X+B+(mt)A$. We apply the theorem with fibre degree $e/m$. In this way, we get more examples with non-reduced polarization. This works whenever $m$ divides $e$, $\alpha\in m^2\mathbb Z$ and $\alpha\geq-3em$. In particular, $m=e$ gives the endpoint $-3e^2$. For $\beta=1$, these statements do not show that the space is non-empty for every integer strictly between $-3e^2$ and $-3e$.

In \cite{Bir2}, the author also asks about the Kodaira dimension and the birational geometry of the components \cite[Question 11.16 and the following paragraph]{Bir2}. We put $\Phi_c=c\mathbb Z^{\geq0}$ for $c>0$.
\begin{ques}\label{quebir}
Can we determine the Kodaira dimension of each component of $M_{\mathrm{slc}}(d,\Phi_c,\mathcal V_{\mathrm{fib}},\sigma)$? When is a component uniruled, rationally connected or rational?
\end{ques}

For $e\geq10$, Theorem \ref{threllcmp} gives irreducible components of $M_{\mathrm{slc}}(2,\mathbb Z_{\geq0},\{e\},\sigma)$ for the polynomials $\sigma(t)=2\beta et+(2ke-re-w)t^2$ listed there. They come from elliptic ruled surfaces. These components are unirational, and the elliptic base curve is allowed to vary. For example, $(\beta,r,w,k)=(1,3,0,0)$ gives $\sigma(t)=2et-3et^2$. We prove these statements in Subsection \ref{sec63}.

\paragraph{\normalfont\bfseries Structure of the paper.} Section \ref{sec2} recalls the preliminaries. Section \ref{sec3} gives the geometric constructions and proves Theorems \ref{thrconst} and \ref{thrfib} in Subsections \ref{sec32} and \ref{sec33}, respectively. Section \ref{sec4} constructs the families and coarse moduli spaces, with the proof of Theorem \ref{thrcoarse} in Subsection \ref{sec46}. Section \ref{sec5} proves the unirationality and rationality results of Theorems \ref{thrunirat} and \ref{thrrat} in Subsections \ref{sec52} and \ref{sec53}, respectively. Finally, Section \ref{sec6} proves the bounds and non-emptiness statements of Theorems \ref{thrlower} and \ref{threxist} in Subsection \ref{sec61}, and the statements on the elliptic components in Subsection \ref{sec63}.

\section*{Acknowledgements}
I thank my supervisors Frank Kutzschebauch and Jan Draisma for their support. This paper extends my master's thesis. I thank Andreas H\"oring, who supervised that thesis, for his support throughout the project. I thank Thomas Dedieu for suggesting the thesis topic and answering my questions about Epema's thesis.

I thank Caucher Birkar for encouraging me to study the effectiveness of moduli spaces of stable minimal models and for suggesting the use of adjunction to find lower bounds. I thank Stefano Filipazzi, Xiaowei Jiang, Junpeng Jiao and Jia Jia for useful discussions. Part of this work was carried out with support from the Yau Mathematical Sciences Center, Tsinghua University. I thank the members of the Birational Geometry Group for their support and helpful discussions.

Claude Opus was used for proofreading and suggested changes for the inaccuracies. I independently checked all the mathematics contents and take full responsibility for the paper. 

\section{Preliminaries}
\label{sec2}

Unless we say otherwise, all varieties are projective over $\mathbb C$. For divisors on a normal variety, $\sim$, $\sim_{\mathbb Q}$ and $\equiv$ mean linear, $\mathbb Q$-linear and numerical equivalence. For a birational morphism $f:X\longrightarrow Y$, we denote by $f^{-1}_*B$ the strict transform of $B$. In a pair $(X,\Delta)$, the variety $X$ is reduced and equidimensional, and $\Delta$ is an effective divisor such that no component of $X$ lies in $\Supp\Delta$. When $X$ is normal and we consider discrepancies, we assume that $K_X+\Delta$ is $\mathbb Q$-Cartier. Every log Calabi-Yau pair $(X,D)$ in this paper satisfies $K_X+D\sim0$.

\subsection{Ruled surfaces and elementary transformations}
\label{sec21}

We refer to \cite{Har1} for ruled surfaces. Let $\Gamma$ be a smooth projective curve of genus $q$, and let $\xi:X\longrightarrow\Gamma$ be a geometrically ruled surface. Since the fibres of $\xi$ are isomorphic to $\mathbb P^1$ and $\xi$ has a section, there is a locally free sheaf $\mathcal E$ of rank two on $\Gamma$ such that
\begin{center}
$ X \simeq \mathbb P_{\Gamma}(\mathcal E).$
\end{center}
We use the convention of Hartshorne. Hence we have $\mathbb P(\mathcal E)=\Proj_{\Gamma}\Sym\mathcal E$, and $\mathcal O_{\mathbb P(\mathcal E)}(1)$ is the tautological quotient.

Since twisting $\mathcal E$ by a line bundle does not change the ruled surface, we choose $\mathcal E$ normalized, so that $H^0(\Gamma,\mathcal E)\neq 0$ and
\begin{center}
$ H^0(\Gamma,\mathcal E\otimes L)=0 \ \text{for every line bundle } L \text{ with } \deg L<0.$
\end{center}
Let $C_0$ be the section given by $\mathcal O_{\mathbb P(\mathcal E)}(1)$. Let $f$ be the numerical class of a fibre, and denote by $D$ a divisor such that
\begin{center}
$ \det \mathcal E \simeq \mathcal O_{\Gamma}(-D) \ \text{with} \ e:=\deg D.$
\end{center}
Then
\begin{center}
$ \Pic(X)\simeq \mathbb Z[C_0]\oplus \xi^*\Pic(\Gamma)\ \text{and} \
N^1(X)\simeq \mathbb Z[C_0]\oplus \mathbb Z[f],$
\end{center}
and the intersection numbers are
\begin{center}
$ C_0^2=-e, \ C_0\cdot f=1 \ \text{and} \ f^2=0.$
\end{center}
The canonical divisor is
\begin{center}
$ K_X\sim -2C_0+\xi^*(K_\Gamma-D).$
\end{center}
If $\mathcal E\simeq \mathcal O_\Gamma\oplus \mathcal O_\Gamma(-D)$, there is a second section
\begin{center}
$C_\infty \sim C_0+\xi^*D.$
\end{center}
It is disjoint from $C_0$, and we have $C_\infty^2=e$. When $e>0$, we call $C_0$ the negative section and $C_\infty$ the positive section.

On the decomposable ruled surface $X_{\min}=\mathbb P_\Gamma(\mathcal O_\Gamma\oplus \mathcal O_\Gamma(-D))$ with $\deg D=e>0$, the sections satisfy $C_0^2=-e$, $C_\infty^2=e$, and $C_\infty\sim C_0+\xi^*D$. In the elliptic case, we denote this divisor by $A$. Epema denotes the same surface by $\mathbb P(\mathcal O_\Gamma\oplus\mathcal O_\Gamma(D))$ with $\deg D<0$ \cite{Epe}. We use the opposite sign for $D$.

For $a\geq 0$ and a divisor $\Delta$ on $\Gamma$, we have
\refstepcounter{equation}\label{eqpush}
\begin{center}
$ R^i\xi_*\mathcal O_X(aC_0+\xi^*\Delta)=0\ (i>0) \ \text{and} \
\xi_*\mathcal O_X(aC_0+\xi^*\Delta) \simeq \Sym^a\mathcal E\otimes \mathcal O_\Gamma(\Delta).$\enspace\textnormal{(\theequation)}
\end{center}
Thus,
\begin{center}
$ H^i\bigl(X,\mathcal O_X(aC_0+\xi^*\Delta)\bigr) \simeq H^i\bigl(\Gamma,\Sym^a\mathcal E\otimes \mathcal O_\Gamma(\Delta)\bigr).$
\end{center}
We use these formulas for the sections of the linear systems below and for their restrictions to the boundary.

We pass from one relatively minimal ruled model to another by elementary transformations. Let $f_P$ be the fibre through a point $P\in X$, and let $\mu:\widetilde X\longrightarrow X$ be the blow-up of $P$. The strict transform $\widetilde f_P$ is a $(-1)$-curve. By contracting it, we get another ruled surface
\begin{center}
$ X' \longleftarrow \widetilde X \longrightarrow X$
\end{center}
over $\Gamma$. This is the elementary transformation at $P$. If $P\in C_0$ lies over $Q\in\Gamma$, then the new surface is $\mathbb P(\mathcal E')$ with
\begin{center}
$\det\mathcal E'\simeq \mathcal O_\Gamma(-D-Q)$ and $ e'=e+1.$
\end{center}
The strict transform of $C_0$ is again a minimal section. Any two relatively minimal ruled surfaces over $\Gamma$ are related by finitely many such transformations. For this reason, we consider the tuples of numerical CHS invariants up to these transformations.

Maruyama describes the automorphism group of a decomposable ruled surface with $\deg D>0$ in \cite{Maruyama71}. We have an exact sequence
\begin{center}
$ 1\longrightarrow \Aut_{\Gamma}(X) \longrightarrow \Aut(X) \longrightarrow \Aut(\Gamma)_{[\mathbb P(\mathcal E)]},$
\end{center}
Here the last group consists of the automorphisms of $\Gamma$ which preserve the ruled surface. We also have
\begin{center}
$\Aut_{\Gamma}(\mathbb P(\mathcal E)) \simeq \Aut(\mathcal E)/\mathbb G_m.$
\end{center}
Now let $\mathcal E=\mathcal O_\Gamma\oplus\mathcal O_\Gamma(-D)$ with $\deg D>0$. Then we have $H^0(\Gamma,\mathcal O_\Gamma(-D))=0$. Hence the relative group is an extension of $\mathbb G_m$ by $H^0(\Gamma,\mathcal O_\Gamma(D))$. When we fix $C_0$ and $C_\infty$ separately, this vector part becomes trivial. In this situation, only the scalars on the fibres remain. We also keep the automorphisms of the base which preserve the boundary choices.
\subsection{Weighted clusters and virtual multiplicities on a smooth surface}
\label{sec22}

For clusters, we refer to \cite[Chapters 3--4]{Cas}. The author studies the Enriques diagrams in \cite[Section 3.9]{Cas}. For the parameter spaces of clusters, they are appeared in \cite{Roe04} and \cite{KPT11}.
\begin{defi}
\label{defprox} Let $X$ be a smooth surface, and let
\begin{center}
$ X_r\xrightarrow{\sigma_r}X_{r-1} \longrightarrow\cdots\longrightarrow X_1\xrightarrow{\sigma_1}X_0=X$
\end{center}
be a sequence of blow-ups at points $p_j\in X_{j-1}$. We call $p_j$ a \textit{proper point} if $X_{j-1}\to X$ is an isomorphism near $p_j$. Otherwise, we call $p_j$ an infinitely near point over its image on $X$. In the terminology of \cite{Roe04}, the sequence $(p_1,\ldots,p_r)$ is an \textit{ordered cluster}. Its \textit{blow-up} is the composite $X_r\to X$. We denote by $E_j$ both the exceptional divisor of $\sigma_j$ and its total transform on $X_r$. The points of the exceptional curve $E_j\subset X_j$ form the \textit{first neighbourhood} of $p_j$. We say that $p_j$ is \textit{proximate} to $p_i$ if $p_j$ lies on the strict transform of $E_i$ in $X_{j-1}$. We write $p_j\succ p_i$ for this relation. In this case, we have $i<j$. A point which is proximate to exactly one point is \textit{free}. A point which is proximate to two points is a \textit{satellite}. The \textit{proximity matrix} of the sequence is the $r\times r$ matrix $P=(P_{ij})$ with
\begin{center}
$ P_{ij} =
\begin{cases}
        1 & \text{if } i=j,\\
        -1 & \text{if } p_i\succ p_j,\\
        0 & \text{otherwise.}
        \end{cases}$
\end{center}
The $E_j$ form an orthogonal basis of the exceptional part of $\Pic X_r$, with $E_j^2=-1$. The strict transforms $\widetilde E_j$ on $X_r$ are given by
\begin{center}
$ (\widetilde E_1,\ldots,\widetilde E_r)=(E_1,\ldots,E_r)P.$
\end{center}
Thus we have $\widetilde E_j=E_j-\sum_{p_i\succ p_j}E_i$, see \cite[Section 1]{Roe04}.
\end{defi}

\begin{defi}
\label{defmult} Let $X$ be a smooth surface and let $C\subset X$ be an effective Cartier divisor. If $p\in X$ and $f\in\mathcal O_{X,p}$ is a local equation of $C$, then
\begin{center}
$\mult_p(C):=\max\bigl\{m\in\mathbb Z^{\geq 0}\mid f\in\mathfrak m_p^m \bigr\}.$
\end{center}

Let $p_\alpha$ be one of the centres of a sequence of blow-ups as in Definition \ref{defprox}, and let $C_\alpha$ be the strict transform of $C$ there. Then the \textit{effective multiplicity} of $C$ at the infinitely near point $p_\alpha$ is
\begin{center}
$ \mult_{p_\alpha}(C):=\mult_{p_\alpha}(C_\alpha).$
\end{center}

A \textit{weighted cluster} is a finite collection
\begin{center}
$ \mathcal K=\{(p_\alpha,m_\alpha)\}_{\alpha\in I}.$
\end{center}
Here the $p_\alpha$ are proper or infinitely near points, and we keep their proximity relations. The weights are $m_\alpha\in\mathbb Z_{>0}$. In \cite[Chapter 4]{Cas}, these weights are called \textit{virtual multiplicities}. To impose them on a curve $C$, we use the chosen order of the blow-ups, and we put
\begin{center}
$ C^{(0)}=C,\
C^{(j)}=\sigma_j^*C^{(j-1)}-m_jE_j.$
\end{center}
Here $E_j$ is the exceptional curve on $X_j$. We call $C^{(j)}$ the virtual transform. So at each step, we subtract the weight $m_j$. We say that $C$ contains $\mathcal K$ if $\mult_{p_j}(C^{(j-1)})\geq m_j$ for every $j$. In the terminology of \cite[Chapter 4]{Cas}, $C$ \textit{goes virtually through} $\mathcal K$. In this case, all these virtual transforms are effective. If the effective multiplicities equal the weights at every centre, then the virtual transforms are the strict transforms. We say that $C$ \textit{goes sharply through} $\mathcal K$ if these equalities hold and $C$ has no singular points outside $\mathcal K$, proper or infinitely near, see \cite[Chapter 4]{Cas}. For a special member with larger multiplicity, a virtual transform can also contain an exceptional curve. We use virtual transforms when we impose the linear conditions on $H^0(X,\mathcal O_X(C))$. A weighted cluster is \textit{consistent} if the proximity inequalities hold at every point. This means that we have
\begin{center}
$m_\alpha \geq \sum_{p_\beta\succ p_\alpha}m_\beta.$
\end{center}
The sum runs over all points $p_\beta$ on the strict transform of the exceptional curve of the blow-up at $p_\alpha$. In particular, the multiplicities do not increase along a proximity chain.
\end{defi}
A finite subscheme of a smooth surface is curvilinear at a point $p$ of its support if its ideal has the form $(v,u^k)$ in local coordinates $(u,v)$ for some $k\geq1$. Such a subscheme lies on the smooth curve $B=(v=0)$. Its local ring is $\mathbb C[u]/(u^k)$, and we get that its length is $k$. Let $k>1$. After the blow-up of $p$, we get the same condition with exponent $k-1$ on the strict transform of $B$. For the resolution, we use the cluster of the first $k$ points of $B$ at $p$. Here, each of these points has weight one and the first point is $p$. Each later point is the intersection of the strict transform of $B$ with the exceptional curve of the preceding point. The point after $p$ is proximate only to its previous point. Hence, it is free in the sense of Definition \ref{defprox}. There is a simple example which can give us a good insight: The curve $C=(v+u^k=0)$ cuts out the divisor $kp$ on $B$, and it is smooth at $p$. We observe that in this case the intersection multiplicity $k$ with $B$ gives $k$ centres of multiplicity one. We use the case $k=2$ in Subsection \ref{sec34}. We now consider a linear system without fixed components. Its \textit{cluster of base points} consists of all its proper and infinitely near base points. Here, we can put a weight on each centre by the multiplicity of a general member at this centre. Then, we blow up this cluster, and subtract the exceptional multiplicities. In the end of the process, we get a base-point-free linear system.

\begin{defi}
\label{defenriq} The \textit{Enriques diagram} $\mathfrak N(\mathcal K)$ of a weighted cluster $\mathcal K$ is a finite rooted forest with extra data. Its vertices are the points of $\mathcal K$. Each non-root is joined to the point from which it was blown up. The extra data are the proximity relation $\succ$ and the weights $m_\alpha$. The points themselves do not appear in the diagram. So two clusters have the same diagram if there is an isomorphism of the forests which preserves the proximity relations and the weights. This abstract tuple of forests, proximity relations and weights is an Enriques diagram in the sense of \cite[Definition 2.1]{KPT11}. A root is proximate to no vertex. A non-root is proximate to the vertex right before it, and to at most one more vertex. This further vertex is a remote predecessor. The weighted cluster is consistent. A vertex has at most two satellite immediate successors, and these are satellites of different vertices. The \textit{dimension} of $\mathfrak N$ is
\begin{center}
$\dim\mathfrak N=\sum_{V}t(V),$
\end{center}
with
\begin{center}
$ t(V)=
\begin{cases}
        2 & \text{if $V$ is a root},\\
        1 & \text{if $V$ is a free non-root},\\
        0 & \text{if $V$ is a satellite};
        \end{cases}$
\end{center}
see \cite[2.6]{KPT11}.
\end{defi}

By Definition \ref{defprox}, the exceptional configuration of the blow-up of $\mathcal K$ depends only on $\mathfrak N(\mathcal K)$. The same holds for all its intersection numbers. If the multiplicities of $C$ at the points of $\mathcal K$ equal the weights, then
\refstepcounter{equation}\label{eqstrict}
\begin{center}
$\widetilde C\cdot\widetilde E_j = m_j-\sum_{p_i\succ p_j}m_i .$\enspace\textnormal{(\theequation)}
\end{center}
Conversely, we can read the diagram from this configuration. So the diagram is a discrete invariant, and we fix it for the families below.

Let $X_{r-1}$ be the variety of ordered clusters of $r\geq1$ points on a smooth surface $S$. Let $\Cl(\mathfrak N)$ be the locus where the centres have the proximity relations and the weights of $\mathfrak N$. Here the positions of the centres vary. The relative version of this locus lies in the $r$-fold iterated blow-up of a smooth surface family $\mathcal S\to Y$ with geometrically irreducible fibres.

Let $E_j$ be the total exceptional divisors on the universal blown-up surface. In the terminology of \cite[Section 2]{Roe04}, the virtual exceptional divisor of the vertex $j$ of $\mathfrak N$ is
\begin{center}
$E_j-\sum_{p_i\succ p_j\text{ in }\mathfrak N}E_i.$
\end{center}
Here the sum uses the fixed diagram, even when a special sequence has more proximities. On $\Cl(\mathfrak N)$, these divisors are the strict exceptional curves. The effective locus below is the locus where all these divisors are effective.

\begin{prop}\label{propclust}
The locus $\Cl(\mathfrak N)$ is smooth and irreducible of dimension $\dim\mathfrak N$, and open in the closed locus of effective virtual exceptional divisors. Relatively, it is smooth over $Y$ with irreducible geometric fibres of dimension $\dim\mathfrak N$.
\end{prop}

\begin{proof}
The absolute assertions follow from \cite[Propositions 2.2, 2.4 and 2.6]{Roe04}, and the relative assertion from \cite[Theorem 3.10]{KPT11}.
\end{proof}

\begin{rema}
\label{remstrata} The loci $\Cl(\mathfrak N)$ do not form a stratification of $X_{r-1}$. One locus may meet the closure of another without being contained in it \cite[Section 3]{Roe04}.
\end{rema}

\subsection{Singularities and birational modifications}
\label{sec23} For singularities and birational modifications, we refer to \cite{KolSing,KM98}.

Let $(Y,\Delta_Y)$ be a normal pair. A log resolution is a birational morphism $f\colon X\longrightarrow Y$ with the following property. The variety $X$ is smooth, and the union of the exceptional locus of $f$ with the strict transform of $\Supp\Delta_Y$ is a simple normal crossing divisor. If $K_Y+\Delta_Y$ is $\mathbb Q$-Cartier, then
\begin{center}
$K_X+f^{-1}_*\Delta_Y \sim_{\mathbb Q} f^*(K_Y+\Delta_Y)+\sum_E a(E,Y,\Delta_Y)E,$
\end{center}
Here we sum over the prime exceptional divisors of $f$. The rational number $a(E,Y,\Delta_Y)$ is the discrepancy of $E$ with respect to $(Y,\Delta_Y)$. The pair is called terminal, canonical, klt, plt or lc when the following inequalities hold:
\begin{center}
$\begin{array}{c|c}
\text{type} & \text{condition} \\
\hline
\text{terminal} & a(E,Y,\Delta_Y)>0 \text{ for all exceptional } E,\\
\text{canonical} & a(E,Y,\Delta_Y)\geq 0 \text{ for all exceptional } E,\\
\text{klt} & a(E,Y,\Delta_Y)>-1 \text{ for all divisors } E,\\
\text{plt} & a(E,Y,\Delta_Y)>-1 \text{ for all exceptional } E,\\
\text{lc} & a(E,Y,\Delta_Y)\geq -1 \text{ for all divisors } E.
\end{array}$
\end{center}
Here ``all divisors'' means all prime divisors over $Y$. This also includes the prime divisors on $Y$ itself. A log canonical pair is divisorial log terminal, or dlt, if the following holds. After we remove a closed subset of codimension at least two, the pair is simple normal crossing, and the exceptional discrepancies are greater than $-1$.

For stable limits, we also need semi-log canonical singularities. A pair $(X,\Delta)$ is semi-log canonical, or slc, if the following conditions hold. The variety $X$ is $S_2$, and it has only nodal singularities in codimension one. No component of $\Delta$ lies in the conductor. The divisor $K_X+\Delta$ is $\mathbb Q$-Cartier. Let $\nu\colon X^\nu\longrightarrow X$ be the normalization with conductor divisor $B^\nu$. Then $(X^\nu,\Delta^\nu+B^\nu)$ is log canonical. These are the singularities which we allow for non-normal stable pairs.

The following birational modifications occur in our constructions.

\begin{defi}\label{defcanmod}
Let $X$ be a normal variety with $K_X$ $\mathbb Q$-Cartier.
\begin{enumerate}
    \item The canonical modification of $X$ is a birational morphism
    \begin{center}
    $f\colon X^{\can}\longrightarrow X$
    \end{center}
    such that $X^{\can}$ has canonical singularities and $K_{X^{\can}}$ is $f$-ample. It is unique when it exists.
    \item A terminal modification, or terminalization, is a birational morphism
    \begin{center}
    $f\colon X^{\ter}\longrightarrow X$
    \end{center}
    such that $X^{\ter}$ is $\mathbb Q$-factorial terminal and $K_{X^{\ter}}$ is $f$-nef. However, it need not be unique.
\end{enumerate}
\end{defi}

\begin{defi}\label{deflcmod}
Let $(X,\Delta)$ be a normal pair with $K_X+\Delta$ $\mathbb Q$-Cartier. The log canonical modification is a birational morphism
\begin{center}
$f\colon X^{\lc}\longrightarrow X$
\end{center}
with the following property. Let $E$ be the reduced exceptional divisor of $f$, and put $\Delta^{\lc}=f^{-1}_*\Delta$. Then $(X^{\lc},\Delta^{\lc}+E)$ is log canonical, and
\begin{center}
$K_{X^{\lc}}+\Delta^{\lc}+E$
\end{center}
is $f$-ample. The log canonical modification is unique.
\end{defi}

In dimension two, we get these modifications from resolutions by the relative minimal model program. The boundary blow-ups below resolve the given clusters. They are different from the canonical and log canonical modifications, since the divisors which define these modifications are relatively ample.

\subsection{Stable pairs and stable minimal models}
\label{sec24} For stable pairs and stable minimal models, we refer to \cite{Kol2,FI21,Bir2}.

\begin{defi}\label{defpair}
Let $I\subset [0,1]\cap \mathbb Q$. A pair $(X,\Delta)$ with coefficients of $\Delta$ in $I$ is a stable pair if
\begin{enumerate}
    \item $(X,\Delta)$ is semi-log canonical;
    \item $K_X+\Delta$ is $\mathbb Q$-Cartier; and
    \item $K_X+\Delta$ is ample.
\end{enumerate}
\end{defi}

For the log Calabi-Yau surfaces below, we have $K_X+D\sim0$. If $C$ is ample, then
\begin{center}
$K_X+D+\epsilon C\sim_{\mathbb Q}\epsilon C$
\end{center}
is ample for $0<\epsilon\ll1$. This leads to the following definition of stable minimal models.

\begin{defi}\label{defmodel}
A stable minimal model is a triple $(X,B),A$ consisting of an slc pair $(X,B)$ and an effective $\mathbb Q$-divisor $A$ such that
\begin{enumerate}
    \item $K_X+B$ is semiample and defines a contraction $f\colon X\longrightarrow Z$;
    \item $K_X+B+tA$ is ample for all sufficiently small rational numbers $t>0$; and
    \item $(X,B+tA)$ is slc for all sufficiently small rational numbers $t>0$.
\end{enumerate}
If $\dim X=d$ and the coefficients of $B$ and $A$ belong to a set $\Phi\subset \mathbb Q^{\geq 0}$, we call it a $(d,\Phi)$-stable minimal model.
\end{defi}

We denote by $\mathcal V_{\mathrm{fib}}\subset\mathbb Q^{\geq0}$ the set of permitted fibre volumes, and we take $\sigma(t)\in\mathbb Q[t]$. A $(d,\Phi,\mathcal V_{\mathrm{fib}},\sigma)$-stable minimal model is a $(d,\Phi)$-stable minimal model $(X,B),A$ with the following properties. For a general fibre $F$ over each irreducible component of $Z$, we have
\begin{center}
$\vol(A|_F)\in \mathcal V_{\mathrm{fib}},$
\end{center}
and
\begin{center}
$\vol(K_X+B+tA)=\sigma(t)$
\end{center}
for every sufficiently small positive rational $t$.

For a family over a reduced base $S$, we take a projective morphism
\begin{center}
$\mathcal X\longrightarrow S$
\end{center}
with $\mathbb Q$-divisors $\mathcal B,\mathcal A$. We require that $(\mathcal X,\mathcal B+t\mathcal A)\to S$ is locally stable for all small $t>0$, and each geometric fibre $(\mathcal X_s,\mathcal B_s),\mathcal A_s$ is a stable minimal model. Locally stable means that the fibres are slc with $K_{\mathcal X/S}+\mathcal B+t\mathcal A$ is $\mathbb Q$-Cartier and these conditions must compatible with base change. If every $(\mathcal X_s,\mathcal B_s)$ is lc, we call the family an lc family of stable minimal models.

For $c\in\mathbb Q^{\geq0}$, we recall Birkar's notation in \cite{Bir2}
\begin{center}
$\Phi_c=c\mathbb Z^{\geq 0}$
\end{center}
for the coefficient set of intersection polynomial.

A family of $(d,\Phi_c,\mathcal V_{\mathrm{fib}},\sigma)$-stable minimal models has $(\mathcal X,\mathcal B),\mathcal A\to S$ with
\begin{center}
$\mathcal B=c\mathcal D, \ \text{and} \
\mathcal A=c\mathcal N.$
\end{center}
Here $\mathcal D,\mathcal N$ are relative Mumford divisors and every geometric fibre has the given numerical invariants.

For $\mathbb Q$-pairs, we follow \cite[Definition 2.7]{FI21}. A $\mathbb Q$-pair $(X,B;A)$ is an slc pair $(X,B+A)$ with the following extra data. The divisor $A$ is an effective $\mathbb Q$-Cartier divisor. The coefficients are fixed, and the polynomial is $p(t)=(K_X+B+tA)^{\dim X}$. There is a birational contraction $\pi$ to a stable pair with
\begin{center}
$K_X+B+A=\pi^*(K_{X^s}+B^s+A^s), \ \pi_*B=B^s,\ \pi_*A=A^s.$
\end{center}
It is a $\mathbb Q$-stable pair if $(X,B+(1-\epsilon)A)$ is stable for every sufficiently small $\epsilon>0$. Thus the log canonical divisor of a $\mathbb Q$-pair is semiample. For $\mathbb Q$-stability, we also need ampleness after we decrease the coefficients of $A$. When $B=0$, we write $(X;A)$ for the pair.

\subsection{Log Calabi-Yau surfaces and anticanonical pairs}
\label{sec25}

\begin{defi}\label{deflcy}
A log Calabi-Yau surface pair is a log canonical pair $(Y,D)$ with $\dim Y=2$ and
\begin{center}
$K_Y+D\sim 0.$
\end{center}
In other words, we have $K_Y\sim -D$. A stable log Calabi-Yau surface is a triple
\begin{center}
$((Y,D),C)$
\end{center}
with the following properties. The pair $(Y,D)$ is a projective log Calabi-Yau surface pair. The divisor $C$ is an effective ample Cartier divisor. The pair $(Y,D+\epsilon C)$ is log canonical for every sufficiently small rational $\epsilon>0$. On a resolution we keep the nef and big pullback of $C$.
\end{defi}

For families of stable log Calabi-Yau surfaces, we impose the conditions of Definition \ref{defmodel}, that is $K_{X_s}+D_s\sim 0$ on every geometric fibre.

For a stable log Calabi-Yau surface, we have
\begin{center}
$K_S+D+tC\sim_{\mathbb Q}tC,$
\end{center}
and thus,
\begin{center}
$\vol(K_S+D+tC)=t^2 C^2$
\end{center}
for $0<t\ll1$. So, we have that the volume polynomial of these stable minimal models is $\sigma(t)=C^2t^2$.

If $X$ is smooth and $D\in|-K_X|$ is reduced, we call $(X,D)$ an anticanonical pair. For rational $X$, the boundary $D$ is often a cycle of smooth rational curves. On the elliptic ruled surfaces here, the boundary $D=C_0+C_\infty$ consists of two disjoint elliptic curves. For a labelled pair, we order the components of $D$ as in \cite[Section 1]{Fried}.

\begin{defi}\label{deftaut}
A labelled anticanonical pair $(Y,D=D_1+\cdots+D_r)$ is called taut if the following holds. Every labelled anticanonical pair $(Y',D'_1+\cdots+D'_r)$ with the same self-intersection sequence
\begin{center}
$(D_1^2,\ldots,D_r^2) =(D_1'^2,\ldots,D_r'^2)$
\end{center}
is isomorphic to $(Y,D)$.
\end{defi}

Friedman distinguishes corner blow-ups and interior blow-ups \cite[Section 2]{Fried}. The difference is whether the centre is a node or a smooth point of the boundary. A corner blow-up has its centre at a node $D_i\cap D_{i+1}$. In this case, we take the reduced inverse image as the new boundary. An interior blow-up has its centre at a smooth point of $D_i$. In this case, we take the strict transform as the new boundary. The inverse operations are the corner blow-downs and the interior blow-downs. Gross, Hacking and Keel call a sequence of corner blow-ups a toric blow-up \cite[Section 1]{GHK1}. Corner blow-ups and corner blow-downs preserve tautness. An interior blow-down of a taut pair also preserves tautness \cite[Lemma 2.15]{Fried}. Hence we can use these operations to replace rational ruled surfaces and their boundaries by taut anticanonical pairs. This does not change the numerical classification. With these conventions for the boundary, a blow-up $\sigma:\widetilde Y\to Y$ at a smooth boundary point or at a node satisfies $K_{\widetilde Y}+\widetilde D=\sigma^*(K_Y+D)$. Conversely, let $K_X+W\sim0$, and let $g:X\to Y$ contract a $(-1)$-curve. In this case, we get $K_X+W=g^*(K_Y+g_*W)$ with compatible canonical divisors. Thus these blow-ups and blow-downs are crepant for the boundaries above.

In the terminology of \cite[Definition 2.1]{Fried}, an interior exceptional curve is a smooth rational $(-1)$-curve which is not contained in $D$. Following \cite[Definition 1.4]{GHK1}, we call a smooth rational $(-2)$-curve disjoint from $D$ an internal $(-2)$-curve. Take a chain with an interior exceptional curve at one end and internal $(-2)$-curves as its other components. Such a chain is called a generalized exceptional curve in \cite[Definition 2.3]{Fried}. We get such a chain in the following way. First, we blow up a smooth boundary point. Then, at each stage, we blow up the intersection of the boundary with the exceptional curve of the blow-up before.

Let $(Y,D)$ be smooth with a boundary cycle $D=D_1+\cdots+D_r$. Then the labelled pair has unobstructed deformations, and we have
\begin{center}
$\dim \Def(Y;D_1,\ldots,D_r) =h^1(Y,T_Y(-\log D)) =\rank\Lambda(Y,D),$
\end{center}
where
\begin{center}
$\Lambda(Y,D):=\bigl\{\alpha\in H^2(Y,\mathbb Z)\mid \alpha\cdot D_i=0\text{ for every }i\bigr\}$
\end{center}
is the lattice orthogonal to $D$. The pair is infinitesimally rigid exactly when $[D_1],\ldots,[D_r]$ span $H^2(Y,\mathbb C)$.

\subsection{Surfaces with canonical hyperplane sections}
\label{sec26}

Let $(S,C)$ be a surface with CHS. By Epema's classification \cite{Epe}, we have $\mathcal O_C(C)\simeq\omega_C$ on a general smooth $C\in|C|$. If $C$ is non-hyperelliptic, then its canonical system gives an embedding, and its genus is at least three. In the hyperelliptic case, the canonical map has degree two. Let $X\to S$ be the minimal resolution, and let $C_X$ be the strict transform. By adjunction, we get $C_X^2=2g(C)-2$ and $K_X\cdot C_X=0$. By Theorem \ref{threpema}, we have $K_S\sim0$, and the singularities are Gorenstein. In this paper, we study the log canonical rational and elliptic ruled cases. We call a curve $E$ $C$-trivial when $C\cdot E=0$.

For the numerical description on a relatively minimal ruled surface, see \cite[Proposition III.1.4 and Remark III.2.2.2]{Epe}. In the ruled cases, let $\xi\colon X_{\min}=\mathbb P_{\Gamma}(\mathcal O_\Gamma\oplus\mathcal O_\Gamma(-A))\longrightarrow \Gamma$ be relatively minimal. We put $q=g(\Gamma)$ and $e=\deg A$. Let $C_\infty$ be the positive section, and let $W\in|-K_{X_{\min}}|$. Let $L_1\subset|aC_\infty|$ be a system with base locus in $W$, and take a general $C_1\in L_1$. Let $m_1,\ldots,m_s$ be the multiplicities of the proper and infinitely near base points. The resolved strict transform $C$ satisfies $C^2=a^2e-\sum_{i=1}^s m_i^2$. It is disjoint from the transform of $W$ if and only if
\refstepcounter{equation}\label{eqmultsum}
\begin{center}
$\sum_{i=1}^s m_i=a(e-2q+2).$\enspace\textnormal{(\theequation)}
\end{center}
In this case, we have $K_X\cdot C=0$, and adjunction gives
\refstepcounter{equation}\label{eqgenus}
\begin{center}
$g(C)=1+\frac{1}{2}\left(a^2e-\sum_{i=1}^s m_i^2\right).$\enspace\textnormal{(\theequation)}
\end{center}
Moreover, the multiplicities of the base points satisfy
\refstepcounter{equation}\label{eqmultbnd}
\begin{center}
$1\leq m_i\leq a-1.$\enspace\textnormal{(\theequation)}
\end{center}
This means that no exceptional curve is a fixed component of the transformed system. We call $(a,e,q;\mathcal K)$ a \textit{tuple of numerical CHS invariants} if it satisfies \eqref{eqmultsum}, \eqref{eqgenus} and \eqref{eqmultbnd}. We also put a condition that $a^2e-\sum_i m_i^2\geq2$, and that the effectivity assumptions for the ruled surface hold. 

Now, let us take such a tuple as described above. Here, we might choose $W$ and $C_1\in|aC_\infty|$ such that $C_1\cap W$ is the given cluster. First, we blow up the cluster. After that, we contract $W$. This gives us the surface with CHS.

In the next section, we note that the cluster data consists of a weighted cluster and the labels of the boundary components. Here, the weighted cluster can have infinitely near centres. For $q=1$, we will also let the elliptic base curve vary.
\begin{rema}\label{remdata}
We recall the meaning of these numbers. The integer $a$ is the horizontal degree of the polarization on the ruled surface. The number $e$ is the self-intersection of the positive section. The integer $q$ is the genus of the base curve. The $m_i$ are the multiplicities of the base points in the cluster. These points lie on the anticanonical divisor. After the blow-ups, we have $C\cdot W=0$. So we can contract the anticanonical divisor. Then the curve $C$ becomes a canonical hyperplane section.
\end{rema}

\section{General constructions of stable log Calabi-Yau surfaces}\label{sec3}
To construct the surfaces in the classification of Epema, we keep the weighted cluster and the labels of the boundary components. We ask that a general member has the given weights as its multiplicities. For surfaces with CHS, we also ask that the whole intersection $C_1\cap W_1$ is resolved. These conditions give the marked clusters of Definition \ref{defmark} and the $\mathcal K$-systems of Definition \ref{defksys}. Construction \ref{constr1} starts from taut rational anticanonical pairs and from elliptic ruled pairs. By adding pencils to these pairs, Construction \ref{constr2} gives surfaces with CHS and $(-2)$-curves. After we describe the families and the moduli of Construction \ref{constr3} in Section \ref{sec4}, we study the stable minimal models of Constructions \ref{constr3} and \ref{constr4} in Section \ref{sec6}.

Let $(X_{\min},W_1)$ be a anticanonical pair with smooth divisor or a taut pair and $L$ be a line bundle on $X_{\min}$. In the specific constructions below, we take $L\sim aC_\infty$ with $X_{\min}$ a minimal ruled surface just for simplicity. For a more general rational ruled surface case, we must take $L\sim aC_\infty+\lambda f$, where $\lambda$ is the coefficient of the fibre class $f$.

\begin{defi}
\label{defmark} A \textit{marked anticanonical cluster}, or a \textit{marked cluster} for short, on $(X_{\min},W_1)$ is a finite collection
\begin{center}
$\mathcal K=\{(p_\alpha,m_\alpha,\prox_\alpha,\lambda_\alpha)\}_{\alpha\in I}$
\end{center}
with the following properties:

\begin{enumerate}
    \item Each $p_\alpha$ is a proper point of $W_1$, or an infinitely near point coning from blowing up earlier points of the cluster.

    \item The point $p_\alpha$ has multiplicity $m_\alpha>0$.

    \item The datum $\prox_\alpha$ is the set of earlier centres to which $p_\alpha$ is proximate. Here the weighted cluster $\mathcal K$ is consistent in the sense of Definition \ref{defmult}.

    \item The label $\lambda_\alpha$ is the boundary component that contains $p_\alpha$. At an infinitely near point, the label is a tangent direction. It is given by the strict transform of the component with a label of $W_1$.
\end{enumerate}
Let $C_1\in|L|$ be a curve. We say that $C_1$ contains $\mathcal K$ if its virtual transforms satisfy the given multiplicity conditions at each stage. This is the notion of Definition \ref{defmult}. When the effective multiplicities are equal to the weights, the virtual transforms are the strict transforms. In this case, the conditions on the sections are linear. We add the labels of the boundary components to the weighted Enriques diagram and call the result the \textit{labelled Enriques diagram}. Here, We let the positions of the proper centres vary.
\end{defi}
\begin{defi}\label{deftuple}
An \textit{admissible tuple} is a collection $(X_{\min},W_1,C_1;\mathcal K,\Sigma,R)$ of the following form. The pair $(X_{\min},W_1)$ is an anticanonical pair, and $C_1$ is a curve. The datum $\mathcal K$ is a marked anticanonical cluster. The equality $C_1|_{W_1}=\Sigma+R$ is a decomposition of the scheme-theoretic restriction to the boundary. Here $\Sigma$ is the part which $\mathcal K$ represents, and $R$ is an effective residual divisor on $W_1$. The numbers $a\geq 2$ and $e>0$ are integers, and $q\in\{-\infty,0,1\}$. The pair $(X_{\min},W_1)$ is one of the following:
\begin{enumerate}
    \item[(i)] a taut rational anticanonical pair, so that $W_1$ is a cycle of rational curves. Here $q=0$ when the chosen surface is a Hirzebruch surface $X_{\min}\simeq\mathbb F_e$, and $q=-\infty$ when it is the plane $X_{\min}\simeq\mathbb P^2$;
    \item[(ii)] a rational pair with irreducible boundary. In the case $q=0$, we have $X_{\min}\simeq\mathbb F_e$ with $e\in\{1,2\}$, and $W_1\in|-K_{\mathbb F_e}|$ is an integral nodal or smooth elliptic member. In the case $q=-\infty$, we have $X_{\min}\simeq\mathbb P^2$, and $W_1$ is an integral nodal or smooth plane cubic;
    \item[(iii)] a decomposable elliptic ruled pair
    \begin{center}
    $X_{\min}=\mathbb P_\Gamma\left(\mathcal O_\Gamma\oplus\mathcal O_\Gamma(-A)\right) \ \text{with} \ \deg A=e,$
    \end{center}
    over a smooth elliptic curve $\Gamma$, in the case $q=1$;
\end{enumerate}
In each case, we have $W_1\in|-K_{X_{\min}}|$ and $C_1\in|aC_\infty|$. In the plane case $q=-\infty$, the positive section $C_\infty$ means a line $H\subset\mathbb P^2$.

Let $\mu_{\mathcal K}:X_{\mathcal K}\longrightarrow X_{\min}$ be the blow-up of $\mathcal K$. Let $E_\alpha^*$ be the total transform on $X_{\mathcal K}$ of the exceptional divisor over $p_\alpha$. For the resolved polarizing divisor
\begin{center}
$C=\mu_{\mathcal K}^*C_1-\sum_\alpha m_\alpha E_\alpha^*,$
\end{center}
we require that $|C|$ is base-point-free and that $C^2>0$.

We put $q=0$ for a Hirzebruch surface and $q=1$ for an elliptic ruled surface. For the plane, which has no ruling, we put $q=-\infty$. We treat the plane and the Hirzebruch surfaces together. For this, we use elementary transformations and the contraction of a section.
\end{defi}

For a line bundle $L$ on $X_{\min}$, we use the blow-up of $\mathcal K$, and we put
\begin{center}
$L_{\mathcal K}:=\mu_{\mathcal K}^*L-\sum_\alpha m_\alpha E_\alpha^*.$
\end{center}
Then the linear system
\begin{center}
$\mathcal L^{\geq \mathcal K} := \mathbb P H^0(X_{\mathcal K},L_{\mathcal K}) \subset |L|$
\end{center}
is the subsystem of the curves which contain $\mathcal K$ with multiplicities at least $m_\alpha$. When the infinitely near points and the variation of the base curves play no role, we write $\mathcal K=(m_\alpha)_{\alpha\in I}$ for the cluster. When we know the number of singular points, we write $\mathcal K=(m_1,\ldots,m_s)$.

We put $F_{\mathcal K}=\sum_\alpha m_\alpha E_\alpha^*$ and $\mathcal I_{\mathcal K}=(\mu_{\mathcal K})_*\mathcal O_{X_{\mathcal K}}(-F_{\mathcal K})$. By the projection formula, we have
\begin{center}
$(\mu_{\mathcal K})_*\mathcal O_{X_{\mathcal K}}(L_{\mathcal K}) =\mathcal O_{X_{\min}}(L)\otimes\mathcal I_{\mathcal K}.$
\end{center}
Next, we will compare two blow-ups of clusters. Here, we assume that their polarizations should have the same pullback on a common smooth surface.
\begin{lemm}\label{lemlinear}
The complete polarizing system is determined by weighted diagram with the marked points and the line bundle $L$.
\end{lemm}
\begin{proof}
By taking sections of the equality of pushforward and tensor by $\mathcal I_{\mathcal K}$ above, we can identify $H^0(X_{\mathcal K},L_{\mathcal K})$ with $H^0(X_{\min},\mathcal O(L)\otimes\mathcal I_{\mathcal K})$. Now, let $h:V\to X_{\mathcal K}$ be a proper birational morphism between smooth surfaces. We have $h_*\mathcal O_V=\mathcal O_{X_{\mathcal K}}$. Thus the projection formula gives an isomorphism
\begin{center}
$h^*:H^0(X_{\mathcal K},L_{\mathcal K}) \xrightarrow{\sim}H^0(V,h^*L_{\mathcal K}).$
\end{center}
Let $h_i:V\to X_{\mathcal K_i}$ be two blow-ups of clusters with $h_1^*L_{\mathcal K_1}=h_2^*L_{\mathcal K_2}$. it is not difficult to see that $(h_2^*)^{-1}\circ h_1^*$ identifies the two section spaces and this finishes the lemma.
\end{proof}
The abstract diagram fixes the form of the conditions. Its marked points determine the ideal. Now we recover the cluster of base points. At each stage, we blow up the base points of the given system. Then we subtract the exceptional divisor with the multiplicity of a general member. We stop when the transformed system is base-point-free. After this, a general member has multiplicity zero at every new centre. So further blow-ups will not make the section space changes, as in \cite[Chapter III, Section 1]{Epe}.

\begin{defi}\label{defksys}
A $\mathcal K$-\textit{system} is a non-empty linear subsystem
\begin{center}
$\mathcal L_1\subset \mathcal L^{\geq \mathcal K}$
\end{center}
with the following properties:

\begin{enumerate}
    \item Every member of $\mathcal L_1$ contains the cluster $\mathcal K$.

    \item General members of $\mathcal L_1$ go sharply through $\mathcal K$.

    \item For a general member $C_1\in \mathcal L_1$, the strict transforms $C$ and $W$ of $C_1$ and $W_1$ on $X_{\mathcal K}$ satisfy $C\cap W= R^{\prime} \in |R|$.
\end{enumerate}
\end{defi}

Here, we need the sharpness condition only for a general member of $\mathcal L_1$, since special members may have higher multiplicity or even more singularities. Thus, the members of $\mathcal L_1$ going sharply through $\mathcal K$ will form a dense open subset $\mathcal L_1^{\mathrm{sh}}\subset \mathcal L_1$.

\begin{lemm}
The polarizing system on a resolution of a surface with CHS has no fixed divisorial part.
\end{lemm}
\begin{proof}
The hyperplane system in Definition \ref{defchs}, or the finite system in case (b), is base-point-free. The same holds for its pullback. A general member avoids the images of the exceptional curves. Thus its pullback is its strict transform.
\end{proof}

Notice that this statement is only about the CHS system. For a general stable log Calabi-Yau pair, we keep the polarizing system and its virtual transforms as in Lemma \ref{lemlinear}. Indeed, if we removed a curve of positive polarizing degree, the polarization would change.

\subsection{General construction for stable log Calabi-Yau surfaces and surfaces with CHS}\label{sec31}
Epema gave a relationship between the singularities of surfaces with CHS and the relatively minimal model of their resolutions \cite[Chapter II]{Epe}. Here, we will reconstruct those surfaces from the marked clusters of Definition \ref{defmark} together with the tuples of form $(a,e,q;m_1,\ldots,m_s)$. Let $\mu:X\longrightarrow X_{\min}$ be the resolution of the cluster and $C_1$ be a divisor on $X_{\min}$. Then, the strict transform of a general member $C_1$ with multiplicity $m_\alpha$ at $p_\alpha$ is
\begin{center}
$C=\mu^*C_1-\sum_\alpha m_\alpha E_\alpha^*.$
\end{center}
Here, the divisors $E_\alpha^*$ are the total exceptional transforms.

For the reverse construction, let $(S,C)$ be a log canonical non-$K3$ surface with CHS. Let $\pi:X\longrightarrow S$ be its minimal resolution. We put $C_X=\pi^*C$. We denote by $W$ the reduced exceptional divisor over the non-canonical singularities. Since $S$ is Gorenstein and log canonical, we have $K_X+W=\pi^*K_S\sim0$. Moreover, a general $C_X$ is the strict transform of $C$. So it is disjoint from $W$. We take the hyperplane system $V\subset|C_X|$, and we follow the given reverse MMP. This is a sequence of blow-downs
\begin{center}
$\phi:X=X_s\longrightarrow X_{s-1}\longrightarrow\cdots \longrightarrow X_1=X_{\min}.$
\end{center}
Let $W_i$ and $V_i$ be the push-forwards to $X_i$. Each centre $p_\alpha$ lies on $W_i$. So we can put $m_\alpha:=\mult_{p_\alpha}(V_i)$. In this way, we keep all the data of the marked cluster $\mathcal K=\{(p_\alpha,m_\alpha,\prox_\alpha,\lambda_\alpha)\}$. We put $\mathcal L_1:=\phi_*V\subset|L|$. A general member $C_{\min}\in\mathcal L_1$ has multiplicity $m_\alpha$ at each centre. After the resolution, the strict transforms of $C_{\min}$ and $W_{\min}$ are disjoint. This implies that $\mathcal L_1$ is a $\mathcal K$-system. Conversely, let $(X_{\min},W_1)$ be an initial pair with $W_1\in|-K_{X_{\min}}|$. Let $\mathcal L_1\subset|L|$ be a $\mathcal K$-system for a cluster $\mathcal K$ marked on $W_1$. For a general $C_1\in\mathcal L_1$, we denote by $C$ and $W$ the strict transforms of $C_1$ and $W_1$ on $\mu:X\longrightarrow X_{\min}$. In this situation, the divisor $C$ is nef and big. Since all centres lie on $W_1$, we have $K_X+W\sim0$. By condition (3) of Definition \ref{defksys}, we have $C\cap W=\varnothing$. For sufficiently divisible $m>0$, the system $|mC|$ defines a birational contraction $\rho:X\longrightarrow Y$. This contraction contracts exactly the $C$-trivial curves. Among them are the components of $W$ of $C$-degree zero. In this way, we get the stable log Calabi-Yau surface $(Y,D),C_Y$, where $D=\rho_*W$ is the remaining divisorial part. A non-hyperelliptic general member of $|C|$ maps to a canonical hyperplane section. In the hyperelliptic case, the polarizing system gives a morphism of degree two.

\begin{construction}\label{constr1}

We have the following diagram.
\begin{center}
\begin{tikzcd}
X \arrow{r}{\pi} \arrow[swap]{d}{\phi} & Y \arrow[dashed]{r}{\Psi} & S \ni x \arrow{r}{i} & \mathbb{P}^g \\X_{\min} \arrow{r}{\xi}& \Gamma
\end{tikzcd}
\end{center}
The notation is as follows.
\begin{enumerate}
    \item $\pi:X\to Y$ resolves the stable pair which we get by contracting the curves of log canonical degree zero. We fix a sufficiently small rational $\tau>0$, and we put $C_Y=\pi_*C$. Then we have $K_X+W+\tau C=\pi^*(K_Y+D+\tau C_Y)$. Thus $(X;W+\tau C)$ is a $\mathbb Q$-pair.
    \item $(Y,D),C$ is the stable log Calabi-Yau minimal model. Here $\pi$ is the minimal resolution of the singularities of $Y$, or $\pi$ is defined by $|mC|$.
    \item $\Psi$ is the birational map from $Y$ to $S$. Here $S$ is a surface with CHS embedded in $\mathbb{P}^g$. The map $\Psi$ is defined by a big divisor.
    \item $X_{\min}$ is the initial surface of the admissible tuple, birational to $X$. $\dim \Gamma$ is $0$ or $1$.
    \item $x \in S$ is a cusp singularity when $q \in\{-\infty, 0\}$ and $W$ is a nodal curve or a cycle of rational curves. It is a simple elliptic singularity when $W$ is smooth elliptic. For $q=1$, there are two simple elliptic points, see Theorem \ref{threllchs}.
    \item We also denote the maps in the sequence of blow-ups by $\phi_t$:
    \begin{center}
    $X = X_s \xrightarrow{\ \phi_{s-1}\ } X_{s-1} \xrightarrow{\ \phi_{s-2}\ } \cdots \xrightarrow{\ \phi_{1}\ } X_1 = X_{\min},$
    \end{center}
    Here $\phi_t$ is the blow-up of the base points of a linear system $\mathcal L_t$ on $X_t$. We define the systems $\mathcal L_t$ step by step. We start with $\mathcal L_1\subset |L|$ on the minimal model, whose base scheme lies on $W_1\in |-K_{X_{\min}}|$. Then we take $\mathcal L_{t+1}$ to be the strict transform of $\mathcal L_t$ under $\phi_t$. Thus $\phi$ resolves the indeterminacy of the rational map defined by $\mathcal L_1$. The strict transform $\mathcal L_s$ on $X$ is base-point-free. A general member of $\mathcal L_s$ is disjoint from $W\in |-K_X|$. The morphism defined by a sufficiently divisible multiple of $\mathcal L_s$ contracts $W$. In the CHS case $R=0$, the composite $\Psi\circ\pi$ maps $X$ onto the surface with canonical hyperplane sections $S\subset\PP^g$.

\end{enumerate}
\end{construction}

We consider the normal log canonical CHS case. We apply \cite[Lemma 2.16]{FI21} to $(S,\tau C_S)$. This gives the canonical modification $f:Z\to S$ with an effective boundary $W_Z$. We have $K_Z+W_Z=f^*K_S$, and $(Z,W_Z+\tau C_Z)$ is log canonical, where $C_Z=f^*C_S$. Since $K_S\sim0$, we have
\begin{center}
$K_Z+(1-\epsilon)(W_Z+\tau C_Z) \sim_{\mathbb Q}\epsilon K_Z+(1-\epsilon)\tau C_Z.$
\end{center}
Here $K_Z$ is $f$-ample and $C_S$ is ample. So this divisor is ample for $0<\epsilon\ll1$. Hence $(Z;W_Z+\tau C_Z)$ is a $\mathbb Q$-stable pair. The pair on the smooth resolution need not be $\mathbb Q$-stable. Indeed, a $C$-trivial internal $(-2)$-curve still has degree zero after this perturbation.

\subsection{Reconstruction from rational and elliptic ruled surfaces}
\label{sec32} We fix $a\geq2$ and a taut rational anticanonical pair $(X_{\min},W_1)$ as in Construction \ref{constr1}. Here for simplicity, we assume $X_{\min}\simeq\mathbb F_e$ with $e\geq1$. In the next proposition, we will also assume that the linear system after blow-ups is base-point-free and that its divisor $C$ satisfies $C^2>0$.

\begin{prop}\label{propdec}
Each decomposition gives a stable log Calabi-Yau surface. On this surface, the polarization meets the boundary along $R$. For $R=0$, we get a surface with CHS. It has only one minimally elliptic singularity.
\end{prop}
\begin{proof}
We divide the proof into several steps.\\
\\
\textbf{Step 1.} First, we reduce the case where the polarization is semi-ample but not base-point-free to the base-point-free case. Let $L_{\mathcal K}=\mathcal O_{X_{\mathcal K}}(C)$ be the resolved polarization. Assume that $L_{\mathcal K}$ is semi-ample, but $|L_{\mathcal K}|$ is not base-point-free. Then $|mL_{\mathcal K}|$ is base-point-free for some $m>0$. The section rings of $L_{\mathcal K}$ and of $mL_{\mathcal K}$ have the same $\Proj$. So they define the same contraction of $X_{\mathcal K}$. It contracts the same curves, since a curve has $L_{\mathcal K}$-degree zero if and only if it has $mL_{\mathcal K}$-degree zero. Hence the stable log Calabi-Yau surface $((Y,D),C)$ of Steps 2 and 3 does not change. We use $|mL_{\mathcal K}|$ only to define the contraction, and we keep $C$ as the polarization. So we can assume that $|L_{\mathcal K}|$ is base-point-free without losing any generality.\\
\\
\textbf{Step 2.} Definition \ref{deftuple} gives directly the smooth log Calabi-Yau pair $(X,W)$ on the blow-up $\phi:X=X_{\mathcal K}\longrightarrow X_{\min}$ of $\mathcal K$, with $K_X+W\sim0$ and $C\cap W=R$, and the contraction of $X$ gives $Y$.\\
\\
\textbf{Step 3.} We contract the curves of $C$-degree zero. The divisor $C$ is nef and big. Hence for $m\gg0$, the linear system $|mC|$ defines a birational morphism
\begin{center}
$\phi_{|mC|}:X\longrightarrow Y\subset \PP H^0(X,mC).$
\end{center}
This morphism contracts the curves of $C$-degree zero and only these curves. Among them are the components of $W$ of $C$-degree zero. The divisor $C$ descends to an ample divisor $C_Y$ on $Y$. We denote by $D$ the image of $W$, and we get $K_Y+D\sim0$. When $R=0$, adjunction gives $\mathcal O_C(C)\simeq\omega_C$ and $C^2=2g(C)-2>0$. In this case, the original system $|C|$ gives the canonical embedding of a general non-hyperelliptic member. In the hyperelliptic case, it gives the canonical map of degree two.\\
\\
\textbf{Step 4.} We use the complete section ring of $C$ to define the contracted surface $Y$. We descend the original system to $|C_Y|$. If $R=0$, we get the CHS system of Definition \ref{defchs}. For $R\neq0$, the components of $W$ of positive $C$-degree map to the boundary $D$ on $Y$. We can also impose the residual divisor as base points. We want to get a surface $S$ with CHS in this way. For this, the new resolved system must again be base-point-free. It must also have positive square.

Now we take a part of $C_1|_{W_1}$. We want a subsystem of $|L|$ such that this part is the fixed part of its restriction to the boundary. Here, we denote the boundary restriction by
\begin{center}
$Q:=C_1|_{W_1}=\Sigma+R.$
\end{center}
In this situation, $\Sigma$ is the given base divisor, and $R$ is the residual part. We choose $R$ such that $|R|$ is base-point-free on $W_1$. Let $s_\Sigma$ be a section of $\OO_{W_1}(\Sigma)$ with zero divisor $\Sigma$, and put
\begin{center}
$H_\Sigma:=s_\Sigma H^0(W_1,\OO_{W_1}(R)) \subset H^0(W_1,\OO_{W_1}(Q)).$
\end{center}
For a section $s\in H^0(W_1,\OO_{W_1}(R))$ defining a divisor $R'\in|R|$, we have that the section $s_\Sigma s$ has zero divisor $\Sigma+R'$.

Consider $X_{\mathcal K}$ as in Definition \ref{deftuple}, we have
\begin{center}
$0\longrightarrow \OO_{X_{\mathcal K}}(K_{X_{\mathcal K}}+L_{\mathcal K}) \longrightarrow \OO_{X_{\mathcal K}}(L_{\mathcal K}) \longrightarrow \OO_{W}(L_{\mathcal K}) \longrightarrow 0 .$
\end{center}
By the admissibility, $L_{\mathcal K}$ is nef and big. Hence, Kawamata--Viehweg vanishing gives us $H^1\!\left(X_{\mathcal K},\OO_{X_{\mathcal K}}(K_{X_{\mathcal K}}+L_{\mathcal K})\right)=0.$ Hence, we get the following short exact sequence on global sections
\begin{center}
$0\longrightarrow H^0\!\left(X_{\mathcal K},\OO_{X_{\mathcal K}}(K_{X_{\mathcal K}}+L_{\mathcal K})\right) \longrightarrow H^0\!\left(X_{\mathcal K},\OO_{X_{\mathcal K}}(L_{\mathcal K})\right) \xrightarrow{\rho_{W}} H^0\!\left(W,\OO_{W}(L_{\mathcal K})\right) \longrightarrow 0 .$
\end{center}
Notice that the part $H_{\Sigma}$ does not change after blow-ups, thus we get 
\begin{center}
$ H^0\!\left(X_{\mathcal K},\OO_{X_{\mathcal K}}(L_{\mathcal K})\right) =  \rho_W^{-1}\!\left(H_\Sigma \right).$
\end{center}
We put $\mathcal L_1:=\mathbb P\bigl(\rho_W^{-1}(H_\Sigma)\bigr)\subset|L|$. Thus, this implies that the base locus of $\mathcal L_1$ along $W_1$ is $\Sigma$.

We apply this with $Q$ in place of $\Sigma$ and with $0$ in place of $R$. Then the base locus along $W_1$ is all of $Q$, so we also impose the residual divisor. After we resolve $Q$, the strict transform of a general member is disjoint from the strict transform of $W_1$. Under the two conditions above, the contraction of the curves of degree zero gives the surface $S$ with CHS.
\end{proof}

Let $X_{\min}=\mathbb P_{\mathbb P^1} (\mathcal O\oplus\mathcal O(-e))$ with $e\geq1$. Let $W_1\in|-K_{X_{\min}}|$ be connected, and let $a\geq2$. For a general smooth $C_1\in|aC_\infty|$, we resolve $C_1|_{W_1}$. Let $C$ be its strict transform. Let $S$ be the surface which we get by contracting the curves of $C$-degree zero.

\begin{theo}\label{thrrulchs}
If $a((a-1)e-2)\geq2$, then $S$ is a surface with CHS, and its only non-rational singularity is minimally elliptic. In the non-hyperelliptic case, its polarizing system embeds $S$ in $\mathbb P^g$. In the hyperelliptic case, it gives a finite morphism of degree two.
\end{theo}

\begin{proof}
We apply Proposition \ref{propdec} with $\Sigma=Q:=C_1|_{W_1}$ and $R=0$. We take $H_Q$ to be the line spanned by a section of $\mathcal O_{W_1}(C_1)$ with zero divisor $Q$, and we put $\mathcal L_1:=\mathbb P(\rho_{W}^{-1}(H_Q))\subset|aC_\infty|.$ So the divisor $Q$ is the base divisor on $W_1$. By Step 2, the morphism $\varphi\colon X\longrightarrow X_{\min}$ resolves the base points of the linear system. For a smooth general member, we need $a(e+2)$ blow-ups to resolve the scheme theoretic base locus $Q=C_1|_{W_1}$. Since $Q$ is also a reduced divisor, we have $C^2=a^2e-a(e+2)=a((a-1)e-2).$ Notice that the inequality condition excludes $(a,e)=(2,1),(2,2),(3,1)$.

By Step 3, for $m\gg0$, the morphism $\varphi_{|mC|}\colon X\longrightarrow S\subset\mathbb PH^0(X,mC)$ contracts exactly the curves of polarizing degree zero. Now, let $C$ and $W$ be the strict transforms of a general member of $\mathcal L_1$ and of $W_1$ from $\varphi$ above. The connected anticanonical divisor $W$ has arithmetic genus one. On the construction stratum, its intersection matrix is negative definite. So its contraction gives one minimally elliptic singularity. Any other contracted rational curves give rational singularities. Thus $S$ has exactly one non-rational singularity.

The resolved system $\widetilde{\mathcal L}_1\subset|C|$ is base-point-free on the construction stratum. So it descends through $\varphi_{|mC|}$. We get a system $\mathcal L$ on $S$. Let $C'$ be a general member. By the adjunction argument in Steps 3 and 4, we have $\widetilde{\mathcal L}_1|_{C'}=|K_{C'}|$. Thus $\mathcal L_1$ gives the embedding of the polarized surface $S$, or the map of degree two, as in the statement.
\end{proof}
\begin{rema}\label{remrecon}
Let $(S,C)$ be a rational log canonical surface with CHS. Let $\pi:X\to S$ be its minimal resolution. Let $W$ be the reduced exceptional divisor over the non-canonical point. Then we have $K_X+W\sim0$. Let $\varphi:X\to X_{\min}$ be a sequence of blow-downs of $(-1)$-curves. We assume that it ends at one of the pairs $(X_{\min},W_1=\varphi_*W)$ of Definition \ref{deftuple}. Now let $E\not\subset W$ be a $(-1)$-curve. Since $K_X+W\sim0$, we have $W\cdot E=-K_X\cdot E=1$. So all centres of $\varphi$ lie on the boundary. For each centre $p_\alpha$, let $E^*_\alpha$ be the total transform in $X$ of its exceptional curve. We put $m_\alpha=\pi^*C\cdot E^*_\alpha$. Since $\pi^*C$ is nef, we have
\begin{center}
$\pi^*C\cdot\widetilde E_\alpha =m_\alpha-\sum_{\beta\succ\alpha}m_\beta\geq0.$
\end{center}
Hence the centres with $m_\alpha>0$ form a consistent marked cluster $\mathcal K$ on $W_1$. We observe that $\mathcal K$ is the cluster of base points of $\varphi_*|\pi^*C|$. A general member must have multiplicity $m_\alpha$ at each point of $\mathcal K$. Its strict transform is disjoint from $W$. Thus $\varphi_*|\pi^*C|$ is a $\mathcal K$-linear system with $R=0$. Notice that by Lemma \ref{lemlinear}, the centres with $m_\alpha=0$ do not change the polarizing system. Hence, we can apply Construction \ref{constr1} to $(X_{\min},W_1,\varphi_*|\pi^*C|,\mathcal K)$ to recover $(S,C)$.
\end{rema}
Now we come to the case $g(\Gamma)=1$. Here we have $-K_{X_{\min}}=C_0+C_\infty$. Every blow-up centre lies on the anticanonical boundary. So we have $K_{X_i}+W_i\sim0$ at each stage. Thus we can use the same argument as in the rational case. This gives us a similar statement for elliptic log canonical surfaces with CHS. Let $\Gamma$ be a smooth elliptic curve, and let $A$ be a divisor on $\Gamma$ with $e=\deg A\geq1$. We put $X_{\min}=\mathbb P_\Gamma(\mathcal O_\Gamma\oplus\mathcal O_\Gamma(-A))$. Let $C_0$ and $C_\infty$ be its negative and positive sections. Then we have $-K_{X_{\min}}\sim W_1:=C_0+C_\infty$. Now, We take a general smooth $C_1\in|aC_\infty|$ and resolve $C_1|_{W_1}$. Now denote by $C$ the strict transform of $C_1$. We then can contract the curves that does not intersect $C$. Let $S$ be the surface which we get from the process above, then we have
\begin{theo}\label{threllchs}
Let $a\geq2$ and $e\geq1$. The surface $S$ has canonical hyperplane sections. Its non-rational singular locus consists of two simple elliptic singularities.
\end{theo}

\begin{proof}
We put $W_1:=C_0+C_\infty\in|-K_{X_{\min}}|$, and we fix $a\geq2$ as in Construction \ref{constr1}. We choose a general smooth irreducible curve $C_1\in|aC_\infty|$. Since $C_0\cdot C_\infty=0$, this curve is disjoint from $C_0$. So $C_1\cap W_1$ is the divisor $Q:=C_1|_{W_1}=C_1|_{C_\infty}$. The restriction sequence
\begin{center}
$0 \longrightarrow \mathcal O_{X_{\min}}(C_1-C_\infty) \longrightarrow \mathcal O_{X_{\min}}(C_1) \longrightarrow \mathcal O_{C_\infty}(C_1) \longrightarrow 0$
\end{center}
gives a surjection $\rho_{C_\infty}$ from $H^0(X_{\min},\mathcal O_{X_{\min}}(aC_\infty))$ onto $H^0(C_\infty,\mathcal O_{C_\infty}(aC_\infty))$. We have
\begin{center}
$C_1-C_\infty\sim K_{X_{\min}}+\Delta+(a-1)C_\infty.$
\end{center}
Here $\Delta:=C_0+C_\infty$ is log canonical, and $(a-1)C_\infty$ is semi-ample, since $\mathcal O_{X_{\min}}(C_\infty)|_{C_0}$ is trivial. We apply the injectivity theorem \cite[Theorem 1.4]{FujinoInjectivity} to $D:=C_\infty$. We get that the map from $H^1(X_{\min},\mathcal O_{X_{\min}}(C_1-C_{\infty}))$ to $H^1(X_{\min},\mathcal O_{X_{\min}}(C_1))$ is injective. As in Theorem \ref{thrrulchs}, let
\begin{center}
$H_Q\subset H^0(C_\infty,\mathcal O_{C_\infty}(aC_\infty))$
\end{center}
be the one-dimensional subspace generated by a section with zero divisor $Q$. We put $\mathcal L_1:=\mathbb P(\rho_{C_\infty}^{-1}(H_Q))\subset|aC_\infty|$. Here, the base locus of $\mathcal L_1$ along $W_1$ is $Q$ and $\Supp Q \subset C_\infty$.

Now, we resolve $Q$ as in Proposition \ref{propdec}. On the resolved surfaces, let $C$ be the strict transform of a general member of $\mathcal L_1$ and $\widetilde C_\infty$ be the strict transform of $C_\infty$. Here, we put $W:=C_0+\widetilde C_\infty$ with all centres lie on $C_\infty$ and $C_0$ does not change. As in Theorem \ref{thrrulchs}, we have $K_X+W\sim0$ and $C\cap W=\varnothing$. There are $ae$ centres and each of them has multiplicity one. Hence, we get
\begin{center}
$C^2=a(a-1)e,\ g(C)=1+\frac{a(a-1)e}{2}.$
\end{center}
In particular, for $(a,e)=(2,1)$, we have $C^2=2$ and $g(C)=2$. Now take $m\gg0$. Then the polarizing system contracts exactly the curves of degree zero. We have $C\cdot C_0=C\cdot\widetilde C_\infty=0$. The two components of $W$ are disjoint smooth elliptic curves. Their self-intersections are negative. So by contracting them, we get two simple elliptic singularities. Since $C\cap W=\varnothing$ and $K_X+W\sim0$, adjunction gives $\mathcal O_C(C)\simeq\omega_C$. Moreover, the restriction to $C_\infty$ is surjective. Thus we have
\begin{center}
$h^0(X,\mathcal O_X(C)) =h^0(X_{\min},\mathcal O(aC_\infty))-ae+1 =2+\frac{ea(a-1)}2=g(C)+1.$
\end{center}
The restriction to $C$ has a one-dimensional kernel. This kernel is generated by the section which defines $C$. Hence its image has dimension $g(C)$. This means that the image is all of $H^0(C,\omega_C)$. So the system on $S$ restricts to the complete canonical system on $C$. Thus the morphism to $\mathbb P^g$ is an embedding when a general member is non-hyperelliptic. Otherwise, it has degree two.
\end{proof}
In the next proposition, $\Gamma$ is an elliptic curve, and $A$ is an effective divisor on $\Gamma$ of degree at least one.
\begin{prop}\label{prop2ell}
Every surface with CHS and exactly two simple elliptic singularities comes from $\mathbb P_\Gamma(\mathcal O_\Gamma\oplus\mathcal O_\Gamma(-A))$.
\end{prop}
\begin{proof}
Let $(S,C')$ be such a surface, and let $X\to S$ be its minimal resolution. By Construction \ref{constr1}, or by Epema's description of the ruled cases, the pair $(X,C')$ has a birational morphism to a relatively minimal ruled surface
\begin{center}
$X_{\min}=\PP_\Gamma\!\left(\OO_\Gamma\oplus \OO_\Gamma(-A)\right)$
\end{center}
over an elliptic curve $\Gamma$. Here $A$ is effective with $\deg(A)\geq1$. Under this birational morphism, the hyperplane system corresponds to a subsystem $\mathcal L_1\subset |aC_\infty|$. The base scheme of $\mathcal L_1$ is supported on the anticanonical divisor. It is exactly the cluster of Construction \ref{constr1}.

Conversely, we take such a ruled surface and a general curve $C_1\in |aC_\infty|$, and we apply Theorem \ref{threllchs}. We resolve the given anticanonical cluster, and we contract its transformed boundary. In this way, we get a surface with CHS and two simple elliptic singularities.
\end{proof}

\begin{proof}[Proof of Theorem \ref{thrconst}]
This is a direct consequence of Proposition \ref{propdec} and Theorem \ref{threllchs}.
\end{proof}

\subsection{Log Calabi-Yau surfaces with fibration structure}
\label{sec33} We consider a one-dimensional subsystem $\mathcal P\subset|C_\infty|$ on $X_{\min}$ which contains the boundary component $C_\infty$. Suppose that $\Bs(\mathcal P)\subseteq\Bs(\mathcal L_1)$. Then by Construction \ref{constr1}, the reconstruction steps also give a base-point-free pencil. We denote this pencil again by $\mathcal P$. Thus we get a natural fibration over $\mathbb P^1$.

\begin{construction}\label{constr2}
We have the following diagram.
\begin{center}
\begin{tikzcd}
X_{\min} & \arrow{l}{\phi} X \arrow{r}{\pi} \arrow[swap]{d}{\phi_{\mathcal P}} & Y \arrow[dashed]{r}{\Psi} & S \ni x \arrow{r}{i} & \mathbb{P}^g \\& \mathbb{P}^1
\end{tikzcd}
\end{center}
The notation is as follows.
\begin{enumerate}
    \item $X_{\min}$, $X$, $\pi$, $\phi$, $\Psi$, $Y$, $i$ and $x \in S \subset \mathbb{P}^g$ are as in Construction \ref{constr1}.
    \item $\phi_{\mathcal P}$ is the morphism defined by the pencil $\mathcal P$. It is an elliptic fibration when $\mathcal P$ is the anticanonical pencil on a rational surface, or the positive-section pencil on an elliptic ruled surface. It is a $\mathbb{P}^1$-fibration when $\mathcal P\subset|C_\infty|$ is the rational pencil of Case 2 below.
\end{enumerate}
\end{construction}

Let $X_{\min}=\mathbb P_\Gamma(\mathcal O_\Gamma\oplus\mathcal O_\Gamma(-A))$, where $g(\Gamma)\in\{0,1\}$ and $\deg A\geq1$. We fix an admissible tuple and a pencil $\mathcal P\subset|M|$. Here $M=-K_{X_{\min}}$ in the rational anticanonical case, or $M=C_\infty$. We require three things. The general member of $\mathcal P$ is smooth and irreducible. Its cluster of base points is contained in $\mathcal K$. Every curve contracted by the polarization is vertical for the resolved pencil.

\begin{prop}\label{proppenc}
The pencil induces a fibration $f:Y\to\mathbb P^1$ on the stable log Calabi-Yau surface $((Y,D),C)$. The fibration is elliptic for an anticanonical pencil on a rational surface or a positive-section pencil on an elliptic ruled surface, and rational for a pencil in $|C_\infty|$ on $\mathbb F_e$. If $R=0$, then $Y$ is a surface with CHS.
\end{prop}
\begin{proof}
We divide the proof into two cases.

\textbf{Case 1.} We take $\mathcal L$ to be the system of Construction \ref{constr1}. We denote the boundary part of its base divisor by $Q$. In the rational ruled case, we choose a smooth elliptic $E\sim-K_{X_{\min}}$. For $X_{\min}=\F_n$ with $n=1,2$, we have $E^2=8$, and $|E|$ is base-point-free. The del Pezzo case is the same. Let $s_E\in H^0\!\left(X_{\min},\OO_{X_{\min}}(E)\right)$ be the section whose zero divisor is $E$. Then the restriction to $E$ gives the exact sequence
\begin{center}
$0\longrightarrow \mathbb C s_E \longrightarrow H^0\!\left(X_{\min},\OO_{X_{\min}}(E)\right) \xrightarrow{\rho_E} H^0\!\left(E,\OO_E(E)\right) \longrightarrow 0 .$
\end{center}
We choose a divisor $T_E\in |\OO_E(E)|$ with $T_E\leq Q$. Let $H_{T_E}\subset H^0(E,\OO_E(E))$ be the line generated by a section with zero divisor $T_E$. Then
\begin{center}
$\mathcal P_{T_E}:=\PP\!\left(\rho_E^{-1}(H_{T_E})\right)\subset |E|$
\end{center}
is a pencil of genus-one curves with base scheme $T_E$. In the elliptic ruled case, we denote the initial surface by
\begin{center}
$X_{\min}=\PP_E(\OO_E\oplus\OO_E(-A))\ \text{with} \deg A=e.$
\end{center}
Its positive section is $C_\infty$, and we have $A=\OO_{C_\infty}(C_\infty)$. Let $s_\infty$ be the section with divisor $C_\infty$. Then the restriction sequence is
\begin{center}
$0\longrightarrow \mathbb C s_\infty \longrightarrow H^0\!\left(X_{\min},\OO_{X_{\min}}(C_\infty)\right) \xrightarrow{\rho_\infty} H^0(C_\infty,A) \longrightarrow 0 .$
\end{center}
Let $T\in |A|$ be a divisor contained in $Q$, and let $H_{T}\subset H^0(C_\infty,A)$ be the line of $T$. We put
\begin{center}
$\mathcal P_T:=\PP\!\left(\rho_\infty^{-1}(H_{T})\right)\subset |C_\infty|.$
\end{center}
This is a pencil with base divisor $T$. Its base points have the multiplicities of $T$. In both cases, we resolve the cluster of base points. Then we get an elliptic fibration. The contracted curves are vertical for this fibration. So it descends through the chosen contractions.\\
\\
\textbf{Case 2.} Let $X_{\min}\simeq\F_e$. We choose an effective divisor
\begin{center}
$T=\sum_i k_i p_i\in |\OO_{C_\infty}(C_\infty)|.$
\end{center}
We choose it such that its cluster of base points is contained in $\mathcal K$. Here we have $\sum_i k_i=e$. Let $s_\infty$ be the section which defines $C_\infty$. We take $t_T\in H^0(C_\infty,\mathcal O_{C_\infty}(C_\infty))$ with zero divisor $T$. Since $H^1(\mathbb F_e,\mathcal O)=0$, we can lift $t_T$ by the restriction sequence. We get a lift $\widetilde t_T\in H^0(X_{\min},\mathcal O(C_\infty))$. Now we put
\begin{center}
$\mathcal P_T:=\mathbb P\langle s_\infty, \widetilde t_T\rangle\subset|C_\infty|.$
\end{center}
Notice that if we change the lift, we only add a multiple of $s_\infty$. Thus the pencil does not change, and we have that $\Bs(\mathcal P_T)=T$. We resolve this base scheme. In this way, we get the rational pencil. The fibres of $\pi$ are connected. We also have $\pi_*\mathcal O_X=\mathcal O_Y$. Thus the pencil factors uniquely through a morphism $Y\to\mathbb P^1$.

\end{proof}

\begin{proof}[Proof of Theorem \ref{thrfib}]
For $g\geq3$ and $e=g-1$, Theorem \ref{threllchs} gives a surface $S$ with CHS of genus $g$. The disjoint elliptic curves $C_0$ and $\widetilde C_\infty$ both have square $-e$. They are contracted to the two simple elliptic points of $S$. By Proposition \ref{proppenc}, we can construct these surfaces together with an elliptic pencil on $S$. This pencil gives the elliptic fibration.
\end{proof}
\subsection{Log Calabi-Yau surfaces with $(-2)$-curves}\label{sec34}
\subsubsection{Boundary clusters}\label{sec341}

Let $B\subset X$ be a curve which is smooth at $p$. Let $(u,v)$ be local coordinates such that $B=(v=0)$. As a subscheme of $X$, the divisor $2p\subset B$ has the ideal $I_{2p/B}=(v,u^2)$. Its resolution gives two centres of multiplicity one $(p,1),(p',1)$ with $p'\succ p$. We call this a cluster of type $(1,1)$.

Let $X$ be a smooth projective surface, and let $B\subset X$ be a smooth elliptic curve or a nodal curve of arithmetic genus one. For a line bundle $L$, we put $N=L|_B$ and $d=\deg N>0$, and we assume the following two conditions.
\begin{enumerate}
    \item The map
    \refstepcounter{equation}\label{eqrestr}
    \begin{center}
    $\rho_L:H^0(X,L)\longrightarrow H^0(B,N)$\enspace\textnormal{(\theequation)}
    \end{center}
    is surjective.
    \item The map
    \refstepcounter{equation}\label{eqnormal}
    \begin{center}
    $H^0(X,L(-B))\longrightarrow L(-B)|_p$\enspace\textnormal{(\theequation)}
    \end{center}
    is nonzero for every $p\in B^\circ$, where $B^\circ\subset B_{\mathrm{sm}}$ is a fixed open subset.
\end{enumerate}
We choose distinct points $p_1,\ldots,p_k\in B^\circ$, and we put $D=\sum_i p_i$. We assume that $N(-2D)$ has a section which does not vanish at any $p_i$. We put $\Lambda_D=\mathbb P(\rho_L^{-1}H^0(B,N(-2D)))\subset|L|$.
\begin{prop}\label{propdblpt}
At each $p_i$, the system $\Lambda_D$ has the base ideal $(v,u^2)$ and the cluster $(p_i,1),(p_i',1)$ with $p_i'\succ p_i$. Its members which are singular at $p_i$ form a hyperplane $\Sigma_i$. A general member is smooth at every $p_i$.
\end{prop}

\begin{proof}
Let $s_0\in H^0(B,N(-2D))$ be a section which does not vanish at any $p_i$. Inside $N$, this section has an equation of the form $u^2\cdot(\text{unit})$ on $B$. By \eqref{eqrestr}, it lifts to a section on $X$ with an equation $u^2a(u)+v\,b(u,v)$, where $a(0)\neq0$. By \eqref{eqnormal}, we also have a section of $L(-B)$ with an equation $v\,h(u,v)$ such that $h(0,0)\neq0$. Thus the system generates the ideal $(u^2,v)$. The condition $\mult_{p_i}(C)\geq2$ is a nontrivial linear condition on $\Lambda_D$. This means that $\Sigma_i$ is a hyperplane and a finite union of proper hyperplanes cannot cover $\Lambda_D$. Hence, the complement
\begin{center}
$\Lambda_D \setminus \bigcup_{i=1}^k\Sigma_i$
\end{center}
is nonempty and open. This implies that every member of it is smooth at all points $p_i$.
\end{proof}
The assumption on $N(-2D)$ is a condition on the position of the points. We describe it in the next lemma. Let $B_j$ be the components of $B$. Let $d_j=\deg(N|_{B_j})$ be the degree of $N$ on $B_j$. We denote by $n_j$ the number of chosen points on $B_j$, and let $k=\sum_jn_j\geq1$. In the case $2n_j=d_j$ for all $j$, we impose $2D\sim N$. When $k\geq2$, we can choose such collections of points in $B^\circ$. For an irreducible curve $B$, this gives the range $1\leq k\leq\lfloor d/2\rfloor$, except possibly for $k=d/2=1$.
\begin{lemm}\label{lem-blt}
If a section of $N(-2D)$ does not vanish at any point $p_i$, then we have $2n_j\leq d_j$ for every $j$. Conversely, let $p_i$ be general distinct points of $B^\circ$ satisfying the extreme-case condition above. If those inequalities above hold, then $N(-2D)$ has such section.
\end{lemm}
\begin{proof}
A section that does not vanish at a point of $B_j$ is nonzero on this component. Hence the bundle $N(-2D)|_{B_j}$ has degree $d_j-2n_j\geq0$.

For the converse direction, suppose that $d_j-2n_j>0$ for some $j$. If $B$ is irreducible, we must have that a general member of $|N(-2D)|$ avoids general points $p_i$. If $B$ is reducible, then all its components are rational. On the components with $d_j-2n_j=0$, we can see that the bundle $N(-2D)$ is trivial. Hence a section that is nonzero at the nodes extends across these components. On a component of positive degree, we can choose the values at the two nodes freely. Thus the sections glue around the cycle, and we can choose them nonvanishing at general $p_i$.

If $d_j=2n_j$ for every $j$, then the bundle $N(-2D)$ has multidegree zero. In this case, a nonvanishing section exists when $N(-2D)\simeq\mathcal O_B$, i.e., when $2D\sim N$. Because multiplication by two is surjective on the generalized Jacobian of $B$, we see that such collections exist. When $k\geq2$, the finitely many choices for the last point vary with the other points and hence we can take all the points distinct and inside $B^\circ$.
\end{proof}

\subsubsection{Application to rational and elliptic ruled surfaces}\label{sec342}

Let $a\geq2$. We consider the following three cases and put $d:=L\cdot B$ in each of them.\\
\\
\textup{(RNS) Rational ruled case with reducible boundary.} Let $X=\mathbb F_e$ with $e\in \mathbb N$, and let $H=C_\infty$ be the positive section. Let $B$ be a reducible nodal anticanonical curve. We take $B=H+C_0+F_1+F_2$ or $B=H'+F+C_0$. Here $F$ is a fibre, and $H'$ is a general member of $|H+F|$. We put $L=\mathcal O_X(aH)$. Then we have $d=a(e+2)$.\\
\\
\textup{(RS) Rational ruled case with irreducible boundary.} Let $X=\mathbb F_e$ with $e\in\{1,2\}$, and let $H=C_\infty$ be the positive section. Let $B\in|-K_X|$ be an integral anticanonical curve. We will assume that $B$ is smooth elliptic or nodal. Here, put $L=\mathcal O_X(aH)$. Then we have $d=a(e+2)$.\\
\\
\textup{(E) Elliptic ruled case.} Let $X=\mathbb P_E(\mathcal O_E\oplus\mathcal O_E(-A))$, where $E$ is an elliptic curve and $\deg A=e\geq 1$. The two sections give us an anticanonical divisor $-K_X\sim C_0+C_\infty$. We take $B=C_\infty$ and $L=\mathcal O_X(aC_\infty)$. This gives us that $d=ae$.

We will give the ranges for the possible cases. In case \textup{(RNS)}, we exclude $(a,e)=(2,1)$. On each component $B_j$, we allow $n_j\geq0$ clusters when $2n_j\leq\deg(L|_{B_j})$. In particular, on the positive section, we allow $1\leq k\leq\lfloor ae/2\rfloor$. In cases \textup{(RS)} and \textup{(E)}, $1\leq k\leq\lfloor d/2\rfloor$ is possible. The only exception can be $k=d/2=1$. In case \textup{(RS)}, we will also assume $a+e\geq4$, of which will be needed for \eqref{eqnormal}.

\begin{prop}\label{proprulcl}
In cases \textup{(RNS)}, \textup{(RS)} and \textup{(E)}, in the ranges above, there is a system in $|L|$ with $k$ clusters of base points of type $(1,1)$ on $B_{\mathrm{sm}}$. Its general member is smooth at the proper points of these clusters.
\end{prop}

\begin{proof}
We check the assumptions of Proposition \ref{propdblpt}. In the rational cases, condition \eqref{eqrestr} follows from the restriction sequence in Step 4 of the proof of Proposition \ref{propdec}. In the elliptic ruled case, it follows from the restriction sequence in the proof of Theorem \ref{threllchs}. The ranges for $k$ come from Lemma \ref{lem-blt}.\\
\\
\textbf{Case 1.} We consider the rational ruled surface with irreducible boundary. If $e=2$, then we have
\begin{center}
$L(-B)\sim(a-2)H.$
\end{center}
This bundle is generated away from the negative section. It is trivial for $a=2$. If $e=1$, then we have
\begin{center}
$L(-B)\sim C_0+(a-3)H.$
\end{center}
Now let $a+e\geq4$. For $e=1$, this means $a\geq3$. Then $L(-B)$ has a section which does not vanish at a general point of $B_{\mathrm{sm}}\setminus C_0$. So \eqref{eqnormal} holds with $B^\circ$ a dense open subset of $B_{\mathrm{sm}}$.\\
\\
\textbf{Case 2.} We consider the elliptic ruled surface. We use
\begin{center}
$L(-C_\infty)=\mathcal O_X((a-1)C_\infty).$
\end{center}
We put $\mathcal A:=\mathcal O_{C_\infty}(C_\infty)$. This bundle has degree $e$. We have $\xi_*\mathcal O_X((a-1)C_\infty)\simeq\bigoplus_{i=0}^{a-1}\mathcal A^{i}$. So the restriction to $C_\infty$ is surjective onto $H^0(C_\infty,\mathcal A^{a-1})$. The degree of $|\mathcal A^{a-1}|$ is $(a-1)e$. This degree is at least two, unless $(a,e)=(2,1)$. If it is at least two, then the system is base-point-free. So \eqref{eqnormal} holds with $B^\circ=C_\infty$. For $(a,e)=(2,1)$, let $p_0$ be the unique member of $|\mathcal A|$. Then \eqref{eqnormal} holds with $B^\circ=C_\infty\setminus\{p_0\}$.\\
\\
\textbf{Case 3.} We consider the reducible nodal boundary. Both choices of $B$ lie in $|-K_X|$. We use the basis $C_0,f$ with $H\sim C_0+ef$. Then we have
\begin{center}
$L(-B)\sim aH+K_X\sim(a-2)C_0+(ae-e-2)f.$
\end{center}
This divisor is effective if and only if $(a,e)\neq(2,1)$. Now we take the sections $s_{C_0}^{\,a-2}\cdot\xi^*t$ with $t\in H^0(\mathbb P^1,\mathcal O(ae-e-2))$. They vanish only along $C_0$ and along finitely many fibres. So they are nonzero at a general point of every component of $B$ other than $C_0$. Hence \eqref{eqnormal} holds with $B^\circ$ a dense open subset of $B_{\mathrm{sm}}\setminus C_0$. Notice that we have $\deg(L|_{C_0})=aH\cdot C_0=0$. Hence, there is no cluster on $C_0$ by Lemma \ref{lem-blt}.
\end{proof}
In the following theorem, let $a\geq2$, $e\geq2$ and $q\in\{0,1\}$ be integers with $a^2e-a(e-2q+2)\geq2$.
\begin{theo}\label{thrm2crv}
There is a surface with CHS such that its minimal resolution has $e$ disjoint internal $(-2)$-curves. Each of these curves is disjoint from a general polarizing member.
\end{theo}
\begin{proof}
First, we choose reduced divisors $T=p_1+\cdots+p_e\in|\mathcal A|$ and $R_a\in|\mathcal A^{a-2}|$. They can be chosen such that $R_a$ is disjoint from $T$. For $a=2$, we have $R_2=0$. Put $Z_a=2T+R_a\in|\mathcal A^a|$. Now we consider the system
\begin{center}
$\Lambda_{Z_a}=\mathbb P(\rho_a^{-1}(\Bbbk s_{Z_a})) \subset|aC_\infty|.$
\end{center}
We observe that it has $e$ chains $(1,1)$ and $(a-2)e$ proper simple base points. Indeed, by Proposition \ref{propdblpt}, the local ideal at each $p_i$ is $(v,u^2)$. Since $T$ is general and $R_a\cap T=\varnothing$, a general member is smooth at every $p_i$ and has an infinitely near simple base point for each $2p_i$. In the elliptic ruled case, there are $2e+(a-2)e=ae$ simple centres. Hence we have $\sum_\alpha m_\alpha=ae$ and $m_\alpha=1\leq a-1$. Thus, the resolved strict transform $C$ satisfies
\begin{center}
$C^2=(aC_\infty)^2-ae=a(a-1)e$
\end{center}
and $g(C)=1+C^2/2=1+a(a-1)e/2=\dim\Lambda_{Z_a}$. For a rational reducible boundary, we also add the reduced base points on its other components. The degrees of $L$ give the numbers of these points. In this case, we have $C^2=a^2e-a(e+2)$. By our assumption on $a$ and $e$, we get $g(C)=1+C^2/2\geq2$.

We lift a section of $\mathcal A$ which vanishes on $T$. We get a pencil $\mathcal P=\langle C_\infty,H_T\rangle\subset|C_\infty|$ with $\Bs(\mathcal P)=T$. This implies that the length-two clusters lie at the base points of $\mathcal P$ for every $a\geq2$. Let $F$ be the fibre class of the resolved pencil. Over $p_i$, we have two blow-ups. Let $E_i^*$ and $E_i'^*$ be their total exceptional transforms. The strict transform of the first exceptional curve is $G_i=E_i^*-E_i'^*$. We have
\begin{center}
$G_i^2=-2,\ C\cdot G_i=1-1=0,\ F\cdot G_i=1-0=1.$
\end{center}
Since $G_i$ is disjoint from $C$, we have $|C|$ contracts the $(-2)$-curves to $A_1$ points and  finish the proof.
\end{proof}
\begin{rema}
Notice that both centres occur in $\mathcal K$ with multiplicity one for $C$. Only the first one is a base point of the pencil. Thus $(C+F)\cdot G_i=1$. Thus the polarization $C+F$ does not contract $G_i$. Since $F^2=0$ and $K_X+W\sim0$, adjunction gives
\begin{center}
$(C+F)\cdot W=F\cdot W=2-2q.$
\end{center}
For $q=1$, the system $|C+F|$ still contracts the boundary. It gives us a surface with CHS with $(-2)$ curves on it. For $q=0$ case, we get a stable log Calabi-Yau polarization with non-trivial boundary. For example, let $W_1=C_0+C_\infty+F_1+F_2$. Then the strict transforms of $F_1$ and $F_2$ have intersection number one with $C+F$. Hence, they are not contracted. 
\end{rema}
\subsection{Elliptic fibrations}\label{sec35}

Let $(S,C,\mathcal P)$ be such a surface. We denote by $H$ the class of the pencil on $X_{\min}$. We denote by $\mathcal K_T$ its cluster of base points. On the blow-up of $\mathcal K$ with total exceptional transforms $E_\alpha^*$, we have
\begin{center}
$F_T=\mu^*H-\sum_{\alpha\in\mathcal K_T}E_\alpha^*, \ C=\mu^*L-\sum_\alpha m_\alpha E_\alpha^*.$
\end{center}
Thus the fibre degree is
\begin{center}
$N_{\mathrm{ell}}:=C\cdot F_T =L\cdot H-\sum_{\alpha\in\mathcal K_T}m_\alpha>0.$
\end{center}
Here $H=C_\infty$ on an elliptic ruled surface, and $H=-K_{X_{\min}}$ for the anticanonical pencil on a rational surface. For $L=aC_\infty$, the first intersection number is $ae$ in the elliptic ruled case. It is $a(e+2)$ on $\mathbb F_e$ in the anticanonical case. In particular, in the elliptic ruled case we have $N_{\mathrm{ell}}=(a-1)e$ only when the polarization has multiplicity one at all pencil centres.

We choose a smooth fibre $E\in\mathcal P$ such that $D:=E\cap C$ is reduced of degree $N_{\mathrm{ell}}$. We blow up $D$ by $\mu:S'\longrightarrow S$. Let $G_1,\ldots,G_{N_{\mathrm{ell}}}$ be the exceptional curves. Then we have
\begin{center}
$K_{S'}\sim\sum_{i=1}^{N_{\mathrm{ell}}}G_i.$
\end{center}
Moreover, the strict transform $\widetilde E$ satisfies
\begin{center}
$\widetilde E+\sum_{i=1}^{N_{\mathrm{ell}}}G_i \in\mathcal P'.$
\end{center}

The contraction defined by the transformed CHS polarization
\begin{center}
$\pi_{|aC'|}:S'\longrightarrow S_{\mathrm{ell}}$
\end{center}
contracts $\widetilde E$. We denote this contraction by $\pi=\pi_{|aC'|}$. Let $f_S:S\to\mathbb P^1$ be the morphism defined by $\mathcal P$. The pencil $\mu^*\mathcal P$ gives $f_{\mathrm{ell}}:S_{\mathrm{ell}}\to\mathbb P^1$ with $f_{\mathrm{ell}}\circ\pi=f_S\circ\mu$. Let $F$ be the fibre class of $f_{\mathrm{ell}}$. The contraction formula gives
\begin{center}
$\pi^*K_{S_{\mathrm{ell}}} =K_{S'}+\widetilde E \sim\sum_iG_i+\widetilde E=\mu^*E=\pi^*F.$
\end{center}
Hence we have $K_{S_{\mathrm{ell}}}\sim F$. The fibres are connected. Hence we get $h^0(S_{\mathrm{ell}},mK_{S_{\mathrm{ell}}})=m+1$ for $m\geq1$. Thus $\kappa(S_{\mathrm{ell}})=1$, and the canonical system defines $f_{\mathrm{ell}}$.

More generally, we can choose fibres $E^{(1)},\ldots,E^{(n)}\in\mathcal P$, and we can blow up their intersections with $C$. Then we get a surface $S_{\mathrm{ell}}$ with $K_{S_{\mathrm{ell}}}\sim nF,$ where $F$ is the fibre class of $f_{\mathrm{ell}}$.
\begin{construction}\label{constr3}
We have the following diagram.
\begin{center}
\begin{tikzcd}
& S' \arrow{r}{\pi} \arrow[swap]{d}{\mu} & S_{\mathrm{ell}}    \\& (S,E),C \ni x \arrow{r}{i} & \mathbb{P}^g
\end{tikzcd}
\end{center}
The notation is as follows.
\begin{enumerate}
    \item $(S,C)$ is a surface with CHS, $E$ is a fibre of the elliptic fibration, $S^{\prime}$ is the blow-up of $S$, and $S_{\mathrm{ell}}$ is the resulting elliptic surface with simple elliptic singularities.
    \item $S$, $x$ and $i$ are as in Construction \ref{constr1}. The blow-up $\mu$ and the contraction $\pi$ are defined above.
    \item $\phi_{\mathcal P}$ is the morphism defined by the elliptic pencil $\mathcal P$ of Construction \ref{constr2}.
\end{enumerate}
\end{construction}
We can replace $C$ by an effective divisor $A$ which meets $E$ transversely in the smooth locus. Then we use the section ring of $K_{S'}+\widetilde E+tA'$. Here $A'$ is the strict transform of $A$. For this, we need that $(S',\widetilde E+tA')$ is log canonical and that this divisor is nef and big for small $t>0$.
\begin{prop}
The surface $S_{\mathrm{ell}}$ of Construction \ref{constr3} is projective.
\end{prop}

\begin{proof}
First, we blow up the reduced divisor $D=E\cap C$ of degree $N_{\mathrm{ell}}$ by $\mu:S'\longrightarrow S$. We write $G_1,\ldots,G_{N_{\mathrm{ell}}}$ for the exceptional curves. We write $\widetilde E$ and $C'$ for the strict transforms. We have $K_{S'}=\mu^*K_{S}+\sum_iG_i\sim\sum_iG_i$ and $\mu^*E=\widetilde E+\sum_iG_i$. These two identities give $K_{S'}+\widetilde E\sim\mu^*E$. Here the curve $\widetilde E$ is still smooth elliptic and we have $\widetilde E^2=E^2-N_{\mathrm{ell}} =-N_{\mathrm{ell}}<0$. By Artin's theorem, the contraction $\pi:S'\longrightarrow S_{\mathrm{ell}}$ gives a simple elliptic point. Hence, we have $\pi^*K_{S_{\mathrm{ell}}} =K_{S'}+\widetilde E\sim\mu^*E$. After the blow-ups, $\widetilde E$ is disjoint from $C'$. Thus, $C'$ descends to $C_{\mathrm{ell}}$ with $\pi^*C_{\mathrm{ell}}=C'$.

The pullback of $K_{S_{\mathrm{ell}}}+\epsilon C_{\mathrm{ell}}$ is $K_{S'}+\widetilde E+\epsilon C' \sim F'+\epsilon C'$. Here $F':=\mu^*E$, and we have $F'\cdot C'=\mu^*E\cdot C'=E\cdot C=N_{\mathrm{ell}}$. Thus we get
\begin{center}
$(F'+\epsilon C')^2 = 2\epsilon F'\cdot C'+\epsilon^2(C')^2 = 2N_{\mathrm{ell}}\epsilon+\epsilon^2(C')^2>0$
\end{center}
for small $\epsilon>0$. Now let $B\subset S_{\mathrm{ell}}$ be an irreducible curve with strict transform $B'\ne\widetilde E$. If $B'$ is horizontal, then $F'\cdot B'>0$. So $(F'+\epsilon C')\cdot B'>0$ for small $\epsilon$. If $B'$ is vertical, then $F'\cdot B'=0$. In this case, either $B'=G_i$ with $C'\cdot G_i=1$, or $B'$ lies in another fibre. In the second case, we have $C'\cdot B'=C\cdot\mu_*B'>0$, because $C$ is ample on $S$. Thus we have $(F'+\epsilon C')\cdot B'>0$ in every case. By Nakai--Moishezon on normal surfaces, we get that $K_{S_{\mathrm{ell}}}+\epsilon C_{\mathrm{ell}}$ is ample for $0<\epsilon\ll1$. Hence $S_{\mathrm{ell}}$ is projective.
\end{proof}
\begin{rema}
For $q=0$ or $q=-\infty$ with reducible boundary, this construction is not possible. Let $S$ be a rational surface with CHS which has only a non-simple elliptic singularity. Then $S$ has no fibration structure. Indeed, the pencil has to meet the boundary before the blow-down step. So it is not base-point-free.
\end{rema}
In the following construction, we let the boundary vary in its complete system. The CHS case is one of these choices.

\begin{construction}\label{constr4}
Let $Y_0$ be a smooth rational surface with a morphism $f_0:Y_0\to\mathbb P^1$. We assume that $f_0$ is the ruling of a Hirzebruch surface. We also allow that $f_0$ comes from a pencil on a Hirzebruch surface, after we resolve its base points. First, we choose a reduced boundary $B_0$ and an effective Cartier divisor $A_0$ with $K_{Y_0}+B_0\sim\beta F$. Here $\beta>0$, and $F$ is a fibre. Then, we blow up some given transverse intersection points of $A_0$ with the smooth locus of $B_0$ and denote the strict transforms by $A_Y$ and $B_Y$. Notice that here we also allow infinitely near boundary centres. Their weights are the effective multiplicities of the strict transform of $A_0$ at each step and their proximity relations are part of the cluster data. We require: The pair $(Y,B_Y+tA_Y)$ is lc for small $t>0$, $d=A_Y\cdot F>0$ and $A_Y\cdot V\geq0$ for every vertical curve $V$. In the end, we get the birational contraction $\rho:Y\to S$ which is defined by $\beta F+tA_Y$. 
\end{construction}

\begin{prop}\label{proprulf}
Construction \ref{constr4} gives a stable minimal model with $K_S+B\sim_{\mathbb Q}\beta F$. Its polynomial is $2\beta d\,t+A_Y^2t^2$. The contraction $\rho$ contracts the curves which do not meet $A_Y$.
\end{prop}
\begin{proof}
The blow-up centres lie on the boundary. So we have $K_Y+B_Y\sim\beta F$. Let $t>0$ be small. We check that $\beta F+tA_Y$ is nef. Its intersection numbers with vertical curves are nonnegative. Now let $G$ be a curve with $A_Y\cdot G<0$. Then $G$ is one of the finitely many horizontal components of $A_Y$. On such a curve, the term $\beta F$ is positive. So the divisor $\beta F+tA_Y$ is nef. Its square is $2\beta dt+t^2A_Y^2>0$. Hence by surface log abundance \cite[Theorem 6.1]{FujinoSurfaces12}, this divisor is semiample. Its zero curves are the vertical curves which do not meet $A_Y$. Let $\rho:Y\to S$ be the contraction. Then the pencil defines a morphism $f:S\to\mathbb P^1$. We have $K_Y+B_Y+tA_Y=\rho^*(K_S+B+tA)$. This implies that $K_S+B+tA$ is ample. It also implies that $(S,B+tA)$ is lc and that $A^2=A_Y^2$.
\end{proof}

\paragraph{Smooth rational surfaces from Hirzebruch surfaces.}\label{parhirz} Let $C_0$ be the negative section of $\mathbb F_n$, and let $f$ be a fibre of the ruling. We choose the boundary in the linear system
\begin{center}
$B_0\in|2C_0+(n+\beta+2)f|.$
\end{center}
With this choice, we have $K_{\mathbb F_n}+B_0\sim\beta f$. Among the admissible reduced members, we choose the following two smooth families. For $n=0$, we denote the two ruling classes by $F$ and $h$. First, we choose $d\geq2$ distinct sections $A_i\in|h|$ and a smooth $B_0\in|3F+2h|$ and then blow up all the $3d$ transverse intersections of $B_0$ with $\sum_iA_i$. By this way, we get a smooth surface $Y$ with $K_Y+B_Y\sim F$. The divisor $A_Y=\sum_iA_i'$ consists of disjoint $(-3)$-sections. Now for $d=1$, we can take $Y=\mathbb F_3$ and $A_Y=C_0$ together with a smooth $B_Y\in|2C_0+6f|$ which is disjoint from $C_0$. Hence, in both cases, we have $A_Y\cdot F=d$ and $A_Y^2=-3d$. Thus, we get
\begin{center}
$\vol(K_Y+B_Y+tA_Y)=2dt-3dt^2.$
\end{center}
\section{Moduli spaces of stable log Calabi-Yau surfaces}
\label{sec4}

In this section, we construct coarse moduli spaces for the stable log Calabi-Yau surfaces which come from an admissible tuple $\Upsilon$. These surfaces are stable minimal models in the sense of \cite{Bir2}. So we identify their locus in the strongly embedded space $\mathrm E^s\mathrm{SLCY}_{\mathrm{lc}}\Xi$. This is our notation for $M^e_{\mathrm{lc}}\Xi$ in \cite[Section 9 and Proposition 9.5]{Bir2}. Here
\begin{center}
$\Xi=(2,\Phi,\mathcal V_{\mathrm{fib}},\sigma(t),\tau,r,N)$
\end{center}
consists of the dimension, the coefficient set $\Phi$, the permitted fibre volumes $\mathcal V_{\mathrm{fib}}$, the volume polynomial $\sigma(t)=\vol(K_X+D+tC)$ and the embedding data $\tau,r,N$. We fix a sufficiently small rational $\tau>0$ and an integer $r>0$ such that $r(K_X+D+\tau C)$ is Cartier and very ample with $N+1$ sections. Thus the embedding is in $\mathbb P^N$.

The group $G_N=\mathrm{PGL}_{N+1}$ changes the projective coordinates. The quotient
\begin{center}
$q_\Xi:\mathrm E^s\mathrm{SLCY}_{\mathrm{lc}}\Xi \longrightarrow\mathrm{SLCY}_{\mathrm{lc}}\Xi$
\end{center}
forgets these coordinates, see \cite[Section 10.7]{Bir2}. We prove that the construction locus has a finite partition into $G_N$-invariant locally closed strata. The image of each stratum under $q_\Xi$ is its quasi-projective coarse moduli space. We denote this space by $\mathrm{SLCY}_{\mathrm{lc}}(\Upsilon,\nu)$, where $\nu$ is the label of the stratum. Then we parametrize the admissible pencils. In this way, we get the coarse spaces $\mathrm{LCYF}_{\mathrm{lc}}\Upsilon$ with a chosen fibration. We apply the same method to the elliptic construction, and we get $\mathrm{CHS}^{\Upsilon_\ell}_\ell$.

In this section, $Y_0$ is the initial surface $X_{\min}$ of Section \ref{sec3}. It comes with the anticanonical boundary $D_0$ and the polarizing line bundle $L$. In the ruled CHS construction, we have $L=\mathcal O_{Y_0}(aC_\infty)$. We denote by $Y^{\mathrm{lc}}$ the extra log canonical modification. We denote by $\widetilde{\mathcal Y}_\Upsilon$ the blow-up of the marked cluster. The contraction defined by the transformed polarization gives the stable surface. We give the tuple $\Upsilon$ in Subsection \ref{sec42}. We denote $\mathcal D_\Upsilon$ for the boundary of the stable family and $J_\Upsilon$ for the image of the construction in the ambient strongly embedded space.

For a fixed tuple $\Upsilon$, we will get the parameter spaces are finite unions of locally closed strata. Thus by Chevalley's theorem, their images in the ambient space are constructible. In the decomposable elliptic ruled case, we get a stronger result that these images are locally closed by Proposition \ref{propelloc}. For the rational case, we use the finite partition of Definition \ref{defstrata}. The key idea for us is the use Construction \ref{constr1} with $C_1|_{D_0}=\Sigma+R$. Notice that here $R$ is the residual intersection on the resolved surface and we have $R=0$ in the CHS case.

\subsection{Ambient strongly embedded moduli spaces}
\label{sec41}

The space $\mathrm E^s\mathrm{SLCY}_{\mathrm{slc}}\Xi$ parametrizes the strongly embedded $\Xi$-stable slc surfaces. Its normal lc locus $\mathrm E^s\mathrm{SLCY}_{\mathrm{lc}}\Xi$ is locally closed. We denote the universal family by
\begin{center}
$\pi_\Xi:\bigl((\mathcal X_\Xi,\mathcal D_\Xi),\mathcal C_\Xi\bigr) \longrightarrow\mathrm E^s\mathrm{SLCY}_{\mathrm{slc}}\Xi.$
\end{center}

Let $f:((\mathcal X,\mathcal D),\mathcal C)\to S$ be a construction family. The sheaf $f_*\mathcal O_{\mathcal X}(r(K_{\mathcal X/S}+\mathcal D+\tau\mathcal C))$ is locally free. So it is trivial on the members of a finite open cover of $S$. On each of these open sets, we choose projective coordinates. This gives a local classifying morphism to the ambient strongly embedded space. The images of these morphisms under $G_N$ give the strongly embedded construction image.

Let $f:\bigl((\mathcal X,\mathcal D),\mathcal C\bigr)\to\mathcal S$ be a strongly embedded $\Xi$-family over a reduced scheme of finite type. We say that it is a family of stable log Calabi-Yau surfaces of type $\Upsilon$ if every geometric fibre satisfies the following two conditions: First, we can reverse Construction \ref{constr1} on a stratum with the fixed labelled Enriques diagram. In this situation, the extra member goes through the marked cluster. We get back $(Y_0,D_0,C_\infty,L,C_1,\Sigma,R)$ with $C_1\in|L|$ and $C_1|_{D_0}=\Sigma+R$ and
denote the recovered marked cluster by $\mathcal K_\Upsilon$. Here, We also require that its labelled Enriques diagram and the type of $C_1|_{D_0}$ are the data of Subsection \ref{sec42}. For the relation between these data, see Lemma \ref{lemclust}. Second, $\mathcal C_t$ lies in the given subsystem and also satisfies Definition \ref{defvarmem}. We denote this strongly embedded functor by $\mathcal P_\Upsilon$.

\subsection{The construction space attached to \texorpdfstring{$\Upsilon$}{Upsilon}}
\label{sec42}

We fix the weights, the proximity relations and the boundary component labels of the cluster. We let the positions of its proper centres vary. We take the strata on which the extra member $C_1$ goes sharply through $\mathcal K$ and has the fixed decomposition of $C_1|_{D_0}$. This decomposition includes the residual divisor $R$. We impose these conditions in Definition \ref{defstable}. The $\mathcal K$-system of Definition \ref{defksys} also requires $R=0$. We choose $C_1\in|L|$, and we stratify $C_1|_{D_0}$. Now we look at the log canonical modification morphism of $(Y_0,D_0+C_1)$. We take its weighted dual graph with the component labels, the intersection numbers and the multiplicities. We also take $C_1|_{D_0}$ with its points labelled by component. These two data determine the centres of the cluster. Here we mean the centres over the points where $C_1$ meets $D_0$ with multiplicity at least two. We also count the infinitely near centres. We fix the tuple
\begin{center}
$\Upsilon=\bigl(\Xi,q,a,[Y_0\to\Gamma,D_0,C_\infty,L, \mathfrak b_\Upsilon,\mathfrak N_\Upsilon]_{\mathrm{elt}}\bigr).$
\end{center}
Here $\mathfrak b_\Upsilon$ is the type of $C_1|_{D_0}$. It gives the multiplicities on each labelled component and also gives the number of points there. Notice that the datum $\mathfrak N_\Upsilon$ is the labelled Enriques diagram and we can get it back from the log canonical modification by Lemma \ref{lemclust}.

We write $[\ ]_{\mathrm{elt}}$ for the equivalence class under elementary transformations. Notice that, we have already compared two tuples on a common smooth resolution, as in Lemma \ref{lemlinear}. At there, we compared the boundary transforms, the pulled-back polarizing systems and any given subsystem.

Now, we recall toric models. Gross, Hacking and Keel define a toric model as a morphism to a smooth toric pair \cite[Definition 1.2]{GHK1}. We get the pair from the toric pair by blowing up smooth boundary points. Some of these points can be infinitely near. By \cite[Lemma 1.3]{GHK1}, every smooth rational pair with a nodal anticanonical cycle has a toric model after corner blow-ups. Alexeev, Arg\"uz and Bousseau use the same diagram for polarized pairs \cite[Definition 2.6]{AAB}. Let $(Y,D)$ be a smooth maximal log Calabi-Yau pair. Let $L$ be ample such that $C\in|L|$ defines a stable pair. We fix a toric model $(Y,D)\xleftarrow{g}(Y^\sharp,D^\sharp)\xrightarrow{h}(Y_0,D_0)$. Here $g$ is a sequence of corner blow-ups, and $h$ is a sequence of interior blow-ups.

\begin{lemm}
\label{lemtoric} From the boundary and the polarizing system, we get back the initial surface, the boundary, the line bundle and the marked cluster of this toric model, with $L_0=h_*g^*L$.
\end{lemm}
\begin{proof}
For a general pair, we take a toric model by \cite[Lemma 2.8]{AAB}. Its labelled toric pair is taut. By \cite[Lemma 2.15]{Fried}, it stays taut under corner blow-ups and blow-downs. Now the lemma follows from Lemma \ref{lemlinear}.
\end{proof}

\subsubsection{Initial surfaces and the universal curve}\label{sec421}

We choose a parameter space $\overline R_\Upsilon$ with universal family
\begin{center}
$\overline\rho_\Upsilon: (\overline{\mathcal Y}_\Upsilon,\overline{\mathcal D}_{\Upsilon,0}, \overline C_{\infty,\Upsilon},\overline{\mathcal L}_\Upsilon) \longrightarrow\overline R_\Upsilon$
\end{center}
for the tuples of initial surfaces, boundaries, sections and line bundles in the fixed class under elementary transformations.

For $q=0$ or $q=-\infty$ with reducible labelled boundary, we use the taut representatives of Definition \ref{deftuple}(i). Their classification and their behaviour under blow-downs are given in \cite[Lemma 2.15]{Fried}. We fix the self-intersection sequence and the class under elementary transformations. Then we take $\overline R_\Upsilon$ to be the finite disjoint union of these labelled representatives.

Now we consider an integral nodal boundary on $\mathbb F_e$ with $e\in\{1,2\}$. In this case, we have $\rank\Lambda(Y_0,D_0)=\rho(\mathbb F_e)-1=1$, and the pair has one-dimensional moduli. Let $U_e^{\mathrm{nod}}$ be the non-empty $\Aut(\mathbb F_e)$-invariant locally closed locus of integral anticanonical members with a unique node. We put $\overline R_\Upsilon=U_e^{\mathrm{nod}}$ and $\overline{\mathcal Y}_\Upsilon=\mathbb F_e\times U_e^{\mathrm{nod}}$. The universal boundary is $\overline{\mathcal D}_{\Upsilon,0}=\mathcal E_e^{\mathrm{nod}}$.

For a smooth anticanonical boundary and $q=0$, let $U_e\subset|-K_{\mathbb F_e}|$ be the non-empty invariant open locus of smooth members, again with $e\in\{1,2\}$. We put $\overline R_\Upsilon=U_e$ and $\overline{\mathcal Y}_\Upsilon=\mathbb F_e\times U_e$. We take the universal boundary $\mathcal E_e$, the section $C_\infty\times U_e$ and the line bundle $\pr_1^*L$, where $L=\mathcal O_{\mathbb F_e}(aC_\infty)$ and $a\geq2$. Thus we have $\overline\rho_{\Upsilon*}\overline{\mathcal L}_\Upsilon =H^0(\mathbb F_e,L)\otimes\mathcal O_{U_e}$, and the construction space below is $U_e\times|L|$. For the nodal case and for the smooth case, we use the same construction in Subsections \ref{sec42}--\ref{sec44}.

For $q=1$, we choose $N_0\geq3$, and we take the universal elliptic curve with full level structure $f:\mathcal E(N_0)\to Y(N_0)$. We put $\overline R_\Upsilon=\Pic^{-e}_{\mathcal E(N_0)/Y(N_0)}$. The curve is projective and flat with geometrically integral fibres. So its relative Picard scheme is separated and locally of finite type by \cite[Theorem 4.8]{Kle05}. The part of degree $-e$ is open, closed and proper over $Y(N_0)$. Moreover, a Poincar\'e bundle $\mathcal M_{\mathrm{univ}}$ exists by \cite[Exercise 4.3]{Kle05}. Over $[\mathcal M]\in\overline R_\Upsilon$, it represents $\mathcal M$. Rigidification at the zero section $O$ means that we choose an isomorphism
\begin{center}
$(O\times\id_{\overline R_\Upsilon})^* \mathcal M_{\mathrm{univ}}\simeq\mathcal O_{\overline R_\Upsilon}.$
\end{center}
Isomorphisms of these bundles have to preserve this trivialization. We have $f_*\mathcal O_{\mathcal E(N_0)}=\mathcal O_{Y(N_0)}$. Thus the rigidified universal bundle is unique up to unique isomorphism by \cite[Definition 2.8, Lemmas 2.9 and 2.10]{Kle05}. For this universal bundle $\mathcal M_{\mathrm{univ}}$, we put
\begin{center}
$\overline{\mathcal Y}_\Upsilon =\mathbb P_{\mathcal E(N_0)\times_{Y(N_0)}\overline R_\Upsilon} (\mathcal O\oplus\mathcal M_{\mathrm{univ}}).$
\end{center}
Its two universal sections are $\overline C_{0,\Upsilon}$ and $\overline C_{\infty,\Upsilon}$. We put $\overline{\mathcal D}_{\Upsilon,0}=\overline C_{0,\Upsilon} +\overline C_{\infty,\Upsilon}$ and $\overline{\mathcal L}_\Upsilon=\mathcal O(a\overline C_{\infty,\Upsilon})$.

We pass to strata where cohomology commutes with base change. Then $\overline\rho_{\Upsilon*}\overline{\mathcal L}_\Upsilon$ is locally free, and it commutes with arbitrary base change. Hence we can define
\begin{center}
$\mathbb P_\Upsilon =\mathbb P_{\overline R_\Upsilon} \bigl((\overline\rho_{\Upsilon*}\overline{\mathcal L}_\Upsilon)^\vee\bigr).$
\end{center}
By our convention in Subsection \ref{sec21}, its points are $(Y_0,D_0,C_\infty,L;C_1)$ with $C_1\in|L|$. Let $\pi_\Upsilon:\mathbb P_\Upsilon\to\overline R_\Upsilon$ be the projection. Let $\varpi_\Upsilon:\mathcal Y_{\mathbb P}\to\mathbb P_\Upsilon$ be the pulled-back family. We denote the pulled-back boundary, section and line bundle by $\mathcal D_{\Upsilon,0}$, $\mathcal C_{\infty,\Upsilon}$ and $\mathcal L_\Upsilon$. We pull back the tautological inclusion $\mathcal O_{\mathbb P_\Upsilon}(-1)\hookrightarrow \pi_\Upsilon^*\overline\rho_{\Upsilon*}\overline{\mathcal L}_\Upsilon$ to $\mathcal Y_{\mathbb P}$, and we compose it with the evaluation map. In this way, we get a section of $\mathcal L_\Upsilon\otimes\varpi_\Upsilon^*\mathcal O_{\mathbb P_\Upsilon}(1)$. Its zero scheme is the universal divisor $\mathcal C_1$. We also write $C_1$ for this divisor. Notice that the intersection $\mathcal Z=C_1\cap\mathcal D_{\Upsilon,0}$ is finite on the open locus where $C_1$ contains no boundary component.

\subsubsection{The stratum of \texorpdfstring{$C_1|_{D_0}$}{the restriction of C1 to D0}}
\label{sec422}

Let $\mathcal D_\lambda$ be the boundary components. So we have $\mathcal D_{\Upsilon,0}=\sum_{\lambda\in\Lambda}\mathcal D_\lambda$. We put $\mathcal Z_\lambda=\mathcal Z\cap\mathcal D_\lambda$. We write $C_1|_{D_0}=\Sigma^{\mathrm{mult}}+T+R$. The tuple $\mathfrak b_\Upsilon$ fixes the type of this decomposition. Let us explain the three parts. First, the support of $\Sigma^{\mathrm{mult}}$ consists of the points where $C_1$ meets $D_0$ with intersection multiplicity at least two. The log canonical modification gives us the centres of the cluster over these points. Second, the divisor $T$ consists of the transverse intersection points. We will blow up each of them. After finishing the process, the divisor $R$ is this remaining intersection of the strict transforms of $C_1$ and $D_0$.

On each component $\mathcal D_\lambda$, we fix $\Sigma_\lambda^{\mathrm{mult}}=\sum_{i=1}^{h_\lambda}m_{\lambda,i}p_{\lambda,i}$ with $m_{\lambda,i}\geq2$ and a reduced divisor $T_\lambda=q_{\lambda,1}+\cdots+q_{\lambda,t_\lambda}$. For the residual divisor $R_\lambda$, we also fix the number of its points and their multiplicities. These points are distinct smooth points of the labelled component. Notice that the supports of these divisors are pairwise disjoint and are in the smooth locus of $D_0$. In particular, they do not contain singular points of the boundary and of the total space.

Let $\Div_\lambda^{\mathrm{mult}}$, $\Div_\lambda^{\mathrm{tr}}$ and $\Div_\lambda^R$ be the relative divisor spaces over $\overline R_\Upsilon$ for the three divisors above. They parametrize the divisors with the given multiplicities at distinct smooth points, the reduced divisors of degree $t_\lambda$ and the residual divisors of the given type, respectively. Let $I_\Upsilon^{\mathrm{bd}}$ be the incidence correspondence of the tuples $(C_1;\Sigma^{\mathrm{mult}},T,R)$. We define it inside
\begin{center}
$\mathbb P_\Upsilon\times_{\overline R_\Upsilon} \prod_{\lambda\in\Lambda}\Div_\lambda^{\mathrm{mult}} \times\prod_{\lambda\in\Lambda}\Div_\lambda^{\mathrm{tr}} \times\prod_{\lambda\in\Lambda}\Div_\lambda^R$
\end{center}
by two conditions. The first one is the scheme-theoretic equality $C_1|_{D_0}=\Sigma^{\mathrm{mult}}+T+R$. The second one is the Abel--Jacobi--Picard condition $\mathcal O_{D_0}(\Sigma^{\mathrm{mult}}+T+R)\simeq L|_{D_0}$. The points of $I_\Upsilon^{\mathrm{bd}}$ are $(Y_0,D_0,C_\infty,L;C_1;\Sigma^{\mathrm{mult}},T,R)$.

The Abel map sends an effective Cartier divisor $Z$ in the smooth locus of $D_0$ to its line bundle:
\begin{center}
$\Ab_{D_0}:Z\longmapsto [\mathcal O_{D_0}(Z)]\in\Pic(D_0).$
\end{center}
Here $\Pic^0$ is the connected component of the identity in the Picard scheme. For a smooth elliptic curve $E$, it is $\Jac(E)$ \cite[Section 5]{Kle05}. For divisors of fixed multidegree, a reference divisor $Z_0$ gives $Z\mapsto[\mathcal O_{D_0}(Z-Z_0)]\in\Pic^0(D_0)$. Friedman calls this the Abel--Jacobi map for a boundary cycle \cite[Definition 1.12 and Lemma 1.13]{Fried}. In families, we use the relative Picard scheme \cite[Sections 2, 4 and 5]{Kle05}. Thus our Abel--Jacobi--Picard condition is $\Ab_{D_0}(\Sigma^{\mathrm{mult}}+T+R)=[L|_{D_0}]$. This is the fibre over the fixed line bundle. In other words, the class of $\mathcal O_{D_0}(\Sigma^{\mathrm{mult}}+T+R)\otimes(L|_{D_0})^{-1}$ is the identity in $\Pic^0(D_0)$. In the CHS case, we resolve $C_1|_{D_0}$, and we get $\mathcal O_X(C)|_W\simeq\mathcal O_W$, where $W$ is the boundary. This allows an algebraic contraction. For $R\neq0$, we get $\mathcal O_X(C)|_W\simeq\mathcal O_W(R)$ instead. The following conditions are all open: the supports are disjoint, the supports lie in the smooth locus, and $T$ is transverse.

Hacon and Jiang construct counterexamples to the boundedness of generalized log canonical surfaces \cite[Remark 4.3 and Lemma 4.4]{HaconJiang25}. In their constructions, they blow up and then contract an elliptic curve. This gives a family of algebraic spaces whose projective fibres occur on a countable dense subset. A torsion condition on a line bundle of degree zero on the elliptic curve controls the projectivity. Their Section 4.4 also uses blow-ups of elliptic ruled surfaces. The ruled surfaces themselves are projective. The projectivity fails after the contraction. We restrict the polarization to the elliptic curve. The projectivity of the contracted surface depends on the class of this restriction in the Picard group. For this reason, we fix the class of $\mathcal{O}_{D_0}(\Sigma^{\mathrm{mult}}+T+R)$ in $\Pic(D_0)$. This is the Abel--Jacobi--Picard condition above. Let $\rho:X\to S$ be a morphism which contracts $W$, and let $mC$ be the pullback of a Cartier divisor on $S$ for some integer $m>0$. Then we have
\begin{center}
$\mathcal O_X(mC)|_W\simeq\mathcal O_W.$
\end{center}
When $R=0$, the Abel--Jacobi--Picard condition gives $\mathcal O_X(C)|_W\simeq\mathcal O_W$ for the fixed polarization of the admissible tuple. This holds on $X$, before the contraction. On $\mathcal D_{\Upsilon,0}\times_{\mathbb P_\Upsilon}I_\Upsilon^{\mathrm{bd}}$, we have three universal divisors. We denote them by $\Sigma^{\mathrm{mult}}$, $\mathcal T$ and $\mathcal R$.
\begin{lemm}
\label{lembdinc} The space $I^{\mathrm{bd}}_\Upsilon$ is of finite type over $\mathbb C$, and the conditions $C_1|_{D_0}=\Sigma^{\mathrm{mult}}+T+R$ and $\mathcal O_{D_0}(\Sigma^{\mathrm{mult}}+T+R)\simeq L|_{D_0}$ are locally closed.
\end{lemm}
\begin{proof}
By \cite[Section 3]{Kle05}, the relative spaces of effective divisors of fixed degrees are of finite type. On these divisor spaces, the supports and the multiplicities have fixed types. The Fitting ideals of $(p_\lambda)_*\mathcal O_{\mathcal Z_\lambda}$ give a stratification of the finite intersection schemes by their length. On each stratum, the defining equality is a closed condition. The Abel--Jacobi--Picard condition is also closed, since it defines a fibre of the Abel map. The remaining conditions are open. So we get that $I_\Upsilon^{\mathrm{bd}}$ is locally closed and of finite type.
\end{proof}

\subsection{Simultaneous log canonical modification and the recovered cluster}
\label{sec43}

Over $I_\Upsilon^{\mathrm{bd}}$, the extra pair $(\mathcal Y_{\mathbb P},\mathcal D_{\Upsilon,0}+\mathcal C_1)$ determines the points of $C_1|_{D_0}$. The final boundary is the transform of $\mathcal D_{\Upsilon,0}$. The following version of \cite[Theorems 5.16 and 5.52]{Kol2} for surfaces gives the simultaneous modification. Let $f:(Y,\Delta)\to S$ be a projective flat family of pure relative dimension $2$ over a reduced seminormal scheme $S$ of finite type over $\mathbb C$. We assume the following conditions.
\begin{enumerate}
    \item Every geometric fibre $Y_t$ is normal for $t\in S$.
    \item The divisor $\Delta=\sum_{i=1}^r a_i\Delta_i$ is effective with $0\leq a_i\leq1$, and the $\Delta_i$ are relative Mumford divisors. In the applications below, they are relative Cartier divisors.
    \item There is a closed subset $Z\subset Y$ with the following properties. We have that $\Codim_{Y_t}(Z\cap Y_t)\geq2$ for every geometric point $t\in S$ and $\depth_Z\mathcal O_Y\geq2$. The family $f|_{Y\setminus Z}:(Y\setminus Z,\Delta|_{Y\setminus Z})\to S$ is locally stable.
    \item The divisor $K_{Y/S}+\Delta$ is $\mathbb Q$-Cartier.
\end{enumerate}
Let $H$ be an $f$-ample Cartier divisor. For every geometric point $t\in S$, we have $\nu_t:(Y_t^{\mathrm{lc}},\Theta_t) \to (Y_t,\Delta_t)$ be the log canonical modification for $(Y_t,\Delta_t)$. Here, $\Theta_t=(\nu_t)^{-1}_*\Delta_t+E_t$, and $E_t$ is the reduced $\nu_t$-exceptional divisor. Let
\begin{center}
    $A_t=\nu_t^*H_t$ and $L_t=K_{Y_t^{\mathrm{lc}}}+\Theta_t$.
\end{center} We define $\mathcal I(t)=I(A_t,L_t)=(A_t^2,A_t\cdot L_t,L_t^2)$ with the lexicographic order.

\begin{prop}
\label{propsimul} The tuple $\mathcal I(t)$ is constructible. It is lower semicontinuous for the lexicographic order. Take a locally closed locus on which $\mathcal I$ is constant. Then its seminormalization has a unique simultaneous log canonical modification. On the other hand, if such a modification exists, then $\mathcal I$ is locally constant.
\end{prop}

\begin{proof}
We divide the proof into several steps. First, we look at a single fibre. The fibres are normal, and the coefficients lie in $[0,1]$. Moreover, $K_{Y_t}+\Delta_t$ is $\mathbb Q$-Cartier. So by \cite[Definition 11.29 and Theorem 11.30]{Kol2}, the modification of each fibre exists. It is also unique.

\noindent\textbf{Step 1.} Let $W\subset S$ be an irreducible locally closed subset. We extend the generic modification $\nu_\eta:(Y_\eta^{\mathrm{lc}},\Theta_\eta)\to(Y_\eta,\Delta_\eta)$ to a projective birational morphism over a non-empty open subset $W^\circ$. After we shrink $W^\circ$, the source of this morphism is flat, and the pair is locally stable. Moreover, $K_{Y_{W^\circ}^{\mathrm{lc}}/W^\circ}+\Theta_{W^\circ}$ is ample over $Y_{W^\circ}$. Hence it is a simultaneous modification over $W^\circ$, and $\mathcal I$ is constant there. By Noetherian induction, we get a finite constructible stratification of $S$.

The first two entries are already locally constant on $S$. Indeed, we have $A_t^2=H_t^2$. The divisor $L_t-\nu_t^*(K_{Y_t}+\Delta_t)$ is exceptional. So the projection formula gives $A_t\cdot L_t=H_t\cdot(K_{Y_t}+\Delta_t)$. Both numbers are constant by the relative Cartier assumptions. So only $L_t^2$ is left.

\noindent\textbf{Step 2.} We apply the normal-fibre case of \cite[Theorem 5.52]{Kol2} to a smooth curve which maps to $S$, with generic point $\eta$ and closed point $0$. This gives $L_\eta^2\geq L_0^2$. Equality holds only when the pullback has a simultaneous modification. Thus we have $\mathcal I(\eta)\succeq\mathcal I(0)$, with the same criterion for equality. Now suppose that lower semicontinuity fails. By constructibility, we then get such a curve which joins a stratum on which $\mathcal I$ is constant to a closed point of its closure. This is a contradiction. Hence we get lower semicontinuity, as in \cite[5.30]{Kol2}.

\noindent\textbf{Step 3.} We construct the modification over a seminormal base on which $\mathcal I$ is constant. We work \'etale locally near a closed point $t_0\in S$. By uniqueness, the local constructions glue. The divisor $A_t$ is nef, and Step 1 gives finitely many strata. Hence we can choose one $m$ such that $L_t+mA_t$ is ample on every fibre. Here, we can take $m$ larger so that $mH$ is globally generated. Now, let $G\in|mH|$ be a general relative member and we can choose it so that the pair $(Y_{t_0}^{\mathrm{lc}},\Theta_{t_0}+\nu_{t_0}^*G_{t_0})$ is lc. Since $\mathcal I$ is constant, these pairs are stable by \cite[Theorem 5.52]{Kol2} and their volume must be equal to $L_t^2+2m(A_t\cdot L_t)+m^2A_t^2$. Now, we will rigidify these stable pairs. First, we choose finitely many smooth sections in $Y\setminus Z$. They are away from the divisors and their lifts rigidify the pairs. Then, we fix these sections and apply the argument of \cite[5.55--5.56]{Kol2} for stable pairs. In this argument, \cite[Theorem 5.52]{Kol2} gives the curve criterion. By this way, we get a reduced algebraic space $S'$, which is proper over $S$. It parametrizes the modifications and their graphs. The modification is unique on each fibre. So the morphism $S'\to S$ is quasi-finite and geometrically bijective. It is also proper. So it is finite and subintegral. Since $S$ is seminormal, we have $S'\simeq S$. From the universal family and its graph, we get $\nu:Y^{\mathrm{lc}}\to Y$ and the stable family $(Y^{\mathrm{lc}},\Theta+\nu^*G)\to S$. We subtract the Cartier divisor $\nu^*G$. Then $K_{Y^{\mathrm{lc}}/S}+\Theta$ is still relatively $\mathbb Q$-Cartier. Its fibres are lc. Thus $(Y^{\mathrm{lc}},\Theta)\to S$ is locally stable with fibres $(Y_t^{\mathrm{lc}},\Theta_t)$.

Conversely, take a simultaneous modification. Then $A=\nu^*H$ and $L=K_{Y^{\mathrm{lc}}/S}+\Theta$ are relative $\mathbb Q$-Cartier classes on a flat family. So their intersection numbers $A_t^2,A_t\cdot L_t,L_t^2$ are locally constant. This proves the converse.
\end{proof}

We apply Proposition \ref{propsimul} to $(\mathcal Y_{\mathbb P},\mathcal D_{\Upsilon,0}+\mathcal C_1)$ over $I_\Upsilon^{\mathrm{bd}}$. We take a locally closed locus on which $\mathcal I$ is constant, and we denote its seminormalization by $I_\Upsilon^{\mathrm{lc}}$. On this space, we get
\begin{center}
$\nu_\Upsilon^{\mathrm{lc}}: (\mathcal Y_\Upsilon^{\mathrm{lc}}, \mathcal D_\Upsilon^{\mathrm{lc}}+\mathcal C_\Upsilon^{\mathrm{lc}}+\mathcal E_\Upsilon) \longrightarrow (\mathcal Y_{\mathbb P},\mathcal D_{\Upsilon,0}+\mathcal C_1)|_{I_\Upsilon^{\mathrm{lc}}}.$
\end{center}
Here $\mathcal D_\Upsilon^{\mathrm{lc}} =(\nu_\Upsilon^{\mathrm{lc}})^{-1}_*\mathcal D_{\Upsilon,0}$, $\mathcal C_\Upsilon^{\mathrm{lc}} =(\nu_\Upsilon^{\mathrm{lc}})^{-1}_*\mathcal C_1$, and $\mathcal E_\Upsilon$ is the reduced exceptional divisor. We apply Construction \ref{constr1} relatively on the strata of ordered clusters of Proposition \ref{propclust}. Here we include the proper transverse intersection points.

\subsection{The stable construction locus}
\label{sec44} We stratify $I^{\mathrm{lc}}_\Upsilon$, and we take a finite \'etale cover over which the components of the reduced exceptional divisor are defined. Then we can write this divisor as
\begin{center}
$\mathcal E_\Upsilon=\sum_{i=1}^N\mathcal E_i.$
\end{center}
Here each $\mathcal E_i$ is a relative prime exceptional divisor. For a geometric point $b$, we consider the weighted dual graph of the exceptional divisor of the simultaneous log canonical modification. The weight of a vertex $E_i(b)$ consists of its self-intersection number together with its intersection numbers with $D^{\mathrm{lc}}_\lambda(b)$ and $C^{\mathrm{lc}}(b)$. Here, We label each point of $C_{1,b}|_{D_{\Upsilon,0,b}}$ by the component which contains it. By Lemma \ref{lemclust}, this graph and this labelled divisor determine the labelled Enriques diagram $\mathfrak N(b)$ of the cluster of $C_{1,b}\cap D_{\Upsilon,0,b}$ and the singular types $A_n$ of the contracted chains of $(-2)$-curves. Let $I^{\mathrm{dg}}_\Upsilon$ be the locus in $I^{\mathrm{lc}}_\Upsilon$ where $\mathfrak N(b)\simeq\mathfrak N_\Upsilon$.

\begin{lemm}
\label{lemgraph} The locus $I^{\mathrm{dg}}_\Upsilon$ of tuples with the given weighted dual graph and boundary multiplicities from $\Upsilon$ is locally closed in $I^{\mathrm{lc}}_\Upsilon$.
\end{lemm}

\begin{proof}
By Lemma \ref{lemclust}, $I^{\mathrm{dg}}_\Upsilon$ is the locus where the cluster of $C_1\cap D_0$, with points labelled by component, has diagram $\mathfrak N_\Upsilon$. This locus is locally closed by Proposition \ref{propclust}.
\end{proof}

The log canonical modification can contract a $(-2)$-curve which meets the strict transform of $C_1$, if such a curve exists. Such a curve survives on the stable model. For this reason, we first blow up every centre of $\mathcal K$. Only then we contract the curves which are disjoint from the strict transform of $C_1$. On the stratum of ordered centres with the given labels, weights and proximities, Proposition \ref{propclust} gives the universal centres of the blow-up of $\mathcal K$. We blow up the centres which resolve $\Sigma^{\mathrm{mult}}$. This gives $\nu^{\mathrm{mult}}_\Upsilon\colon\mathcal Y^{\mathrm{mult}}_\Upsilon \to\mathcal Y_{\mathbb P}|_{I^{\mathrm{dg}}_\Upsilon}$. The remaining transverse divisor $\mathcal T^{\mathrm{mult}}$ avoids the centres which resolve $\Sigma^{\mathrm{mult}}$. So it agrees with $\mathcal T^{\mathrm{lc}}$ near its support. We take the locus where this divisor is finite \'etale over the base and contained in the smooth locus of the boundary. On this locus, its blow-up commutes with base change. We put
\begin{center}
$ \beta_\Upsilon\colon \widetilde{\mathcal Y}_\Upsilon :=\Bl_{\mathcal T^{\mathrm{mult}}} \mathcal Y^{\mathrm{mult}}_\Upsilon \longrightarrow\mathcal Y^{\mathrm{mult}}_\Upsilon. $
\end{center}
When $T=0$, the morphism $\beta_\Upsilon$ is the identity. We put
\begin{center}
$\mu_\Upsilon:=\nu^{\mathrm{mult}}_\Upsilon\circ\beta_\Upsilon\colon \widetilde{\mathcal Y}_\Upsilon \longrightarrow\mathcal Y_{\mathbb P}|_{I^{\mathrm{dg}}_\Upsilon}. $
\end{center}
We denote by $\widetilde{\mathcal D}_\Upsilon$ and $\widetilde{\mathcal C}_\Upsilon$ the strict transforms of $\mathcal D_{\Upsilon,0}$ and $\mathcal C_1$. The centres of $\Sigma^{\mathrm{mult}}$ and $T$ lie on the boundary. After we resolve them, we have
\begin{center}
$\widetilde{\mathcal C}_\Upsilon|_{\widetilde{\mathcal D}_\Upsilon} =\widetilde{\mathcal R}. $
\end{center}
The divisor $\mathcal R$ does not meet the resolved centres. So $\widetilde{\mathcal R}$ is both the strict transform and the pullback of $\mathcal R$.

\begin{defi}\label{defstable}
Let $I^{\mathrm{st}}_\Upsilon$ be a finite disjoint union of reduced locally closed strata of $I^{\mathrm{dg}}_\Upsilon$ on which the following conditions hold.
\begin{enumerate}
    \item We look at the exceptional centres of the simultaneous log canonical modification. Each of them lies on the strict transform of the labelled anticanonical boundary at the step where it is blown up.
    \item The divisor $\mathcal T^{\mathrm{mult}}$, is finite \'etale over the base and is contained in the smooth locus of the boundary. In particular, its blow-up commutes with base change.
    \item The equality $\widetilde{\mathcal C}_\Upsilon|_{\widetilde{\mathcal D}_\Upsilon} =\widetilde{\mathcal R}$ holds on every geometric fibre.
    \item $\widetilde{\mathcal C}_\Upsilon$ is nef and big on every geometric fibre, and $(\widetilde Y_b,\widetilde D_b+\tau\widetilde C_b)$ is log canonical for a small $\tau>0$.
    \item For $0<t\ll1$, the relative section ring of
    \begin{center}
    $ K_{\widetilde{\mathcal Y}_\Upsilon/I^{\mathrm{st}}_\Upsilon} +\widetilde{\mathcal D}_\Upsilon+t\widetilde{\mathcal C}_\Upsilon $
    \end{center}
    is finitely generated, and its relative Proj is compatible with geometric fibres. On each fibre, this divisor is $\mathbb Q$-linearly equivalent to $t\widetilde C$.
    \item The pair we get is log canonical, and its strongly embedded tuple of invariants is $\Xi$.
\end{enumerate}
\end{defi}

By conditions (1) and (2), all centres of $\mu_\Upsilon$ lie on the boundary. The blow-up formula gives
\begin{center}
$ K_{\widetilde{\mathcal Y}_\Upsilon/I^{\mathrm{st}}_\Upsilon} +\widetilde{\mathcal D}_\Upsilon =\mu_\Upsilon^*\bigl( K_{\mathcal Y_{\mathbb P}/I^{\mathrm{st}}_\Upsilon} +\mathcal D_{\Upsilon,0}\bigr). $
\end{center}
The initial boundary on every geometric fibre $b$ is anticanonical. So we have $K_{\widetilde{\mathcal Y}_{\Upsilon,b}} +\widetilde{\mathcal D}_{\Upsilon,b}\sim0$.

We denote by $h_\Upsilon\colon\widetilde{\mathcal Y}_\Upsilon\to I^{\mathrm{st}}_\Upsilon$ the induced family. For sufficiently divisible $m>0$, we put
\begin{center}
$ \mathcal X_\Upsilon :=\Proj_{I^{\mathrm{st}}_\Upsilon} \bigoplus_{\ell\geq0}h_{\Upsilon*} \mathcal O_{\widetilde{\mathcal Y}_\Upsilon} (\ell m\widetilde{\mathcal C}_\Upsilon). $
\end{center}
Let $\rho_\Upsilon\colon\widetilde{\mathcal Y}_\Upsilon\to \mathcal X_\Upsilon$ be the contraction. We put $\mathcal D_\Upsilon=(\rho_\Upsilon)_*\widetilde{\mathcal D}_\Upsilon$ and $\mathcal C_\Upsilon=(\rho_\Upsilon)_*\widetilde{\mathcal C}_\Upsilon$. This gives the family
\begin{center}
$ ((\mathcal X_\Upsilon,\mathcal D_\Upsilon),\mathcal C_\Upsilon) \longrightarrow I^{\mathrm{st}}_\Upsilon. $
\end{center}

Let $b\in I^{\mathrm{lc}}_\Upsilon$ be a geometric point, and let
\begin{center}
$ \nu_b\colon (Y^{\mathrm{lc}}_b,D^{\mathrm{lc}}_b+C^{\mathrm{lc}}_b+E_b) \longrightarrow(Y_b,D_{\Upsilon,0,b}+C_{1,b}) $
\end{center}
be its log canonical modification.
\begin{lemm}
\label{lemclust} We take the modification morphism $\nu_b$. We also take the decomposition $C_{1,b}|_{D_{\Upsilon,0,b}}=\Sigma^{\mathrm{mult}}(b)+T(b)+R(b)$. Here the points are labelled by boundary component. These data determine the marked cluster of Definition \ref{defmark}, its labelled Enriques diagram and its residual divisor. The weighted dual graph with its component labels and intersection numbers gives the multiplicities, the proximities and the contracted curves. The divisor $C_{1,b}|_{D_{\Upsilon,0,b}}$ with its points labelled by component gives the positions of the centres and $R(b)$. Conversely, the initial surface, the marked cluster and the residual divisor give the stable family.
\end{lemm}

We can make the conditions in this lemma explicit as follows:
\begin{enumerate}
    \item Under condition (1) of Definition \ref{defstable}, the centres over $p=p_{\lambda,i}$ in $\Sigma^{\mathrm{mult}}_\lambda(b)$ form a sequence $p=p^{(1)},\ldots,p^{(k_p)}$. Here, each $p^{(j+1)}$ is the intersection point between the boundary with the exceptional curve of the blow-up at $p^{(j)}$. In particular, we have $p^{(j+1)}\succ p^{(j)}$. Away from $D_{\Upsilon,0,b}$, the curve $C_{1,b}$ has at worst ordinary nodes.
    \item Under condition (1), the reduced points of $T(b)$ and the multiplicities $m_p^{(j)}=\mult_{p^{(j)}}(C_{1,b})$ of Definition \ref{defmult} form a marked anticanonical cluster $\mathcal K(b)$ on $(Y_b,D_{\Upsilon,0,b})$, as defined in Section \ref{sec3}. The curve $C_{1,b}$ contains this cluster with the given multiplicities. Its intersection multiplicity with $D_{\Upsilon,0,b}$ at $p$ is $\sum_{j=1}^{k_p}m_p^{(j)}=m_{\lambda,i}$. Thus on the blow-up $\mu_{\mathcal K(b)}$ of the cluster, the strict transforms of $C_{1,b}$ and $D_{\Upsilon,0,b}$ meet exactly along the preimage of $R(b)$. This preimage is isomorphic to $R(b)$. For $R=0$, this is a $\mathcal K(b)$-system.
    \item The diagram $\mathfrak N(b)$ is the Enriques diagram of $\mathcal K(b)$. We add to it the labels and the types $A_n$ of the contracted chains of $(-2)$-curves along $E_b$. The relative log canonical section ring defines a morphism which contracts these chains to $A_n$ points. The weighted dual graph and its intersection numbers determine these points. If every $m_p^{(j)}$ is equal to one, then we can have the maximal chains of consecutive $(-2)$-curves. On the other hand, we get a stable log Calabi-Yau surface from a general member of a $\mathcal K_\Upsilon$-system. 
\end{enumerate}
\begin{proof}
By condition (1), every centre lies on the boundary at the step where we blow it up. Away from the boundary, we can see the log canonicity implies that the curve has at worst ordinary nodes there. Hence, at each step, we focus on the centre lying the boundary. This gives us the equality $\sum_jm_p^{(j)}=m_{\lambda,i}$ and the remaining restriction $R(b)$. Now Equation \eqref{eqstrict} gives us back the weighted Enriques diagram. By adjunction, we get the self-intersection numbers of the exceptional curves. Hence, we get the contracted chains of $(-2)$-curves and their $A_n$ labels. Finally, by using Construction \ref{constr1} and conditions (4)--(6) of Definition \ref{defstable}, we get that the two constructions are inverse to each other and this finish the proof.
\end{proof}
\begin{prop}
\label{propstloc} The space $I^{\mathrm{st}}_\Upsilon$ is a finite disjoint union of locally closed construction strata. The geometric fibres of its universal family are exactly the stable log Calabi-Yau surfaces which we get from $\Upsilon$.
\end{prop}
\begin{proof}
By Lemma \ref{lembdinc}, the space $I^{\mathrm{bd}}_\Upsilon$ is locally closed and of finite type. On the strata of Proposition \ref{propsimul}, Lemma \ref{lemgraph} gives the locally closed locus $I^{\mathrm{dg}}_\Upsilon$. Over this locus, we have a simultaneous log canonical modification with the given graph. Now we check the conditions of Definition \ref{defstable} one by one. Condition (1) only uses the labelled graph and its intersections with the labelled boundary components.

In condition (2), we ask for reducedness and finite \'etaleness. These are open conditions. The same holds for the containment in the relative log smooth locus. Under these conditions, the blow-up commutes with base change. When $\mathcal R$ is fixed, condition (3) is closed.

For condition (4), we first notice that log canonicity is an open condition. Since the divisor $\widetilde{\mathcal C}_\Upsilon$ is effective, it must be nef on a fibre if it has nonnegative degree on each of its own components.

Fix the Cartier index. We then have these components have bounded degree and finitely many Hilbert polynomials. Now we take the relative Hilbert schemes of $\widetilde{\mathcal C}_\Upsilon$ for the curves of negative degree. We have that their images are proper. Because any limiting cycle has a component of negative degree, it implies that the non-nef locus must be closed. Since $\widetilde C_b^2$ is constant and positive, every nef fibre is big.

For (5), we refine to smooth locally closed strata on which the resolved pair is locally stable. Since the restriction of its log canonical divisor to each fibre is $\tau\widetilde C_b$, $|C_b|$ gives a contraction to a stable pair on each fibre. By \cite[Theorem 2.22]{FI21}, the simultaneous contractions form a locally closed partial decomposition of the base. By Noetherian induction, we get finitely many strata which carry these contractions, with base change on geometric fibres.

For (6), we take the strata with invariants $\Xi$. Notice that the divisor $K_{\widetilde Y_b}+\widetilde D_b$ is principal. Hence its push-forward gives a Cartier divisor $K_{X_b}+D_b\sim0$. Thus, the fibres are the stable log Calabi-Yau surfaces of Construction \ref{constr1}.
\end{proof}

\begin{prop}
\label{propfam} Every fixed tuple $\Upsilon$ has a family of finite type. Its geometric fibres are exactly the stable log Calabi-Yau surfaces constructed from $\Upsilon$.
\end{prop}

\begin{proof}
This is the universal family of Proposition \ref{propstloc}.
\end{proof}

We put $B=I^{\mathrm{st}}_\Upsilon$. We denote this family by $p_\Upsilon\colon\mathcal X_\Upsilon\to B$, and we put $\mathcal A_\Upsilon=\mathcal O_{\mathcal X_\Upsilon}(\mathcal C_\Upsilon)$.

\begin{lemm}
\label{lemdesc} Let $\widetilde{\mathcal C}_\Upsilon$ be relatively Cartier with no contracted irreducible component on any geometric fibre. Then its line bundle descends to the relatively ample invertible sheaf $\mathcal A_\Upsilon$ on $\mathcal X_\Upsilon$.
\end{lemm}

\begin{proof}
We denote the contraction by $\rho=\rho_\Upsilon$ and the polarizing line bundle by $\mathcal M=\mathcal O(\widetilde{\mathcal C}_\Upsilon)$. By condition (5), the contraction contracts exactly the fibre curves which do not meet $\widetilde{\mathcal C}_\Upsilon$. These curves are not irreducible components of the effective divisor $\widetilde{\mathcal C}_{\Upsilon,b}$ on their geometric fibre, and they are disjoint from this divisor. Let $Z\subset\mathcal X_\Upsilon$ be the locus of positive-dimensional fibres of $\rho$. The open set
\begin{center}
$U=\mathcal X_\Upsilon\setminus \rho(\Supp\widetilde{\mathcal C}_\Upsilon)$
\end{center}
contains $Z$. The section which defines $\widetilde{\mathcal C}_\Upsilon$ trivializes $\mathcal M$ on $\rho^{-1}U$. Over $\mathcal X_\Upsilon\setminus Z$, the morphism $\rho$ is an isomorphism. These two open sets cover the target. Since $\rho_*\mathcal O=\mathcal O$, the transition functions on their inverse images come from transition functions on the target. Thus $\mathcal O_U$ and $(\rho_*\mathcal M)|_{\mathcal X_\Upsilon\setminus Z}$ glue to an invertible sheaf. We get $\mathcal A_\Upsilon=\rho_*\mathcal M =\mathcal O_{\mathcal X_\Upsilon}(\mathcal C_\Upsilon) $ and $ \rho^*\mathcal A_\Upsilon\simeq\mathcal M.$ The relative Proj makes a positive tensor power of $\mathcal A_\Upsilon$ relatively ample. So $\mathcal A_\Upsilon$ is relatively ample.
\end{proof}

After we blow up all given cluster centres, $\widetilde{\mathcal C}_\Upsilon$ is relatively Cartier. Hence for the linear-system construction, we restrict to the strata of Lemma \ref{lemdesc}. There its line bundle descends to $\mathcal A_\Upsilon$. On a stratum of $B$ where cohomology commutes with base change, $\mathcal V_\Upsilon=p_{\Upsilon*}\mathcal A_\Upsilon$ is locally free and compatible with arbitrary base change. For a given proper movable subsystem, we take its locally free subbundle instead. Let $\pi_{\Lin}\colon \Lin_\Upsilon=\mathbb P_B(\mathcal V_\Upsilon^\vee)\to B$ and $\mathcal X_{\Lin}= \mathcal X_\Upsilon\times_B\Lin_\Upsilon$. We compose the tautological inclusion $\mathcal O_{\Lin_\Upsilon}(-1) \hookrightarrow\pi_{\Lin}^*\mathcal V_\Upsilon$ with the evaluation map. In this way, we get the universal divisor $\mathcal C^{\mathrm{var}}$ on this family. Its points are the triples $((X_b,D_b),C_b^{\mathrm{aux}},C_b^{\mathrm{var}})$ with both members in $|\mathcal A_b|$. The extra member $C_b^{\mathrm{aux}}=\mathcal C_{\Upsilon,b}$ is smooth and general. It determines the cluster and the contraction defined by its linear system. The variable member is the polarization of the stable pair. It has to satisfy the conditions of the following definition.
\begin{defi}
\label{defvarmem} Let $\Lin^{\mathrm{st}}_\Upsilon \subset\Lin_\Upsilon$ be the locus of geometric points $(b,C_b^{\mathrm{var}})$ satisfying the following conditions.
\begin{enumerate}
    \item $C_b^{\mathrm{var}}$ contains no component of $D_b$, and $(X_b,D_b+\tau C_b^{\mathrm{var}})$ is log canonical, where $\tau$ is determined by $\Xi$.
    \item For every labelled component $D_{b,\lambda}$ which is not contracted by $\rho_{\Upsilon,b}$, the restriction $C_b^{\mathrm{var}}|_{D_{b,\lambda}}$ lies in the image of $\Div^R_\lambda$ under $\rho_{\Upsilon,b}$. Thus the type of $R$ on the boundary is the one fixed by $\Upsilon$.
\end{enumerate}
\end{defi}

The stable log Calabi-Yau family $((\mathcal X_{\Lin},\mathcal D_{\Lin}), \mathcal C^{\mathrm{var}})\to\Lin^{\mathrm{st}}_\Upsilon$ has tuple of invariants $\Xi$. We choose projective coordinates on a finite open cover $\{U_\alpha\}$ of $\Lin^{\mathrm{st}}_\Upsilon$, and we forget the extra member. Then we have the classifying morphisms
\begin{center}
$ \widehat\Phi_{\Upsilon,\alpha}\colon U_\alpha \longrightarrow\mathrm E^s\mathrm{SLCY}_{\mathrm{lc}}\Xi. $
\end{center}
The extra member defines a section $\delta_\Upsilon$ over $I^{\mathrm{st}}_\Upsilon$ by $\mathcal C^{\mathrm{var}}=\mathcal C^{\mathrm{aux}}$. The restrictions of the local classifying morphisms along $\delta_\Upsilon$ are the morphisms of Proposition \ref{propfam}.

\begin{rema}
For $R=0$, we have $\widetilde{\mathcal C}_\Upsilon\cap\widetilde{\mathcal D}_\Upsilon =\varnothing$.
Now let $R\neq0$. Then the contracted surface is a stable log Calabi-Yau minimal model with non-trivial boundary and the polarization divisor meets the boundary along $R$.
\end{rema}

Now we come to the length-two clusters of Subsection \ref{sec34}. Fix $a\geq2$, $e\geq1$, $q\leq1$ and $0\leq w\leq e$. We assume that $\Upsilon$ has the boundary decomposition $C_1|_{D_0}=\Sigma^{\mathrm{mult}}+R$ with $\Sigma^{\mathrm{mult}}=2p_1+\cdots+2p_w$ and $\deg R=(a-2)e+(e-w)$. Here, we blow up twice over each $p_i$ to resolve the clusters, as in Subsection \ref{sec34}. Both multiplicities are equal to one. The other centres are reduce points and we blow up each of them once.

\begin{lemm}\label{lemm2pt}
Each length-two boundary point resolves to a generalized exceptional curve with two components. Their self-intersections are $-2$ and $-1$. The boundary and the polarization meet only the last curve. This is the $(-1)$-curve, and they meet it transversely. The log canonical modification contracts the $(-2)$-curve to an $A_1$ point.
\end{lemm}
\begin{proof}
We apply Lemma \ref{lemclust} with $m^{(1)}=m^{(2)}=1$. By Equation \eqref{eqstrict}, the intersection numbers on $\widetilde E_1$ and $\widetilde E_2$ are zero and one. Then we contract the curve with intersection zero to the polarize divisors, which then gives us $A_1$ points.
\end{proof}

\begin{prop}\label{propresol}
The blow-up of the cluster is smooth and compatible with base change. It has $w$ disjoint internal $(-2)$-curves  which are disjoint from the polarizing divisor and their contraction gives $A_1$ points on the surface of the log canonical modification.
\end{prop}
\begin{proof}
Notice that the first centre $\Supp\Sigma^{\mathrm{mult}}$ is finite \'etale in the smooth boundary locus. Also, we can see that the second centre is the intersection of the first exceptional divisor with the boundary as it is a section over the first centre. Hence, both centres are finite \'etale over the base and lying in the smooth locus of the boundary. Thus, blowing them up commutes with base change.

From all the steps above, the remaining transverse centres are the centres of $\beta_\Upsilon$. Blowing up all these gives us a smooth family which commutes with base change. Here, the fibre of the family is $\mu_{\mathcal K(b)}$ of Construction \ref{constr1}. On each geometric fibre, Construction \ref{constr1} gives $K_{\widetilde Y_b}+\widetilde D_b\sim0$. By Lemma \ref{lemm2pt}, we get the $(-2)$-curves and their contractions.
\end{proof}

\subsection{The strongly embedded locus}
\label{sec45} We use these classifying morphisms. We define the subset
\begin{center}
$ \mathrm E^s\mathrm{SLCY}_{\mathrm{lc}}\Upsilon :=\bigcup_\alpha G_N\cdot\widehat\Phi_{\Upsilon,\alpha}(U_\alpha) \subset\mathrm E^s\mathrm{SLCY}_{\mathrm{lc}}\Xi.$
\end{center}
We denote it by $J_\Upsilon$. When $J_\Upsilon$ is locally closed, we give it the reduced locally closed structure.
\begin{prop}
\label{propimage} The strongly embedded image $J_\Upsilon$ is constructible. It has a finite partition into reduced locally closed subspaces. Moreover, its points are the strongly embedded surfaces which we construct from $\Upsilon$.
\end{prop}
\begin{proof}
For each $\alpha$, we consider the map $G_N\times U_\alpha\to\mathrm E^s\mathrm{SLCY}_{\mathrm{lc}}\Xi$ given by $(g,u)\mapsto g\cdot\widehat\Phi_{\Upsilon,\alpha}(u)$. These maps are of finite type, and there are finitely many of them. Hence the union of their images is constructible by Chevalley's theorem . Also by the construction and by Definition \ref{defvarmem}, we have that its points are the surfaces we want. By the fact that a constructible subset of a Noetherian space has a finite partition into locally closed subsets, we get the partition we aim for.
\end{proof}
In the next lemma, we fix a tuple in the elliptic ruled construction. We represent it with $0<m_\alpha<a$.
\begin{lemm}
\label{lemcanon} The intersection numbers for the canonical modification are constant. Near its image, the universal family over its reduced closure has a simultaneous canonical modification after seminormalization. After a finite local analytic base change, it also has a smooth projective resolution.
\end{lemm}
\begin{proof}
Let $V$ be the blow-up of the ordered cluster. Let $\rho:V\to X$ be the contraction of Construction \ref{constr1}, and let $C$ be the polarization on $V$. The morphism $\rho$ contracts exactly the curves which do not meet $C$. Since the base curve is elliptic, we have that every rational curve on $V$ is exceptional, or it is the strict transform of a fibre of the ruling. Since all the centres must lie on the section with positive self-intersect, we must have that: all the $C$-trivial curves which do not meet the boundary are $(-2)$-curves. Here, $\rho$ contracts them to rational double points. On the other hand, the $C$-trivial curves in the boundary are elliptic. Their images are simple elliptic singularities. The divisor $K_X+D$ is Cartier, and $(X,D)$ is log canonical. Hence, the canonical modification $f:(Z,D_Z)\to(X,D)$ extracts exactly these elliptic curves, with discrepancy $-1$. Thus, we have $K_Z+D_Z\sim-E$ with $E$ is the reduced sum of extracted elliptic curves. For the ambient hyperplane class $H_X$, we get
\begin{center}
$I(f^*H_X,K_Z+D_Z)=(H_X^2,0,E^2).$
\end{center}
The $(-2)$-curves do not meet $E$. Hence the number $E^2$ is the sum of the self-intersection numbers on $V$ of the contracted elliptic boundary components. This number depends only on $\Upsilon$.

Let $B$ be the reduced closure of the image of the parameter space $I^{\mathrm{st}}$. By \cite[Theorem 5.16]{Kol2}, the locus of $B$ where $I$ takes this value is open. It contains this image. Thus over the seminormalization of this locus, the family has a simultaneous canonical modification. The fibres of this modification have only rational double points. By the theorem of Brieskorn \cite{Brieskorn71} and Tjurina \cite{Tjurina70}, the versal deformation of a rational double point has a simultaneous resolution after a finite Galois base change. The Galois group of this base change is the Weyl group, see \cite[Section 0 and Theorem 1.3]{Wahl79}. Local analytically, our family is the pullback of the versal deformations of the singular points of a fixed fibre. Hence it has a simultaneous resolution after a finite local analytic base change. Moreover, the resolved fibres are ruled surfaces, so that $H^2(\mathcal O)=0$ and this implies that an ample line bundle on a resolved fibre extends locally in the relative Picard space. Since ampleness is an open condition, the resolved family is locally projective.
\end{proof}
In the next lemma, we first refine the strata into finitely many invariant locally closed strata.
\begin{lemm}
\label{lemrevMMP} On these refined strata, the universal family is, \'etale locally, the target of the following contraction. We take a family of initial surfaces, we blow up the ordered cluster, and we contract by the polarization. Every fixed finite sequence of smooth $(-1)$-contractions extends locally.
\end{lemm}
\begin{proof}
Let $P$ be the reduced disjoint union of the schemes $G_N\times U_\alpha$. The scheme $P$ has a $G_N$-equivariant surjection onto the strongly embedded image. By generic smoothness in characteristic zero and by Noetherian induction, the image has a finite partition into invariant locally closed strata with the following property. An open part of $P$ maps smoothly and surjectively onto each stratum. A smooth surjection has sections \'etale locally on its base by \cite[Tag~055U]{Stacks}. So we can pull back the initial surfaces and the marked clusters along such sections. We also pull back the blow-ups and the contractions. The classifying maps agree. So this gives the universal family. The pulled-back blow-ups of the clusters form a family of smooth surfaces. On each fibre, this family is the minimal resolution of the singularities.

Let $E$ be a $(-1)$-curve on a fibre $Y$. We have $N_{E/Y}\simeq\mathcal O_{\mathbb P^1}(-1)$. Hence we get $H^0(N_{E/Y})=H^1(N_{E/Y})=0$. This implies that the relative Hilbert scheme is \'etale at $[E]$. So we make an \'etale base change, and we shrink the base. Then $E$ extends uniquely to a relative $(-1)$-curve. Moreover, the relative Castelnuovo contraction commutes with base change, and its target is smooth. We repeat this argument again for others $(-1)$ exceptional curves. By this way, we can extend the chosen sequence of contractions algebraically. Finally, by Lemma \ref{lemlinear}, the polarized linear spaces are the same. Thus, we get back the given polarizing system fiberwise and finish the proof.
\end{proof}
In the next proposition, we consider the decomposable elliptic ruled construction.
\begin{prop}
\label{propelloc} The strongly embedded image of the tuple $\Upsilon$ is locally closed.
\end{prop}
\begin{proof}
We will use the elementary transformations of Subsection \ref{sec42}. By Lemma \ref{lemlinear}, the sections do not change under blow-ups of weight zero points. Hence, we can take a representative with $0<m_\alpha<a$. By Lemma \ref{lemcanon}, we can make a finite local analytic base change to get a smooth projective resolution near the image point. By Lemma \ref{lemrevMMP}, the sequence of $(-1)$-contractions to the ruled model and the polarizing classes extend over a neighbourhood. 

Now take a resolved fibre $V$. The projection of the initial ruled surface induces a morphism $\alpha:V\to\Gamma$. On the dense image, the morphism $\alpha$ and the two labelled elliptic boundaries determine the ruled contraction and its cluster up to $[\ ]_{\mathrm{elt}}$. Hence, the boundary labels and the proximity relations agree with the given ones. These conditions are closed. On the other hand, there are other two conditions are open. The first one is that the effective multiplicities are equal to the weights. The second one is that there are no further centres. We have two disjoint sections. So we get a decomposable elliptic ruled family. By restricting it to the negative section, we get $L\simeq\mathcal O(aC_\infty)$. The boundary restriction is surjective by \eqref{eqpush}. So we can lift $C_1|_{D_0}$ as the elliptic curve varies. The map $(p_i)\mapsto\mathcal O_\Gamma(\sum_i m_i p_i)$ is smooth. Indeed, multiplication by each positive integer $m_i$ has the nonzero differential on the elliptic curve over $\mathbb C$. Thus, the recovered data gives us an open neighbourhood inside the image. Since the Openness passes through the finite cover, we must have that the constructible image is open in a closed subset. This means that it is locally closed.
\end{proof}

\begin{defi}\label{defstrata}
We put $J_\Upsilon:=\mathrm E^s\mathrm{SLCY}_{\mathrm{lc}}\Upsilon$, and we choose a finite invariant partition
\begin{center}
$J_\Upsilon=\coprod_\nu J_{\Upsilon,\nu},\
\mathrm E^s\mathrm{SLCY}_{\mathrm{lc}}(\Upsilon,\nu):=J_{\Upsilon,\nu}$
\end{center}
into reduced locally closed subschemes. Such partition exists by Proposition \ref{propimage}. Here, we get such invariant pieces, labelled by $\nu$, in the following way: First, we take the interior in the closure and then we repeat this on the closed complement. For $q=1$, we omit $\nu$ and let $J_{\Upsilon,0}=J_\Upsilon$. In the rational cases, we keep the finite partition.
\end{defi}

We define the stratum functors on reduced bases of finite type. For such a base, we ask that its classifying map factors through the chosen piece.

\begin{prop}
\label{proprepr} Each $\mathrm E^s\mathrm{SLCY}_{\mathrm{lc}}(\Upsilon,\nu)$ represents the strongly embedded stratum functor of the piece $J_{\Upsilon,\nu}$.
\end{prop}
\begin{proof}
For this problem, we consider the restriction of the universal family to $J_{\Upsilon,\nu}$. Since the bases of the families are reduced, we have that every morphism with image in this reduced locally closed subscheme factors uniquely through it.
\end{proof}

\subsection{The coarse moduli space}
\label{sec46}

\begin{theo}
\label{thrmoduli} Each chosen stratum has a quasi-projective coarse moduli space $\mathrm{SLCY}_{\mathrm{lc}}(\Upsilon,\nu)$.
\end{theo}
\begin{proof}
We denote by $\mathrm{SLCY}_{\mathrm{lc}}(\Upsilon,\nu)$ the image of the stratum above in the coarse moduli space defined by $\Xi$. The points of $\mathrm{SLCY}_{\mathrm{lc}}(\Upsilon,\nu)$ are the isomorphism classes of the surfaces in the stratum. Since stabilizers are the automorphism groups of stable pairs $(X,D+\tau C)$ with ample $K_X+D+\tau C$, we can see that they must be finite.

For a family $f\colon((\mathcal X,\mathcal D),\mathcal C)\to\mathcal S$, the sheaf
\begin{center}
$ f_*\mathcal O_{\mathcal X} \bigl(r(K_{\mathcal X/\mathcal S}+\mathcal D+\tau\mathcal C)\bigr)$
\end{center}
is locally free of rank $N+1$ on the locus where cohomology commutes with base change. We trivialize it on an open cover. In this way, we get strong embeddings in $\mathrm E^s\mathrm{SLCY}_{\mathrm{lc}}(\Upsilon,\nu)$. On the overlaps, these embeddings are related by $G_N$-valued transition functions. Thus their quotient maps glue uniquely to a morphism $\mathcal S\to\mathrm{SLCY}_{\mathrm{lc}}(\Upsilon,\nu)$. This gives the coarse moduli property.

Moduli spaces of stable minimal models are quasi-projective \cite{Kol2,Bir2}. The ambient coarse quotient is geometric. So the invariant closure of the locally closed locus and its invariant closed complement have closed images. The image of the locus is the difference of these two closed images. Hence $\mathrm{SLCY}_{\mathrm{lc}}(\Upsilon,\nu)$ is quasi-projective.
\end{proof}
\begin{proof}[Proof of Theorem \ref{thrcoarse}]
By Theorem \ref{thrmoduli}, each chosen locally closed piece represented by $\Upsilon$ has a quasi-projective coarse moduli space.
\end{proof}

\subsection{Admissible pencils}\label{sec47}

Let $\Upsilon_F=(\Upsilon,H,\sigma_F)$ be the tuple of Construction \ref{constr2}. Here $H$ is the class of the pencil. We work on a strongly embedded stratum $I=J_{\Upsilon,\nu}$, refined as in Lemma \ref{lemrevMMP}. We make an \'etale base change, and we still write $I$ for the new base. Then we have the family $p:\mathcal Y_0\to I$ of initial surfaces with its marked cluster. We also have the blow-up of the cluster and the contraction defined by the polarization. Let $\mathcal E\subset\mathcal Y_0$ be the labelled smooth boundary member of class $H$. For a pencil in $|C_\infty|$, we have $H=C_\infty$. For an anticanonical pencil, we have $H=-K_{Y_0}$. We put $q=p|_{\mathcal E}$, $\mathcal A=\mathcal O_{\mathcal E}(\mathcal E)$ and $\delta=H^2$. In the first case, we get $\delta=e$. In the second case, we get $\delta=K_{Y_0}^2$. After one more \'etale base change, we denote by $p_1,\ldots,p_m$ the proper points of $\mathcal K$ on $\mathcal E$. For each $p_i$, the points of $\mathcal K$ along $\mathcal E$ are the first $r_i$ points of $\mathcal E$ at $p_i$. We put $\mathcal B=\sum_i r_i p_i$. The support of $\mathcal B$ is finite \'etale over $I$. The integers $r_i$ are constant over the stratum. The sheaf $q_*\mathcal A$ is locally free and its formation commutes with base change. As in Subsection \ref{sec42}, the projective bundle $\mathbb L_{\mathcal A}=\mathbb P_I((q_*\mathcal A)^\vee)$ parametrizes the divisors in $|\mathcal A|$. We denote its universal divisor by $\mathcal T_{\mathcal A}$. This divisor is finite and flat over $\mathbb L_{\mathcal A}$. This gives us that the condition $\mathcal T_{\mathcal A}\leq\mathcal B$ is a closed condition. Indeed, it is the vanishing of the canonical section of $\mathcal O(\mathcal B)$ on $\mathcal T_{\mathcal A}$. This implies that the subscheme
\begin{center}
$\widehat I_{\Upsilon_F} =\{\mathcal T_{\mathcal A}\leq\mathcal B\}_{\mathrm{red}} \subset\mathbb L_{\mathcal A}$
\end{center}
is projective over $I$ with its geometric points are the divisors $T\in|\mathcal A_b|$ with $T\leq\mathcal B_b$.

After a stratification for cohomology and base change, the restriction of sections to $\mathcal E$ gives the exact sequence
\begin{center}
$ 0\longrightarrow\mathcal O_I s_{\mathcal E} \longrightarrow p_*\mathcal O_{\mathcal Y_0}(\mathcal E) \xrightarrow{\rho}q_*\mathcal A\longrightarrow0.$
\end{center}
Over a point $(b,T)$, the pencil $P_T=\mathbb P(\rho_b^{-1}(\mathbb C s_T))\subset|\mathcal E_b|$ has base scheme $T$. Let $k_i$ be the multiplicities of $T$, so that $T=\sum_i k_i p_i$ with $0\leq k_i\leq r_i$ and $\sum_i k_i=\delta$. In this situation, the cluster $\mathcal K_T$ of base points of $P_T$ consists of the first $k_i$ points of $\mathcal E$ at $p_i$, each of weight one. This means that $\mathcal K_T\subset\mathcal K$. On the blow-up $\mu$ of $\mathcal K$, the resolved pencil has the fibre class
\begin{center}
$F_T=\mu^*H-\sum_{\alpha\in\mathcal K_T}E_\alpha^*,$
\end{center}
and we have
\begin{center}
$C\cdot F_T=L\cdot H-\sum_{\alpha\in\mathcal K_T}m_\alpha.$
\end{center}
We say that $P_T$ is admissible if $F_T\cdot V=0$ for every curve $V$ contracted by the polarization. The admissible pencils form a locally closed subset $I^{\mathrm{st}}_{\Upsilon_F}\subset\widehat I_{\Upsilon_F}$. Since $T$ determines $P_T$, this subset is the scheme of admissible pencils. Hence the universal admissible pencil descends through the contraction defined by the polarization and we get a family
\begin{center}
$((\mathcal X_{\Upsilon_F},\mathcal D_{\Upsilon_F}), \mathcal C_{\Upsilon_F})\longrightarrow I^{\mathrm{st}}_{\Upsilon_F}.$
\end{center}
Here, we also get a morphism $\varphi_{\Upsilon_F}:\mathcal X_{\Upsilon_F}\to\mathcal Z_{\Upsilon_F}$ and on each fibre, this morphism is induced by $P_T$.

The reconstructed stable pair over $I$ is strongly embedded. So we have the forgetful classifying morphism $\Phi_{\Upsilon_F}:I^{\mathrm{st}}_{\Upsilon_F}\to \mathrm E^s\mathrm{SLCY}_{\mathrm{lc}}\Xi$. Let $J_{\Upsilon_F}$ be the constructible union of these forgetful images. The points of $J_{\Upsilon_F}$ are the surfaces which have an admissible pencil. We choose a finite invariant partition $J_{\Upsilon_F}=\coprod_\nu J_{\Upsilon_F,\nu}$ into reduced locally closed pieces. As in Definition \ref{defstrata}, we put $\mathrm E^s\mathrm{LCYF}_{\mathrm{lc}}(\Upsilon_F,\nu) :=J_{\Upsilon_F,\nu}$. Notice that each divisor $\mathcal B_b$ has finitely many subdivisors of degree $\delta$. Thus, the projective morphism $\widehat I_{\Upsilon_F}\to I$ is finite. We then refine the reduced base using the allowed multiplicity vectors and the admissibility conditions. Then, the schemes of admissible pencils are finite \'etale on each stratum.

Now we work with the case when $q=1$, $R=0$ and the weights satisfy \eqref{eqmultbnd}. Notice that on the minimal resolution of a stable fibre, the fibre class $F_T$ of an admissible pencil is the pullback $\rho^*L_T$ of a unique line bundle $L_T$ on the stable fibre. Then, we have a lemma
\begin{lemm}\label{lempencil}
The classes $[L_T]$ of the admissible pencils form a closed subscheme of $\Pic_{\mathcal X/I}$, which is finite \'etale over $I$ and invariant under $G_N$.
\end{lemm}
\begin{proof}
Let $V$ be the minimal resolution of a stable fibre. Let $V\to\Gamma$ be its Albanese fibration. Over each $p_i$, the fibre of $V\to\Gamma$ is a chain of $r_i+1$ smooth rational curves. This chain joins the two elliptic boundary curves. The curve $F_T$ meets exactly one curve of this chain, namely the $(k_i+1)$-th curve. The weights satisfy \eqref{eqmultbnd}. So by contracting every chain, we get a presentation of type $\Upsilon$. We can start from the $(-1)$-curve which meets $C_0$, or from the $(-1)$-curve which meets $C_1$. Hence there are exactly two minimal models. They differ by the exchange of the two elliptic curves. This exchange replaces $A$ by $\mathcal{O}_\Gamma(\sum_ir_ip_i)\otimes A^{-1}$. It reverses each chain. So it replaces $T$ by $T'=\sum_i(r_i-k_i)p_i$, and it does not change the class $F_T$. This implies that the classes $F_T$ of the admissible pencils depend only on the stable fibre.

The morphism defined by the pencil on the minimal resolution is constant on the curves contracted by $\rho$. Thus it factors through the stable fibre. Hence we have $\mathcal{O}(F_T)=\rho^*L_T$ for a unique line bundle $L_T$. In this way, we get a morphism to the relative Picard scheme. By \cite[Theorem~4.8]{Kle05}, this scheme is separated. Notice that distinct pencils in one presentation have distinct fibre classes. Also, the scheme of admissible pencils is finite \'etale over the cover. Thus this morphism is proper, unramified and universally injective. By \cite[Tag~04XV]{Stacks}, it is a closed immersion. Its image depends only on the stable family. Hence, this image is the same for every cover.
\end{proof}
Next, we have a property on the \'etale covers. 
\begin{prop}\label{proppen1}
Over each of these \'etale covers, the scheme of admissible pencils is locally closed and quasi-finite over the cover. Moreover, it becomes finite \'etale after a refinement of the cover.
\end{prop}

For example, let $\mathcal E_b\simeq\mathbb P^1$, and let all divisors $T\leq\mathcal B_b$ of degree $\delta$ be admissible. Then the fibre has $[z^\delta]\prod_{i=1}^m(1+z+\cdots+z^{r_i})$ points. For $r_i=1$, this number is $\binom{m}{\delta}$. When $\mathcal E_b$ is elliptic, we also impose $\mathcal O_{\mathcal E_b}(T)\simeq\mathcal A_b$.

\begin{proof}
Notice that the multiplicity vectors are bounded. Hence, they give us finitely many reduced loci with the given linear equivalence. On each of them, we have that the admissibility conditions above are locally closed. The restriction sequence gives the universal pencil with universal base locus, hence we can resolve it by blowing up its cluster of base points. Finally, we can refine the base by the allowed vectors above. By this way, we get the finite \'etaleness after refinement.
\end{proof}
From Proposition \ref{proppen1}, we will get an important description for the moduli space of stable log Calabi-Yau surfaces with pencil structure which is an analogue of the non-fibration case.
\begin{prop}\label{proppen2}
$J_{\Upsilon_F}$ is a finite union of invariant locally closed subsets of the ambient fine moduli space.
\end{prop}
\begin{proof}
By Proposition \ref{proppen1}, the schemes of admissible pencils are of finite type. Hence by Chevalley's theorem and the same argument as in Proposition \ref{propimage}, the image $J_{\Upsilon_F}$ is constructible.
\end{proof}
When we study a particular stratum, we can drop $\nu$ from the notation.

In the next theorem, we assume that our $\Upsilon$ has $q=1$, $R=0$ and the cluster weights satisfy $0<m_\alpha<a$.
\begin{theo}\label{thrpenmod}
For the given $\Upsilon$ with $q=1$, $R=0$ and the cluster weights satisfy $0<m_\alpha<a$, each chosen fibration stratum has a quasi-projective coarse space $\mathrm{LCYF}_{\mathrm{lc}}(\Upsilon,\nu)$.
\end{theo}
\begin{proof}
By Lemma \ref{lempencil}, we have descent data for the admissible pencils. Hence, their schemes glue to a finite $G_N$-equivariant scheme $P\to I$. We look at the quotient stack $[I/G_N]$. It is a separated Deligne--Mumford stack of stable pairs. Its stabilizer groups are finite. So they are linearly reductive over $\mathbb C$. By \cite[Theorem 3.2 and Corollary 3.3(a)]{AOV08}, this stack is a quotient $[U/H]$ \'etale locally over its coarse space $M$. Here $H$ is a finite group, and $U\to M$ is finite. The pullback of the scheme of admissible pencils to this chart is a finite $H$-equivariant scheme $V\to U$. On an affine coarse chart $\Spec R$, we write $V=\Spec A$, where $A$ is the coordinate ring of $V$. Since $A$ is finite over $R$, the ring $A^H$ is also finite over $R$. The quotients $\Spec(A^H)$ glue to the coarse space of $[P/G_N]$. This coarse space is finite over $M$. Theorem \ref{thrmoduli} gives the quasi-projectivity of $M$. Hence the finite coarse space is also quasi-projective.
\end{proof}

\subsection{Moduli spaces of elliptic fibrations}
\label{sec48}

For this construction and for the following ones, we take finite invariant locally closed partitions of the final strongly embedded images. For a chosen piece, we add $\nu$ to the moduli notation. When a piece has been fixed, we drop this index.

We fix $\Upsilon_F$, and we apply Construction \ref{constr3}. For this, we choose a smooth fibre $E$, and we blow up the reduced divisor $E\cap C$. We put
\begin{center}
$N_\ell=C\cdot F_T=L\cdot H- \sum_{\alpha\in\mathcal K_T}m_\alpha>0.$
\end{center}
Here $H$ is the pencil class, and $\mathcal K_T$ is its cluster of base points of Subsection \ref{sec47}. Let $v_\ell$ be the square of the final polarization. We put
\begin{center}
$\sigma_\ell(t)=\vol(K_{S_\ell}+tC_\ell) =2N_\ell t+v_\ell t^2.$
\end{center}
We denote the fibration tuple by $\Upsilon_\ell=(\Upsilon_F,N_\ell,v_\ell,\sigma_\ell)$. We denote its tuple of embedding invariants by
\begin{center}
$\Xi_\ell=(2,\Phi,\mathcal V_{\mathrm{fib}},\sigma_\ell(t), \tau_\ell,r_\ell,N'_\ell).$
\end{center}
We work over $I^{\mathrm{st}}_{\Upsilon_F}$ with its universal pencil. Let $E$ be a fibre in the smooth open part of the relative pencil base. We choose it such that $E\cap C$ is reduced of length $N_\ell$. First, we take the universal blow-up of $E\cap C$. Then we add a small positive multiple of the transformed polarization to the log canonical divisor. We take the contraction defined by this divisor. In this way, we get the stable family with polynomial $\sigma_\ell$. Local projective coordinates give classifying morphisms to the ambient space for $\Xi_\ell$. We move their images by $\mathrm{PGL}_{N'_\ell+1}$. This gives the locus $\mathrm E^s\mathrm{CHS}^{\Upsilon_\ell}_\ell$.
\begin{prop}\label{propelim}
The locus $\mathrm E^s\mathrm{CHS}^{\Upsilon_\ell}_\ell$ is constructible, and it has a finite invariant partition into locally closed subsets.
\end{prop}
\begin{proof}
We have that the smoothness of the fibre is an open condition. Then, the reducedness and the length of $E\cap C$ are locally closed conditions which gives us the precise locally close stratum. On this stratum, the centre $E\cap C$ is finite \'etale over the base which lies in the smooth locus of $E$. Thus, its blow-up commutes with base change. By taking the relative contraction, we get a family of finite type. By Chevalley's theorem for the local classifying maps and for the changes of projective coordinates, the locus is constructible. Hence, we get the partition in the usual way into the union of finitely many locally closed subsets.
\end{proof}
\begin{prop}\label{propelmod}
Each chosen elliptic construction stratum has a quasi-projective coarse moduli space $\mathrm{CHS}^{\Upsilon_\ell}_\ell(\nu)$.
\end{prop}
\begin{proof}
First, we take a chosen stratum, of which is reduced, invariant and locally closed. We restrict the universal family to this stratum. As in Proposition \ref{proprepr} and Theorem \ref{thrmoduli}, the considering family represents the strongly embedded stratum functor. We take its geometric quotient by $\mathrm{PGL}_{N'_\ell+1}$. It is locally closed in the coarse moduli space for $\Xi_\ell$. Notice that we can trivialize the embedding locally in families. Hence the quotient maps agree on the overlaps. This gives the coarse moduli property and thus gives the quasi-projectivity.
\end{proof}

\section{Rationality of moduli spaces}\label{sec5}

For birational questions, we only need the dense locally closed stratum which contains the generic point of each irreducible construction image of Section \ref{sec4}. Indeed, this stratum and the reduced closure of the image have the same function field. For this reason, we drop the index $\nu$ of this stratum.

Let $M$ be a reduced algebraic space of finite type with a finite partition into locally closed subspaces.
\begin{lemm}\label{lemstrata}
Every irreducible component $Z$ of $M$ contains an irreducible part $T$ of one stratum as a dense open subset. In particular, if $T$ is rational or unirational, then so is $Z$.
\end{lemm}
\begin{proof}
We divide each stratum into irreducible locally closed parts. We choose the part which contains the generic point of $Z$. The closure of this part is $Z$. The part is locally closed. So it is open in $Z$. Hence the two have the same function field.
\end{proof}

The quotients in this section are birational to the unmarked moduli spaces of Section \ref{sec4}. We denote these rational quotients by $\mathrm M(\Upsilon)$.

By \cite[Lemma 2.8]{AAB}, every generic polarized maximal log Calabi-Yau surface in the setting of Alexeev, Arg\"uz and Bousseau comes from a toric pair by the construction of that lemma. Lemma \ref{lemtoric} gives the boundary, the polarization and the marked cluster on this toric model. Here the conditions of this cluster also determine the complete section space.

\subsection{Quotients}\label{sec51}
For the slice method and the no-name lemma, let $G$ acts on $X$. We denote its rational quotient by $X/G$. We give the notation $X_{\min}$ for the initial ruled surface as in Section \ref{sec3} and denote its universal family by $\mathcal{X}_{\min}$.
\begin{defi}\label{defquot}
For a dominant rational map $\phi:X\dasharrow Y$, we denote the pullback by $\phi^*:k(Y)\to k(X)$. We call an action $\alpha$ of $G$ on $X$ \textit{regular} if $\alpha$ is a morphism. We also call it \textit{almost free}, or \textit{generically free}, if it is free on a dense open subset. 

Now write $g.x=\alpha(g,x)$. A homomorphism $G\to\Bir X$ defines a \textit{rational action} when $G\times X\dasharrow X$ is rational. A dominant rational map $\pi:X\dasharrow X/G$ is called a \textit{rational quotient} if $\pi^*k(X/G)=k(X)^G$. On a nonempty invariant open subset, the action is regular, and the quotient is geometric.
\end{defi}
We use the following form of the slice method \cite[Proposition 2.3]{Ma}.
\begin{lemm}\label{lemslice}
Let $f:U\to V$ be dominant and $G$-equivariant, and let the action on $V$ be almost transitive. Then for general $v\in V$ with stabilizer $G_v$, we have $U/G\sim_{\mathrm{bir}}U_v/G_v$.
\end{lemm}
\begin{lemm}[{\cite[Proposition 2.4]{Ma}}]\label{lemnoname}
For a $G$-linearized vector bundle $E\to U$ with almost free action on $U$, the quotient $E/G\to U/G$ is birationally a vector bundle.
\end{lemm}
We also need the projective version of the no-name lemma. Let an algebraic group $G$ act on an irreducible variety $U$. Let $E\to U$ be a $G$-linearized vector bundle of rank $N+1$. We put $G_0=\ker(G\to\Aut(U))$ and $\overline G=G/G_0$. So every element of $G_0$ fixes every point of $U$. We assume that $\overline G$ acts almost freely on $U$. We also assume that $G_0$ acts by scalars on the fibres of $E$ through a character $\alpha:G_0\to\mathbb G_m$. Moreover, we assume that there is a $G$-linearized line bundle $L\to U$ with the same $G_0$-character. This means that $g\cdot v=\alpha(g)v$ on the fibres of both bundles for every $g\in G_0$. In this setting, we have the following result from \cite[Proposition 2.5]{Ma}.
\begin{prop}\label{propnname}
Under these assumptions, we have $\mathbb P_U(E)/G\sim_{\mathrm{bir}}\mathbb P^N\times(U/G)$.
\end{prop}

\begin{proof}
We put $F=E\otimes L^{-1}$. The bundles $E$ and $L$ have the same $G_0$-character. This tells us that $g\in G_0$ acts on the fibres of $F$ by $\alpha(g)\alpha(g)^{-1}=1$. Hence, the linearization of $F$ factors through $\overline G$. Notice that over the open subset where $\overline G$ acts freely, the quotient $F/\overline G\to U/\overline G$ is a vector bundle of rank $N+1$ and its generic fibre is a vector space of dimension $N+1$ over $k(U/G)$. By choosing a basis, we can identify its projectivization with $\mathbb P^N_{k(U/G)}$. Since $\mathbb P_U(E)=\mathbb P_U(F)$, this gives the statement.
\end{proof}
For linear actions, we use the following theorem of B\"ohning on invariant fields \cite[Corollary 1.2.2.2]{Boeh}.
\begin{theo}\label{thrbohn}
Let an abelian group $G$ of semisimple elements act linearly on a finite-dimensional complex vector space $V$. Then $\mathbb C(V)^G$ is purely transcendental over $\mathbb C$.
\end{theo}
In particular, this theorem applies when $G$ is a torus. An action on a torus with coordinates $z_1,\ldots,z_n$ is called monomial if every group element sends each $z_i$ to $c_i\prod_jz_j^{a_{ij}}$, where $c_i\in\mathbb C^*$ and $(a_{ij})\in\mathrm{GL}_n(\mathbb Z)$.
\begin{lemm}\label{lemmoninv}
Let $\tau$ be an involution of $\mathbb C(x_1,\ldots,x_n)$ which sends each $x_i$ to a nonzero constant multiple of a Laurent monomial. Then the invariant field $\mathbb C(x_1,\ldots,x_n)^{\langle\tau\rangle}$ is rational.
\end{lemm}
\begin{proof}
Let $M$ be the lattice of Laurent monomials in $x_1,\ldots,x_n$. The involution $\tau$ acts on $M$ through its linear part. By \cite{Reiner57}, $M$ is a direct sum of lattices of rank one and of rank two. On a lattice of rank one, $\tau$ acts by $\pm1$. On a lattice of rank two, $\tau$ exchanges a basis. The involution $\tau$ acts multiplicatively, and $\tau^2=1$. Hence, the coordinates transform as $x\mapsto \pm x$, $y\mapsto c/y$ or $(z,u)\mapsto (du,d^{-1}z)$ with $c,d\in\mathbb C^*$. After rescaling $y$ and $u$, we can assume that $c=d=1$. Now, the functions $w=(y-1)/(y+1)$, $s=z+u$ and $t=z-u$ satisfy $w\mapsto -w$, $s\mapsto s$ and $t\mapsto -t$. Thus $\tau$ multiplies each generator of $\mathbb C(M)$ by $\pm1$. Let $v_1,\dots,v_r$ be the generators which are multiplied by $-1$. Then we can see that the fixed generators together with $v_1^2$ and $v_i/v_1$ for $i\ge 2$ generate the invariant field. Hence, we get the invariant field is purely transcendental.
\end{proof}
\subsection{Unirationality of moduli spaces of log Calabi-Yau surfaces}
\label{sec52} In this subsection, we study the unirationality of the moduli spaces $\mathrm{SLCY}_{\mathrm{lc}}(\Upsilon,\nu)$ for $q=0$ and of $\mathrm{SLCY}_{\mathrm{lc}}\Upsilon$ for $q=1$.

We first consider the case $q=0$. For a fixed taut pair, the group $\Aut(X_{\min},D)$ acts on $I^{\mathrm{st}}_\Upsilon$. For the unmarked moduli space, we will also compare the initial surfaces in the elementary transformation class $\Upsilon$ of Section \ref{sec4}.

For simplicity, we assume that $I^{\mathrm{st}}_\Upsilon$ is irreducible. To turn isomorphisms of the constructed surfaces into a group action, we add markings to $I^{\mathrm{st}}_\Upsilon$ in the sense of Gross, Hacking and Keel \cite[Definition 1.11]{GHK1}. On the minimal resolution of the constructed stable pair, we can see that such a marking consists of a Picard marking and a smooth point on each labelled boundary component. With these data, $I^{\mathrm{st}}_\Upsilon$ becomes the space $\widetilde I^{\mathrm{st}}_\Upsilon$. The changes of markings are described in the marked families of \cite[Constructions 5.12--5.13]{GHK1}. Now we look at the birational maps between the initial surfaces and also the corner blow-ups. It is natural for us to take those which preserve the $\Upsilon$ structure. Let $G_\Upsilon$ be the group generated by these maps and by the changes of the markings. For a rational surface with a nodal anticanonical cycle, we use the maps of those marked families with the chosen section included. Further corner blow-ups give other initial surfaces. We compare their boundaries, polarizing systems and clusters on this minimal resolution. We forget the markings in
\begin{center}
$ \mathrm M(\Upsilon):= \widetilde I^{\mathrm{st}}_\Upsilon//_{\mathrm{rat}}G_\Upsilon. $
\end{center}
For a fibration tuple, the parameter space contains the chosen pencil. We require that $G_\Upsilon$ preserves it. In the CHS constructions, we require that we can recover some data in families, as in Section \ref{sec4}. These data are the resolution, the initial surface, the boundary, the polarization, the chosen pencil and any chosen fibre divisor.

\begin{prop}\label{propquot}
The quotient $\mathrm M(\Upsilon)$ is birational to the closure of the irreducible unmarked moduli stratum $\SLCY_{\mathrm{lc}}(\Upsilon,\nu)$.
\end{prop}
\begin{proof}
We work on a dense locus. On this locus, the resolutions commute with base change. Notice that a change of the initial surface or of the boundary markings preserves the boundary, the complete polarizing linear system and its chosen member, so it induces an isomorphism of the stable pairs of which we get by these contractions. Conversely, an isomorphism of the constructed stable pairs lifts to their minimal resolutions and thus identifies the boundaries, the polarizing systems and the chosen members. Hence, we can see that resolutions with their Picard and boundary markings can differ by an element of $G_\Upsilon$.

Thus, the above arguments tell us that the geometric fibres of the classifying map are the $G_\Upsilon$-orbits. This gives a dominant and generically injective classifying map. Hence, we have that Rosenlicht's quotient has the same function field as the moduli stratum. In particular, we get
\begin{center}
$ \mathbb C(\SLCY_{\mathrm{lc}}(\Upsilon,\nu)) =\mathbb C(\widetilde I^{\mathrm{st}}_\Upsilon)^{G_\Upsilon}. $
\end{center}
\end{proof}
In the next theorem, we take $q\in\{0,-\infty\}$ with nodal anticanonical boundary. The index $\nu$ is the index of the dense open set of $I^{\mathrm{st}}_{\Upsilon}$.
\begin{theo}\label{thrnodal}
The space $\mathrm{SLCY}_{\mathrm{lc}}(\Upsilon,\nu)$ is unirational.
\end{theo}

\begin{proof}
We fix the taut anticanonical pair for a reducible boundary. By Section \ref{sec42}, $I^{\mathrm{st}}_\Upsilon$ is a section bundle over the boundary configurations. We group the points which lie on the same labelled boundary component and which have the same multiplicity and the same labelled Enriques diagram. Elementary symmetric functions give coordinates on each such unordered collection. In these coordinates, the Abel--Jacobi--Picard condition is one monomial equation. So its base and its total space are both rational.

Now let $W_1$ be a cycle in $|-K_{X_{min}}|$ as in Construction \ref{constr1}. For an integral nodal boundary on $\mathbb P^2$, $\mathbb F_1$ or $\mathbb F_2$, we use the marked period family of \cite[Construction 5.7]{GHK1} over the torus $T_P=\Hom(P,\mathbb G_m)$ with $P=\Pic(X_{\min})$, and write $t^L$ for the value $t(L)$ of a period $t\in T_P$. Notice that the boundary marking identifies both $W_1^{\mathrm{sm}}$ and $\Pic^0(W_1)$ with $\mathbb G_m$.  After we order the distinct points of $Z=C_1|_{W_1}=\sum_{j=1}^N\mu_jp_j$, we choose their coordinates $z_j\in W_1^{\mathrm{sm}}\simeq\mathbb G_m$ such that the Abel--Jacobi--Picard condition becomes
\begin{center}
$\prod_jz_j^{\mu_j}=t^L.$
\end{center}
On $\mathbb F_2$, the curve $C_0$ does not meet $W_1$, so we also impose $t^{C_0}=1$. All these are character equations on the product $T_P\times(\mathbb G_m)^N$ of the period torus and the torus with coordinates $z_1,\ldots,z_N$. Hence, they define a finite union of translates of subtori. We take the nonempty open subset where the construction conditions hold. Over this subset, the cluster conditions of Section \ref{sec42} give a vector bundle of sections. Its projectivization is rational. Thus $I^{\mathrm{st}}_{\Upsilon}$ is rational. This gives the unirationality of $\mathrm{SLCY}_{\mathrm{lc}}\Upsilon$.
\end{proof}
Now let $q\in\{0,-\infty\}$. We will fix $X_{\min}$ and $L$. In this situation, we will also let $W_1\in|-K_{X_{\min}}|$ and $C_1\in|L|$ vary over a nonempty open subset. On this subset, both curves are smooth and meet transversely, so we take the cluster $\mathcal K$ which consists of these distinct points.
\begin{prop}\label{propsmbd}
With the given $\Upsilon$ above, the space $\mathrm{SLCY}_{\mathrm{lc}}\Upsilon$ is unirational.
\end{prop}
\begin{proof}
First, we fix the rational initial surface $X_{\min}$ and its polarizing line bundle $L$. After that, we take the nonempty open subset of $ |-K_{X_{\min}}|\times|L|$ above. Here, the cluster consists of the points of $C_1\cap W_1$, each blown up once and the pair $(W_1,C_1)$ determines the cluster. Hence, the universal blow-up and then the contraction give us a dominant classifying map from this rational parameter space to the image inside $\mathrm{ SLCY}_{\mathrm lc} \Upsilon$.
\end{proof}
For $q=1$, we first fix the elliptic curve $\Gamma$. We denote the tuple we get by $\Upsilon^\Gamma$. We also fix a line bundle $A$ of degree $e>0$ on $\Gamma$, and we put $X_{\min}=\mathbb{P}_\Gamma(\mathcal{O}_\Gamma\oplus A^{-1})$. Let $U\subset|aC_\infty|$ be the nonempty open subset of the smooth members $C_1$ which meet $C_\infty$ transversely at $ae$ distinct points. We let $\mathcal{K}$ consist of these points, each of weight one.
\begin{prop}\label{propfixel}
In this setting, the moduli space $\mathrm{M}(\Upsilon^\Gamma)$ is unirational.
\end{prop}

\begin{proof}
Multiplication by $e$ on $\Gamma$ is surjective. So every line bundle of degree $e$ on $\Gamma$ is a translate of $A$. Thus $\Gamma$ determines $X_{\min}$ with its boundary $C_0+C_\infty$ up to isomorphism. Moreover, $C_1$ and $C_{\infty}$ determine the cluster $\mathcal{K}$. Hence the open subset $U$ of a projective space parametrizes the construction and thus $U$ dominates $\mathrm{M}(\Upsilon^\Gamma)$.
\end{proof}

Now we study the case when $\Gamma$ can be varied in its parameterizing space. In order to adapt the same trick as in the $q=0$ case and extend the $q=1$ case, we will give a more simple version of the construction of Section \ref{sec4} over a birational parameter space.
\begin{construction}\label{constr5}
Over $U=\mathbb A^1_j\setminus\{0,1728\}$, we take the representatives
\begin{center}
$ \mathcal E_j:\ y^2z=x^3-\frac{3j}{j-1728}xz^2 -\frac{2j}{j-1728}z^3.$
\end{center}
By using the zero section, we identify $\mathcal P^{-e}:=\Pic^{-e}_{\mathcal E/U}$ with $\mathcal E$. We choose a normalized Poincar\'e bundle $\mathcal L$ of degree $-e$ \cite{Kle05}. We put
\begin{center}
$ \mathcal X=\mathbb P_{\mathcal E\times_U\mathcal P^{-e}} (\mathcal O\oplus\mathcal L).$
\end{center}
On each fibre, we have $-K_{X_{\Gamma,L}}\sim C_0+C_\infty$. We take the positive-section system and the polarization system on this family. We let $U^{\mathrm{st}}_\Upsilon$ be the stable construction locus of Definition \ref{defstable}. At each point of this locus, we keep the tuple $(\Gamma,L,C_\infty,C_1,\mathcal K)$. Now, we take the universal cluster blow-up and the relative contraction defined by the polarizing section ring in Proposition \ref{propfam}. From these, we get
\begin{center}
$ \widetilde{\mathcal X}\longrightarrow\mathcal S \longrightarrow U^{\mathrm{st}}_\Upsilon.$
\end{center}
On each resolved fibre, we have that $K_{\widetilde X}+\widetilde C_0+\widetilde C_\infty\sim0$ and $\widetilde C\cdot(\widetilde C_0+\widetilde C_\infty)=0$. Thus the polarization contracts the two elliptic curves in the boundary.

Notice that two different minimal elliptic ruled surfaces can give the same surface with CHS. We identify them by the action from the group $\mathcal G_e$. Here, this group consists of the isomorphisms of the base and of the automorphisms which preserve $C_0$ and $|C_\infty|$. We denote by $\mathrm M^{2\mathrm{ell}}_\Upsilon$ the coarse space of $[U^{\mathrm{st}}_\Upsilon/\mathcal G_e]$. It still remembers the elliptic base, since it maps to $U$. Now let $C_1\in|aC_\infty|$ be smooth. We assume that $C_1|_{C_\infty}=p_1+\cdots+p_{ae}$ has distinct points. In this case, the cluster $\mathcal K$ consists of these points. We blow up each of them once, with multiplicity one. Then over a fixed base, the quotient is the one in Proposition \ref{propfixel}.

Notice that multiplication by $e$ on $\Gamma$ is surjective. Hence, after a translation of $\Gamma$, we have $L\simeq\mathcal O_\Gamma(-eO)$. Since the positive section varies in its complete linear system, we get that the divisor $C_1|_{C_\infty}$ varies in $|\mathcal O_{C_\infty}(aC_\infty)|$ as $C_1 \sim aC_{\infty}$. Moreover, $C_1$ also varies among all lifts of this divisor in $|aC_\infty|$. These choices form projective bundles over the rational curve $U$. From this projective bundle, we take the open subset where the intersections are transverse. It is obviously admissible and nonempty, and it dominates the quotient. Hence the moduli image is unirational.
\end{construction}
Let
\begin{center}
$\Phi:\mathrm M^{2\mathrm{ell}}_\Upsilon\to \mathrm{SLCY}_{\mathrm{lc}}\Upsilon$
\end{center}
be the classifying morphism of the universal family above.
\begin{prop}\label{propbir}
For $q=1$, the morphism $\Phi$ is birational.
\end{prop}

\begin{proof}
First, the morphism $\Phi$ exists. Indeed, we have the construction as in \ref{constr5}, which we get $\mathrm M^{2\mathrm{ell}}_\Upsilon$ parametering the stable log Calabi-Yau surfaces with $\Upsilon$ information and this gives a morphism to the coarse moduli space $\mathrm{SLCY}_{\mathrm{lc}}$. Here, we work on a dense open subset where the initial ruled surface is unique. Hence, we can take two polarized surfaces with CHS in this subset and an isomorphism between them. This isomorphism lifts to their minimal resolutions. Thus, we can see that it identifies the two exceptional curves which are ellitptic. Now, we perform the MMP on the smooth surfaces and get $X_{\min}=\mathbb P_\Gamma(\mathcal O_\Gamma\oplus L)$ with $C_\infty$. This also gives us $C_1\in|aC_\infty|$ and the marked anticanonical cluster. Thus, we can see that the two tuples differ by an automorphism. This implies that two points are the same after taking the quotient and defines a birational map to the image of $\Phi$.
\end{proof}

\begin{proof}[Proof of Theorem \ref{thrunirat}]
The nodal-boundary case follows from Theorem \ref{thrnodal}. Now let the boundary consist of smooth elliptic curves. When the initial surface is rational, we apply Proposition \ref{propsmbd}. When the initial surface is elliptic ruled, Construction \ref{constr5} gives a rational parameter space which dominates the moduli image. By Proposition \ref{propbir}, this image is $\mathrm{SLCY}_{\mathrm{lc}}\Upsilon$. Finally, Proposition \ref{propquot} gives the same conclusion for $\mathrm M(\Upsilon)$.
\end{proof}
Next, we study Conjecture 1.1 of \cite{HKY}. We consider triples $(X,D=E_1+\cdots+E_n,\Theta)$ of one of the following two types.

\noindent\textup{(1)} The complex projective surface $X$ is log Calabi-Yau with at most canonical singularities. It is the canonical modification of a surface $S$ with CHS. The boundary $D$ has no $0$-strata. The semiample divisor $\Theta$ defines the morphism $X\to S$. We get its admissible tuple by letting $W_1$ and $C_1\in|L|$ vary in their complete linear systems, with
\begin{center}
$C_1|_{W_1}=p_1+\cdots+p_N.$
\end{center}
Here we have $N=L\cdot W_1$. The points $p_i$ are distinct smooth points of $W_1$, and all intersections are transverse. The cluster consists of these proper points, each has multiplicity one. In this case, we have $R=0$.

\noindent\textup{(2)} The surface $X$ is rational, $D\in|-K_X|$ is a normal crossing divisor with at least one $0$-stratum, and $\Theta$ is ample.
\begin{prop}\label{prophky}
For each of the two types, the moduli space $\mathrm{SLCY}_{\mathrm{lc}}\Upsilon$ is unirational.
\end{prop}

\begin{proof}
Theorem \ref{thrnodal} implies immediately the type \textup{(2)}. For type \textup{(1)}, we apply Proposition \ref{propsmbd} in the rational case. In the elliptic ruled case, we apply Construction \ref{constr5} together with Proposition \ref{propbir}.
\end{proof}
In case (1), the boundary $D$ is a disjoint union of smooth elliptic curves. The singularities of $S$ on this boundary are simple elliptic and they are extracted by the canonical modification. The Du Val points of Theorem \ref{thrm2crv} lie away from the boundary. The canonical modification does not extract them. Hence $X$ is smooth if and only if $S$ has no Du Val point. In this case, the canonical modification is the minimal resolution, and these triples give a related moduli problem with no $0$-strata and a semiample polarization.

Now we come to fibrations. We study the general strata of the scheme of admissible pencils of Subsection \ref{sec47}, on which the base divisor of the pencil is reduced. 

To begin with, we will look at the rational case with non-smooth boundary.
\begin{prop}\label{propnodpf}
For $q=0$ with non-smooth boundary, $\mathrm{LCYF}_{\mathrm{lc}}\Upsilon_F$ is unirational.
\end{prop}
\begin{proof}
Since the boundary is a cycle of rational curves, the base of the marked period family of Theorem \ref{thrnodal} is a torus. Now we look at two conditions on this torus. The first one is the linear equivalence for $T$ and the second one is the Abel--Jacobi--Picard condition for $C_1|_{W_1}$. Both of them are character equations. Hence, they cut out a finite union of translates of subtori. Over a dense open subset of each translate, the pencils are admissible. There the spaces of pencil sections and of polarizing sections have constant dimension. So these sections form a tower of two projective bundles. The total space of this tower is rational. These finitely many rational varieties together dominate $\mathrm{LCYF}_{\mathrm{lc}}\Upsilon_F$. Since $\mathrm{LCYF}_{\mathrm{lc}}\Upsilon_F$ is irreducible, one of them dominates it and thus implies the unirationality of $\mathrm{LCYF}_{\mathrm{lc}}\Upsilon_F$.
\end{proof}

\begin{construction}
Let $q=1$. We pull the ruled family back to $p:\mathcal X\to I$ with $I=U^{\mathrm{st}}_\Upsilon$. We put $q=p|_{\mathcal C_\infty}$ and $\mathcal A=\mathcal O_{\mathcal C_\infty}(\mathcal C_\infty)$. Let $\mathcal B=\sum_i r_i p_i$. Here $r_i$ is the number of points of $\mathcal K$ over $p_i$, as in Subsection \ref{sec47}. We count $p_i$ itself and the points infinitely near to it. We have that the support of $\mathcal B$ is finite \'etale over $I$. We take the relative complete system $\mathbb L_{\mathcal A}=\mathbb P_I((q_*\mathcal A)^\vee)$. Then we put
\begin{center}
$ \widehat U_{\Upsilon_F}= \mathbb L_{\mathcal A}\times_{\Div^e(\mathcal C_\infty/I)} \Div^e(\mathcal B/I).$
\end{center}
This fibre product has its reduced structure. We restrict it to the admissible locus $U^{\mathrm{st}}_{\Upsilon_F}$. On each geometric point of $U^{\mathrm{st}}_{\Upsilon_F}$, every curve which is contracted by the polarization must be the vertical divisor of the pencil.

From the restriction sequence of sheaves
\begin{center}
$ 0\to\mathcal O_I s_\infty\to p_*\mathcal O_{\mathcal X}(\mathcal C_\infty) \xrightarrow{\rho}q_*\mathcal A\to0$
\end{center}
we have that the universal pencil $\mathbb P(\rho^{-1}(\mathcal H_T))$. Here $\mathcal H_T$ is the line of the chosen divisor. Then, the universal blow-up step resolves the base scheme $T$ of this pencil. Since the contracted curves are vertical, the pencil appears on the contracted variety through the relative contraction defined by the polarization. We take the quotient by $\mathcal G_e$ and get $\mathrm M^{2\mathrm{ell},F}_{\Upsilon_F}$. This space is birational to $\mathrm{LCYF}_{\mathrm{lc}}\Upsilon_F$, because we keep $T$ in the reverse construction. Moreover, forgetting $T$ is quasi-finite because a fixed finite divisor has only finitely many subdivisors of degree $e$ with bounded multiplicities.
\end{construction}
In the next proposition, let $q=1$ and $\mathcal A=\mathcal O_{C_\infty}(C_\infty)$. We let $T\in|\mathcal A|$ and $Z\in|\mathcal A^{a-1}|$ vary among the reduced divisors with disjoint supports and also let $C_1$ vary among the smooth members of $|aC_\infty|$ with
\begin{center}
$C_1|_{C_\infty}=T+Z.$
\end{center}
Now take the cluster $\mathcal K$ which consists of all points of $T+Z$ with each blown up once with multiplicity one. We keep the positive-section pencil with base divisor $T$.
\begin{prop}\label{propellpf}
If this admissible locus is non empty, then by Lemmas \ref{lemadj} and \ref{lemvert} its moduli space $\mathrm{LCYF}_{\mathrm{lc}}\Upsilon_F$ is unirational.
\end{prop}
\begin{proof}
We use the translated representatives of Construction \ref{constr5}. We choose $T$ and $Z$ together with a lifting section as above. These choices form a tower of projective bundles over $U$. The nonempty admissible open subset of this tower dominates the quotient with the chosen fibration. So it also dominates the moduli image.
\end{proof}
\subsection{Rationality for smooth polarized anticanonical pairs}
\label{sec53}

\begin{lemm}\label{lemaffine}
Let $\mathcal T$ be a torus, and let $\sigma$ be an automorphism of $\mathcal T$ of order at most two which is given by a monomial transformation. Let $E\to\mathcal T$ be a $\sigma$-equivariant affine bundle on which $\sigma$ acts by affine maps between the fibres. Then the invariant field $\mathbb C(E)^{\langle\sigma\rangle}$ is rational.
\end{lemm}
\begin{proof}
We put $K=\mathbb C(\mathcal T)$. If $\sigma$ acts nontrivially on $K$, then affine Galois descent gives
\begin{center}
$ \mathbb C(E)^{\langle\sigma\rangle} =K^\sigma(u_1,\ldots,u_m),$
\end{center}
Here the $u_i$ are algebraically independent \cite[Theorem 1]{HajjaKang95}. If $\sigma$ fixes $K$, we choose a $K$-point $p$ of the generic fibre. We take $\frac12p+\frac12\sigma(p)$ as the origin. The linear involution we get has the eigenvalues $1$ and $-1$. Hence we can choose coordinates $x_1,\dots,x_r$ with $\sigma(x_i)=x_i$ and $z_1,\dots,z_s$ with $\sigma(z_j)=-z_j$. If $s>0$, then the fixed coordinates together with $z_1^2,z_2/z_1,\ldots,z_s/z_1$ generate the invariant field. If $s=0$, then every coordinate is fixed. Thus in both cases, the invariant field is purely transcendental over $K^\sigma$. By Lemma \ref{lemmoninv}, $K^\sigma$ is rational.
\end{proof}

Let $(X,D)$ vary in a fix deformation type of smooth rational pairs with $D=\sum_iD_i\in|-K_X|$ is a singular normal crossing divisor with smooth components and $L$ is ample. Let $\Upsilon$ be the admissible tuple of a dense stratum in this deformation space. We fix an initial surface and the blow-ups of the clusters in this class. Here, we also let the cluster centres and the member of the complete system $|L|$ vary. 

We use the Picard and boundary markings of Subsection \ref{sec52}, with the cycle oriented. For the boundary marking, we normalize the coordinate $z_i$ on $D_i^{\mathrm{sm}}\simeq\mathbb G_m$ so that its two nodes are $0,\infty$ and $b_i=1$. This normalization uses the action of the torus $T_\Upsilon=\Aut^0(D)\simeq(\mathbb G_m)^r$ by $z_i\mapsto a_i z_i$. Here $a_i\in\mathbb G_m$, and $r$ is the number of components \cite[Lemma 2.1]{GHK1}. Let $F_\Upsilon$ be the image of $G_\Upsilon$ acting on $P$ and on the labelled oriented cycle. It fixes $L$ with $L^2>0$. Hence it acts on the negative definite lattice $L^\perp$. The cycle has finitely many symmetries. Thus $F_\Upsilon$ is finite. As in \cite[Constructions 5.12--5.13]{GHK1}, an element $g\in G_\Upsilon$ sends the period $t:P\to\mathbb G_m$ to $t\circ g^{-1}$ and it permutes the $z_i$. If $g$ reverses the orientation, then it also inverts the identification $\Pic^0(D)\simeq\mathbb G_m$ and sends $z_i\mapsto z_i^{-1}$. The kernel of $G_\Upsilon\to F_\Upsilon$ consists of the changes of boundary marking. Hence, we have
\begin{center}
$ 1\longrightarrow T_\Upsilon\longrightarrow G_\Upsilon \longrightarrow F_\Upsilon\longrightarrow1. $
\end{center}
Following \cite[Definition 1.6]{GHK1}, a root is a class $\alpha\in P$ which becomes the class of an internal $(-2)$-curve on a deformation equivalent pair. The reflection in a root $\alpha$ is $s_\alpha(\lambda)=\lambda+(\lambda\cdot\alpha)\alpha$. We take the roots $\alpha$ with $\alpha\cdot L=0$ whose reflections belong to $F_\Upsilon$. We denote by $W_\Upsilon$ the group generated by these reflections. These reflections fix $L$ and every $D_i$. Their action on the periods is given by
\begin{center}
$(s_\alpha\cdot t)(\lambda) =t(\lambda)t(\alpha)^{\lambda\cdot\alpha}.$
\end{center}
Conjugation by elements of $F_\Upsilon$ permutes these reflections. Hence $W_\Upsilon$ is a normal subgroup of $F_\Upsilon$. We put $H_\Upsilon=F_\Upsilon/W_\Upsilon$. In the real span of these roots, a chamber is a connected component of the complement of the hyperplanes $\alpha^\perp$. The group $W_\Upsilon$ acts simply transitively on the chambers. We identify $H_\Upsilon$ with the subgroup of $F_\Upsilon$ which preserves a fixed chamber. We get $F_\Upsilon=W_\Upsilon\rtimes H_\Upsilon$. Let $B_\Upsilon$ be the space which parametrizes the triples $(X,D,L)$ with their markings and the divisor $R=C|_D$. Here we forget the chosen member $C$. In the notation of Section \ref{sec42}, we keep the admissible data $(Y_0,D_0,C_\infty,L_0;\mathcal K,R)$, and we forget $C_1$. The boundary divisor varies under the condition $\mathcal O_D(R)\simeq L|_D$. For a lattice $M$, we denote by $\mathbb C(M)$ the fraction field of the group algebra $\mathbb C[M]$. We order the points $p_{ij}\in D_i$ of a general reduced divisor $C|_D$. By \cite[Constructions~5.12--5.13]{GHK1} and the Abel--Jacobi--Picard condition, the version of $B_\Upsilon$ with ordered points has function field $\mathbb C(\Lambda)$, where
\begin{center}
$\Lambda=\Bigl(P\oplus\bigoplus_{i,j}\mathbb Z e_{ij}\Bigr) \Big/\mathbb Z v$ with $v=L-\sum_{i,j}e_{ij},$
\end{center}
Here $e_{ij}$ corresponds to the coordinate of $p_{ij}$. The vector $v$ is primitive. Thus the lattice $\Lambda$ is torsion free. Let $d:\Lambda\to\mathbb Z^r$ be the boundary weight. It is given by $d(\lambda)=(\lambda\cdot D_i)_i$ for $\lambda\in P$ and by $d(e_{ij})=\varepsilon_i$. We put $n_i=L\cdot D_i$, which is the number of points $p_{ij}$ on $D_i$. Every $n_i$ is positive. Hence every $\varepsilon_i$ is the weight of some $e_{ij}$, and the map $d$ is surjective. An element $(a_i)_i\in T_\Upsilon$ multiplies $e^\lambda$ by $\prod_i a_i^{d_i(\lambda)}$. Since $d$ is surjective, $T_\Upsilon$ acts freely by translations on the torus with character lattice $\Lambda$. Therefore, we get $\mathbb C(\Lambda)^{T_\Upsilon}=\mathbb C(\Lambda_0)$ with $\Lambda_0=\ker d$. Let $W'=W_\Upsilon\times\prod_i S_{n_i}$.

\begin{prop}\label{propuprat}
Suppose that the deformation family above satisfies $|H_\Upsilon|\leq2$. Then $\mathrm M(\Upsilon)$ and its unmarked moduli image are rational. Their dimension is $\dim\widetilde I^{\mathrm{st}}_\Upsilon- \dim(G_\Upsilon\cdot\eta)$, where $\eta$ is the generic point.
\end{prop}
\begin{proof}
Since $L$ is ample, Kodaira vanishing gives $H^1(X,K_X+L)=0$. Thus the restriction sequence
\begin{center}
$ 0\longrightarrow\mathcal O_X(K_X+L) \longrightarrow\mathcal O_X(L) \longrightarrow\mathcal O_D(L)\longrightarrow0$
\end{center}
is surjective on sections. The projective sections with boundary divisor $R$ form an affine bundle $E_\Upsilon\to B_\Upsilon$.

We denote by $B_\Upsilon^{\mathrm{ord}}$ the version of $B_\Upsilon$ with ordered boundary points and $E_\Upsilon^{\mathrm{ord}}\to B_\Upsilon^{\mathrm{ord}}$ the pullback of $E_\Upsilon$. To give a proof for rationality, we will compute the quotient in three steps:\\
\\
\textbf{Step 1.} The group $T_\Upsilon$ only changes the boundary marking. It acts freely by translations on the torus with character lattice $\Lambda$. The sections of $K_X+L$ are the sections of $L$ which vanish on $D$. So the group $T_\Upsilon$ acts trivially on them. Thus the marking appears only in the normalization of the boundary section. More exactly, a change of marking multiplies the section of $L|_D$ with divisor $R$ by a character. By Lemma \ref{lemnoname}, this implies that $E^{\mathrm{ord}}_\Upsilon/T_\Upsilon\to B^{\mathrm{ord}}_\Upsilon/T_\Upsilon$ is birationally an affine bundle over the torus with function field $\mathbb C(\Lambda_0)$.\\
\\
\textbf{Step 2.} The boundary weight $d$ is $W'$-invariant. So the group $W'$ preserves $\Lambda_0$. The roots which define $W_\Upsilon$ and the differences $e_{ij}-e_{ik}$ span a root lattice $Q$ in $\Lambda_0$. This lattice injects into $\Lambda$, because $v$ is fixed and the span of the roots has no fixed vector. The group $W'$ acts faithfully on $Q$, hence also on $\Lambda_0$. Thus $W'$ acts generically freely on the torus with character lattice $\Lambda_0$. Over the open subset where this action is free, we have that the quotient of $E_\Upsilon^{\mathrm{ord}}/T_\Upsilon$ by $W'$ is an affine bundle over a variety with function field $\mathbb C(\Lambda_0)^{W'}$.

In order to compute this field, we let $\Phi\subset Q$ be the set of the roots of $W_\Upsilon$ and of the differences $e_{ij}-e_{ik}$ and also let $s_a$ be the reflection of $a\in\Phi$. The roots of $W_\Upsilon$ have square $-2$ in the negative definite lattice $L^\perp$ and the $e_{ij}-e_{ik}$ form root systems of type $A$. Thus, $\Phi$ is a root system of type ADE with Weyl group $W'$. On $\Lambda_0\otimes\mathbb R$, we put $\rho=|W'|^{-1}\sum_{w\in W'}w$ and $\pi=1-\rho$. For $\lambda\in\Lambda_0$ and $a\in\Phi$, we have $s_a(\lambda)-\lambda\in\mathbb Za$. Indeed, $s_\alpha(\lambda)-\lambda=(\lambda\cdot\alpha)\alpha$ for a root $\alpha$ of $W_\Upsilon$, and a transposition changes $\lambda$ by an integral multiple of $e_{ij}-e_{ik}$. Hence $\pi(\lambda)$ lies in the weight lattice $P_\Phi$ of $\Phi$. Since $\lambda=\rho(\lambda)+\pi(\lambda)$, we get that
\begin{center}
$\Lambda_0\subset\widehat\Lambda:=\rho(\Lambda_0)\oplus P_\Phi$
\end{center}
is a sublattice of finite index \cite{Lorenz01}. For $p\in P_\Phi$ and $w\in W'$, we have $w(p)-p\in Q\subset\Lambda_0$. Consequently, $W'$ acts trivially on $\widehat\Lambda/\Lambda_0$. We fix a base of $\Phi$ with fundamental weights $\lambda_1,\ldots,\lambda_r$. Let $\sigma_i\in\mathbb C[P_\Phi]$ be the sum of the $e^\mu$ with $\mu$ in the orbit $W'\lambda_i$. For the dominance order, its unique maximal term is $e^{\lambda_i}$. By \cite[Chapter~VI, \S3, no.~4, Theorem~1]{Bourbaki02}, the $\sigma_i$ are algebraically independent, and they generate $\mathbb C[P_\Phi]^{W'}$. Since $W'$ fixes $\rho(\Lambda_0)$, we get
\begin{center}
$\mathbb C(\widehat\Lambda)^{W'}=\mathbb C(\rho(\Lambda_0))(\sigma_1,\ldots,\sigma_r).$
\end{center}
This is the function field of the torus with character lattice $\rho(\Lambda_0)\oplus\mathbb Z^r$, whose monomials are $e^b\prod_i\sigma_i^{n_i}$.

Now, we see that the finite group $Z=\Hom(\widehat\Lambda/\Lambda_0,\mathbb C^*)$ acts on $\mathbb C(\widehat\Lambda)$ by $z(e^x)=z(x)e^x$. Here, it commutes with $W'$ and its invariant field is $\mathbb C(\Lambda_0)$. Since $W'\lambda_i\subset\lambda_i+Q$, we get that each monomial $e^b\prod_i\sigma_i^{n_i}$ is an eigenvector of $Z$ with character $z\mapsto z(b+\sum_in_i\lambda_i)$. Hence, the $Z$-invariants are spanned by the invariant monomials and we get
\begin{center}
$\mathbb C(\Lambda_0)^{W'}=\bigl(\mathbb C(\widehat\Lambda)^{W'}\bigr)^{Z}=\mathbb C(N),$
\end{center}
where
\begin{center}
$N=\Bigl\{(b,n)\in\rho(\Lambda_0)\oplus\mathbb Z^r:\ b+\textstyle\sum_in_i\lambda_i\in\Lambda_0\Bigr\}.$
\end{center}
Thus $\mathbb C(\Lambda_0)^{W'}$ is the function field of the torus with character lattice $N$.\\
\\
\textbf{Step 3.} The group $H_\Upsilon$ has order at most two. It acts on $\Lambda_0$ and normalizes $W'$. Let $h\in H_\Upsilon$. Since $h$ maps $\Phi$ to itself and the Weyl group $W'$ acts simply transitively on the bases of $\Phi$, there is $w\in W'$ such that $wh$ preserves the base fixed in Step 2. The element $w$ acts trivially on $\mathbb C(\Lambda_0)^{W'}$. The element $wh$ commutes with $\rho$ and $\pi$. It preserves $\Lambda_0$, $\rho(\Lambda_0)$ and $P_\Phi$, and it permutes the $\lambda_i$. Hence it permutes the $\sigma_i$, and it acts on $N$ by a lattice automorphism. So $H_\Upsilon$ acts monomially on the torus of Step 2, and it acts by affine maps between the fibres of the bundle. Lemma \ref{lemaffine} gives the rationality. The group $\prod_i S_{n_i}$ forgets the order of the boundary points. Hence the field we get is
\begin{center}
$\mathbb C(\widetilde I^{\mathrm{st}}_\Upsilon)^{G_\Upsilon} =\left(\left(\mathbb C(E_\Upsilon^{\mathrm{ord}})^{T_\Upsilon} \right)^{W'}\right)^{H_\Upsilon}.$
\end{center}
By Proposition \ref{propquot}, this field is the function field of the unmarked moduli space.
\end{proof}
As in the introduction, let $M_0$ be an irreducible component of the moduli space of smooth rational pairs $(X,D)$ with an ample polarization $C$. Here $D$ is a normal crossing anticanonical boundary with smooth components and a $0$-stratum.
\begin{theo}\label{thrsmooth}
If the general member of $M_0$ has $|H_\Upsilon|\leq2$, then $M_0$ is rational.
\end{theo}
\begin{proof}
Notice that the smooth pairs and the stable members of the complete linear system form a nonempty open subset of the same deformation component. Hence, the theorem follows from Proposition \ref{propuprat}.
\end{proof}
\begin{proof}[Proof of Theorem \ref{thrrat}]
This follows from Theorem \ref{thrsmooth} applying to the deformation type and to the complete polarize linear system.
\end{proof}
We give some examples for the condition $|H_\Upsilon|\leq2$. It holds for many deformation types from non-toric taut pairs and also for many toric types. For them, we need that the labels, the self-intersections and the polarization degrees tell the boundary components apart. Here, we only bound the order of $H_\Upsilon=F_\Upsilon/W_\Upsilon$. We do not need a bound on $|W_\Upsilon|$. Notice that we have $|H_\Upsilon|>2$ for the pairs from $\mathbb P^2$, $\mathbb P^1\times\mathbb P^1$ and cubic surfaces. The same also holds for their corner blow-ups. We treat these cases on a different paper.

Now, let $M_0$ be an irreducible component of the coarse KSBA moduli space $M(Y,D,L)$ of \cite{AAB}. We assume that a general member of $M_0$ is a smooth rational surface with a normal crossing anticanonical cycle. This cycle has smooth components and a $0$-stratum.
\begin{coro}\label{coraab}
If $|H_\Upsilon|\leq2$, then $M_0$ is rational.
\end{coro}
\begin{proof}
The smooth locus is dense in $M_0$. We have that it is also an open part of the HKY moduli component of Theorem \ref{thrsmooth}. So $M_0$ has the same function field as this component.
\end{proof}

\begin{rema}\label{remrigid}
Let $q=1$. We take the smooth members $C_1\in|aC_\infty|$ which meet $C_\infty$ transversely at $ae$ distinct points. The cluster $\mathcal K$ consists of these points. Each of them has multiplicity one, and there are no infinitely near centres. We fix a marking $(C_\infty,\mathcal O_{C_\infty}(C_\infty))\simeq(\Gamma,A)$. In this situation, the group is $\Aut_\Gamma(X_{\min},C_\infty)=T$. By Proposition \ref{propfixel}, the quotient by $T$ is birational to $|A^a|\times\mathbb P^{g-1}$. Hence it is rational.
\end{rema}

\subsection{Unirationality of the elliptic construction}\label{sec54}

In this part, we will study the case when $q=1$. Let $\mathcal A=\mathcal O_{C_\infty}(C_\infty)$. First, We take all reduced divisors $T\in|\mathcal A|$ and $Z\in|\mathcal A^{a-1}|$ with disjoint supports. Then we take all smooth members $C_1\in|aC_\infty|$ with $C_1|_{C_\infty}=T+Z$. The cluster $\mathcal K$ consists of these distinct points with multiplicity one. We take the pencil with its base divisor $T$. Let $E$ be a smooth fibre of this pencil such that $E\cap C$ is reduced. Finally, we apply Construction \ref{constr3} with this $E$.
\begin{prop}\label{propeluni}
If the admissible locus is nonempty, then the moduli space $\mathrm{CHS}^{\Upsilon_\ell}_\ell$ is unirational.
\end{prop}
\begin{proof}
First, we normalize the line bundle by a translatio  as in Construction \ref{constr5}. From the positive section, the base divisor of the pencil and the lifting sections, we get the rational parameter cover of Proposition \ref{propellpf}. On this cover, the relative pencil is a projective bundle with fibre $\mathbb P^1$. By the fact that smooth fibres $E$ with property $E\cap C$ a reduced divisor of degree $N_{\mathrm{ell}}$ form a nonempty open subset of this bundle, we get that this open subset is rational. After universal blow-up and birational contraction, we get a dominant map to the construction image of the moduli space. Hence, we get the unirationality.
\end{proof}
\begin{exam}\label{extang}
We cannot drop the condition that $C_1|_{W_1}$ consists of distinct transverse intersection points from the unirationality statement. We take $q=1$, $a=2$, $e=33$ and $R=0$. Let $p_1,\ldots,p_{11}$ be distinct points on a smooth elliptic curve $E$. We put
\begin{center}
$Z=\sum_{i=1}^{11}i p_i.$
\end{center}
Multiplication by two is surjective on $\Pic^0(E)$. So we can choose $A\in\Pic^{33}(E)$ with $A^{\otimes2}\simeq\mathcal O_E(Z)$. On $\mathbb P_E(\mathcal O_E\oplus A^{-1})$, the restriction map from $H^0(\mathcal O(2C_\infty))$ to $H^0(\mathcal O_{C_\infty}(2C_\infty))$ is surjective. The bundle $A$ is base-point-free. So a general lift of $Z$ is smooth at the points $p_i$, and it meets $C_\infty$ at $p_i$ with intersection multiplicity $i$. We resolve the boundary divisor $Z$. We get $C^2=66$ and $g(C)=34$. The resolved ruling has eleven singular fibres. Over $p_i$, the singular fibre has two $(-1)$-curves at its ends, of $C$-degree one, and $i-1$ intermediate $(-2)$-curves of $C$-degree zero. Hence the contraction defined by $C$ contracts no $(-1)$-curve. So the resolved surface is the minimal resolution of the constructed surface. We denote this resolution by $V$. The induced projection $\alpha:V\to E$ is its Albanese morphism and it gives us back $E$ and the $p_i$. Hence this moduli stratum dominates the moduli space $M_{1,11}$ of eleven-pointed smooth elliptic curves. Notice that a smooth compactification of $M_{1,11}$ has a nonzero holomorphic $11$-form which comes from the weight-$12$ discriminant cusp form \cite[Section 2.3]{FaberPandharipande13}. Hence $M_{1,11}$ is not unirational, and so this stratum is not unirational.
\end{exam}

\subsection{Rationality after forgetting the positive polarization}
\label{sec55}

In this subsection, we keep the initial ruled surface and the labelled boundary component in the objects and in the forgetful maps. Here $G_\Upsilon$ is the group of automorphisms which preserve the surface and the boundary component. We resolve the reduced base divisor $T$ of the pencil as in Construction \ref{constr2}. We denote by $\mathcal E_{\mathcal P}=\sum_iE_i$ the sum of its exceptional sections. These curves are sections of the fibration defined by the resolved pencil. In Construction \ref{constr3}, we keep their strict transforms as negative sections.
\begin{defi}
Let $\Upsilon_F^{\mathrm{neg}}$ be the tuple which keeps the ruled surface, the pencil and $\mathcal E_{\mathcal P}$. It forgets $C_1$. We denote by $\mathfrak B_{\Upsilon_F^{\mathrm{neg}}}$ the space of admissible base divisors $T$, and we define the birational moduli space of $\Upsilon_F^{\mathrm{neg}}$ by $\mathcal M_{\Upsilon_F^{\mathrm{neg}}} :=\mathfrak B_{\Upsilon_F^{\mathrm{neg}}}/\!/_{\mathrm{rat}}G_\Upsilon$.
\end{defi}
For a fixed labelled boundary component $E=C_\infty$, this base is
\begin{center}
$\mathfrak B_{\Upsilon_F^{\mathrm{neg}}}=
\begin{cases}
U_T\subset\Sym^e(E)\simeq\mathbb P^e,&q=0,\\
U_T\subset|A|\simeq\mathbb P^{e-1},&q=1,
\end{cases}$
\end{center}
Here $A=\mathcal O_E(C_\infty)$, and $U_T$ is the open subset where $T$ is reduced and satisfies the transversality condition. Notice that $T$ determines $\mathcal E_{\mathcal P}$ through the universal blow-up. This is equivariant for $G_\Upsilon$. Thus if we forget the choices of the positive sections, we get the rational map
\begin{center}
$\For_C:\mathrm{LCYF}_{\mathrm{lc}}\Upsilon_F \dashrightarrow\mathcal M_{\Upsilon_F^{\mathrm{neg}}}.$
\end{center}
For example, over a fixed elliptic initial pair, the space of tuples is $\{(C_1,T,R):T\in|A|,\ R\in|A^{a-1}|,\ C_1|_E=T+R\}$. In this case, the map $\For_C$ sends $(C_1,T,R)$ to $T$. 

In the following proposition, we study the rationality and unirationality of $M_{\Upsilon_F^{\mathrm{neg}}}$ with a given stratum $\Upsilon_F$.
\begin{prop}\label{propneg}
The moduli space of fibration $\mathcal M_{\Upsilon^{\mathrm{neg}}_F}$ with negative sections is unirational. For $q=0$ and non-smooth boundary, it is rational.
\end{prop}

\begin{proof}
First, we consider the case $q=0$. The diagonalizable part of $G_\Upsilon$ acts diagonally on $H^0(E,\mathcal O_E(e))$. The remaining group for the fixed labelled component has order at most two. It acts monomially on the invariant torus. By Theorem \ref{thrbohn} and Lemma \ref{lemmoninv}, we get that the quotient is rational. Now we come to the case $q=1$. For a fixed elliptic curve and a fixed line bundle, the linear system $|A|$ is a projective space. Thus its quotient is unirational. When the elliptic base varies, we use the translation of $\Gamma$ which gives $L\simeq\mathcal O_\Gamma(-eO)$ in Construction \ref{constr5}. This gives the same rational parameter cover.

For Construction \ref{constr3} with $q=1$, we add one smooth fibre of the pencil on this rational parameter cover. Notice that the choices of this fibre form an open subset of a projective bundle with fibre $\mathbb P^1$. Hence the space we get is rational, and it dominates the moduli image.
\end{proof}
For Construction \ref{constr3}, we forget the positive member, and we get
\begin{center}
$ \mathrm{CHS}^{\mathrm{ell}}_{\Upsilon_{\mathrm{ell}}} \dashrightarrow\mathcal M^{\mathrm{ell}}_{\Upsilon^{\mathrm{neg}}_{\mathrm{ell}}}.$
\end{center}
This map is dominant. Its general fibre consists of the choices of the positive system which we forget. We take them modulo the stabilizer of the ruled surface, the chosen pencil and the exceptional sections.

\section{Effectiveness for moduli spaces of stable minimal models}\label{sec6}
In this section, we give a partial answer to the effectiveness question of \cite[Question 11.13]{Bir2} for
\begin{center}
$\mathcal{M}_{slc}(2,\Phi_c=c\mathbb{Z}^{\geq 0},\mathcal V_{\mathrm{fib}},\sigma(t)=2(K_X+B)At+A^2t^2)$.
\end{center}

\subsection{The base curve of the Iitaka fibration}
\label{sec61} We assume that $\phi_{|n(K_S+B)|}:S\longrightarrow Z$ has one-dimensional fibres. So a general fibre $F$ is an effective Cartier divisor, and $\mathcal V_{\mathrm{fib}}$ is non-empty.
\begin{lemm}
Let $S$ be a normal surface. Let $D$ be a $\mathbb Q$-Cartier prime divisor, and let $C$ be a Cartier divisor. If $|C|$ is non-empty and base-point free, then $D\cdot C\in\mathbb Z$.
\end{lemm}
In the next proposition, we take an slc surface as above. We normalize it, and we take an irreducible component of the normalization. In this way, we get an lc pair $(Y,\Delta)$. Here the conductor is included in $\Delta$. We can assume that $\kappa(K_Y+\Delta)=1$.
\begin{prop}\label{propiit}
Let $c\in\mathbb Q_{>0}$ be fixed. Then there is $m(c)\in\mathbb Z_{>0}$ such that for every such pair $(Y,\Delta)$, the system $|\lfloor m(c)(K_Y+\Delta)\rfloor|$ is birational to the Iitaka fibration.
\end{prop}
\begin{proof}
 We have that the coefficients of $\Delta$ lie in the fixed finite set $(\Phi_c\cap[0,1])\cup\{1\}$ as the normalization adds the conductor with coefficient one. Now, the Proposition follows from \cite[Corollary 1.5]{CHL23}.
\end{proof}
Next are some important lemmas we will use later/.
\begin{lemm}\label{lemindex}
Let $\mathcal F$ be a bounded collection of surface stable minimal models with $\Q$-Cartier $A$. Then there is an integer $r=r(\mathcal F)>0$ such that every $rA$ is Cartier and $A^2\in r^{-2}\Z$ along $\mathcal F$.
\end{lemm}
\begin{proof}
We cover the bounded collection by finitely many projective families $\mathcal X_i\longrightarrow T_i$ with relative $\Q$-Cartier divisors $\mathcal A_i$. We divide each base into finitely many strata. After this, we can assume that the Cartier indices are bounded. Let $r$ be the least common multiple of these indices. Then $rA$ is Cartier on every fibre. So we get $A^2=(rA)^2/r^2\in\frac1{r^2}\Z$.
\end{proof}
\begin{lemm}\label{lemdenom}
Let $(S,B),A\in M_{\slc}(2,c\NN,\mathcal V_{\mathrm{fib}},\sigma(t))$. Then there is an integer $N>0$ such that $A^2\in N^{-1}\Z$. The integer $N$ depends only on the bounded family given by the fixed tuple of invariants.
\end{lemm}
\begin{proof}
These stable minimal models form a bounded family. Hence we can apply Lemma \ref{lemindex}, and we can take $N:=r^2$.
\end{proof}
For the next theorem, we fix $b>0$ and an integer $e>0$. We consider surface stable minimal models $(X,B),A\in S_{\mathrm{slc}}(2,\mathbb Z^{\geq0},\{e\},2b t+\alpha t^2)$ such that $A$ is $\mathbb Q$-Cartier.
\begin{theo}\label{thrbound}
There is a constant $C(b,e)$ which depends only on $b$ and $e$ such that $\alpha=A^2\geq C(b,e)$ for every such stable minimal model.
\end{theo}
Let $(X,B+tA)$ be slc for all sufficiently small $t>0$. We denote by $\nu:\bar X\to X$ the normalization. We denote by $\bar D$ its conductor divisor. We denote by $\bar B$ and $\bar A$ the boundary and the polarization on $\bar X$.
\begin{lemm}\label{lemnorm}
The normalized pair $(\bar X,\bar B+\bar D+t\bar A)$ is lc. The intersection numbers of the pullbacks of $\Q$-Cartier divisors do not change under the normalization. They also do not change under the $\Q$-factorial dlt modifications which we take after it.
\end{lemm}
We have
\begin{center}
$K_{\bar X}+\bar B+\bar D+t\bar A=\nu^*(K_X+B+tA).$
\end{center}
We also have $(\nu^*A)^2=A^2$ and $\nu^*(K_X+B)\cdot\nu^*A=(K_X+B)\cdot A$. The same equalities hold after a modification of each normal component. Let $\pi:Y\to\bar X$ be a dlt modification with boundary $\Delta_Y$. We have
\begin{center}
$K_Y+\Delta_Y =\pi^*(K_{\bar X}+\bar B+\bar D+t\bar A).$
\end{center}
Let $(Y,\Delta)$ be a normal lc surface pair with a boundary $\Delta\geq0$. We assume that $L=K_Y+\Delta$ is nef and semiample. We consider the fibration defined by $L$.
\begin{lemm}\label{lemadj}
Let $C\not\subset\Supp\Delta$ be an irreducible curve which is horizontal for this fibration. Then $C^2$ is bounded below by a constant which depends only on $L\cdot C$.
\end{lemm}

\begin{proof}
After the dlt reduction, $Y$ is $\mathbb Q$-factorial, and $C$ is $\mathbb Q$-Cartier. So by adjunction on $C^\nu$, we get an effective different, and we have $(K_Y+C)|_{C^\nu}=K_{C^\nu}+\Diff_{C^\nu}(0)$. Hence we get
\begin{center}
$(K_Y+C)\cdot C =2g(C^\nu)-2+\deg\Diff_{C^\nu}(0) \geq -2, \ \text{that is,} \
C^2\geq -2-K_Y\cdot C.$
\end{center}
We have $K_Y=L-\Delta$ and $C\not\subset\Supp\Delta$. So we have $\Delta\cdot C\geq0$, and thus $C^2\geq-2-L\cdot C$.
\end{proof}

\begin{lemm}\label{lemvert}
Let $f:X\to Z$ be the contraction defined by $L=K_X+B$. Then every contracted irreducible curve $V$ satisfies $A\cdot V\geq0$. Moreover, if $K_X+B+tA$ is ample, then $A\cdot V>0$.
\end{lemm}

\begin{proof}
The curve $V$ is contracted by $f$. So we have $L\cdot V=0$. If the divisor $L+tA$ is ample for sufficiently small $t>0$, then we get $0<(L+tA)\cdot V=tA\cdot V$. On a dlt modification, the pullback of $L+tA$ is nef, but it need not be ample. This gives the non-strict inequality for the exceptional curves.
\end{proof}
\begin{proof}[Proof of Theorem \ref{thrbound}]
We put $L=K_X+B$. By comparing the volume polynomials, we get $L\cdot A=b$ and $A^2=\alpha$. First, we suppose that $X$ is normal. By Lemma \ref{lemnorm}, we can take a $\mathbb Q$-factorial dlt modification, and we can pull back $L$ and $A$ without changing these intersection numbers. We write $A=A_{\mathrm h}+A_{\mathrm v}$ for the decomposition of $A$ into its horizontal part and its vertical part. Here $A_{\mathrm h}=\sum_i m_iC_i$ with $m_i\in\mathbb Z_{>0}$. The boundary is integral, and $(X,B+tA)$ is lc for small $t>0$. Hence no curve $C_i$ is a boundary component. We have $\sum_i m_i(C_i\cdot F)=e$. Hence we get $\sum_i m_i\leq e$ and $m_i\leq e$. We observe that distinct curves have nonnegative intersection numbers. By Lemmas \ref{lemadj} and \ref{lemvert}, we get
\begin{center}
$A^2\geq A_{\mathrm h}^2\geq -2\sum_i m_i^2-\sum_i m_i^2 L\cdot C_i \geq-2e^2-be.$
\end{center}
If the horizontal part of $A$ is reduced, the same estimate gives $A^2\geq-2e-b$ and hence for $e=b$, we get $A^2\geq -3e$. For $L\equiv F$, we get $A^2\geq-3e^2$ when $A$ can be non-reduced.

Now let $X$ be an slc surface. First, we normalize this surface. Then we take a dlt modification of each component. For each component which maps onto a curve, we take the Stein factorization $f_j:Y_j\to Z_j$. Let $L_j,A_j$ be the pullbacks of $K_X+B,A$. We put $e_j=A_j\cdot F_j$ and $b_j=L_j\cdot A_j$. The condition on the fibre volume gives $1\leq e_j\leq e$. We also have $\sum_j b_j=b$. The computation in the normal case gives $A_j^2\geq-2e_j^2-b_je_j$. By \cite[Corollary 1.5]{CHL23}, we can choose a uniform $m$ for the lc coefficients $\{0,1\}$ such that the moving part of $|mL_j|$ gives the Iitaka fibration. We take a resolution of this system, and we keep the same notation. We denote the divisors in the decomposition by $H_j,E_j$. So we have $mL_j\sim_{\mathbb Q}f_j^*H_j+E_j$ with $\deg H_j\geq1$, where $E_j$ is effective and vertical. The divisor $A_j$ is nef on vertical curves. So we get $mb_j\geq e_j$. Hence we get $\sum_j e_j\leq mb$. On the components which map to points, the polarization is ample before the dlt modification. Hence it has positive square. Thus we get
\begin{center}
$ A^2\geq-2\sum_j e_j^2-\sum_j b_je_j \geq-(2m+1)be.$
\end{center}
Hence we can take $C(b,e)=-(2m+1)be$.
\end{proof}

\begin{proof}[Proof of Theorem \ref{thrlower}]
Theorem \ref{thrbound} gives the lower bound for fixed $b$ and $e$. For a normal surface, its proof gives $A^2\geq-2e^2-be$. When $A$ is reduced, it gives $A^2\geq-2e-b$. If $K_X+B\equiv F$, then we have $b=(K_X+B)\cdot A=e$. This gives the bounds in the statement.
\end{proof}

\begin{rema}
\label{remcoeff} Let $c\in\mathbb Q_{>1}$, and let $(X,B),A\in\Sslc(2,\Phi_c,\{\gamma\},2b t+\alpha t^2)$ with $A$ $\mathbb Q$-Cartier. We assume that $X$ is normal and that the base of the Iitaka fibration is a curve. An slc boundary has coefficients at most one. Every nonzero element of $\Phi_c$ is larger than one. Hence we have $B=0$. Let $N$ be the effective integral Mumford divisor with $A=cN$. We get $N\cdot F=\gamma/c$. In particular, the space is nonempty only if $\gamma/c\in\mathbb Z_{>0}$. By rescaling $t$, we get the stable minimal model conditions for $N$ from the conditions for $A$. We have
\begin{center}
$\vol(K_X+tN) = \frac{2b}{c}t+\frac{\alpha}{c^2}t^2 .$
\end{center}
We apply the proof of Theorem \ref{thrbound} to $(X,0),N$. We get
\begin{center}
$\alpha\ge -2\gamma^2-b\gamma .$
\end{center}
If $K_X\equiv\lambda F$, then we get $b=\lambda\gamma$ and $\alpha\ge-(2+\lambda)\gamma^2$. For $K_X\equiv F$ and $\gamma=ce$, this bound becomes $\alpha\ge-3c^2e^2$. If the horizontal part of $N$ is reduced, the reduced bound gives $\alpha\ge-2c\gamma-cb$. For $K_X\equiv F$ and $\gamma=ce$, this bound becomes $\alpha\ge-3c^2e$. Finally, we scale an integral example with polynomial $2b_0t+\alpha_0t^2$. In this way, we get an example with polynomial $2cb_0t+c^2\alpha_0t^2$.
\end{rema}
In the following examples, we put $d:=A\cdot F\geq1$.
\begin{exam}\label{exell}
Let $\Gamma$ be a smooth elliptic curve. Let $D$ be a divisor on $\Gamma$ of degree $d\geq2$. We put $X_{\min}:=\mathbb P_\Gamma(\mathcal O_\Gamma\oplus\OO_\Gamma(-D))$. The ruling $\xi:X_{\min}\longrightarrow\Gamma$ has a negative section $C_0$ and a positive section $C_\infty$. These two sections are disjoint. We have $C_0^2=-d$, $C_\infty^2=d$ and $-K_{X_{\min}}\sim C_0+C_\infty$. In Construction \ref{constr1}, we take $a=2$, $e=d$ and $q=1$. We put $W=C_0+C_\infty$. We write $\mathcal A$ for $\mathcal O_{C_\infty}(C_\infty)$. Since $d\geq2$, the linear system $|\mathcal A|$ is base-point-free. By Proposition \ref{propdblpt}, we can choose a general reduced divisor $T=t_1+\cdots+t_d\in|\mathcal A|$. We can also choose a general smooth $C_1\in|2C_\infty|$ with $C_1|_{C_\infty}=2T$. Now we take the positive-section pencil with base divisor $T$. By Construction \ref{constr2} and Theorem \ref{thrm2crv}, it gives a surface $S$ with CHS and a fibration $\varphi:S\to\mathbb{P}^1$. The $d$ disjoint $(-2)$-sections $R_i$ survive on $S$. For a general fibre $F$, we have $F^2=0$ and $R_i\cdot F=1$.

The idea is that we use Construction \ref{constr3} with $R_1+\cdots+R_d$ as the polarize divisor $C$. First, we choose a smooth fibre $F_0$ and then blow-up of the $d$ points $R_i\cap F_0$ to get $\mu:Y\to S$. Let $A_i$ and $\widetilde F_0$ be the strict transforms of $R_i$ and $F_0$ respectively. Denote $A_Y$ by the  sum $\sum_i A_i$. We have that $A_i^2=-3$, $\widetilde{F}_0^2=-d$ and $K_Y+\widetilde F_0\sim\mu^*F_0$. Take a small $t>0$. We have that log divisor $\mu^*F_0+tA_Y$ has the intersection number $1-3t$ on each $A_i$. Every other vertical component meets some $A_i$ and every other horizontal curve has positive fibre intersection number. So this divisor is nef and big. The curve $\widetilde F_0$ is the only curve on which it has degree zero. The fibre and the sections are smooth and disjoint. Hence $(Y,\widetilde{F}_0+tA_Y)$ is log canonical for small $t>0$. By surface log abundance \cite[Corollary 6.6]{FujinoSurfaces12}, the log divisor is semiample. The morphism which it defines contracts $\widetilde F_0$ to a simple elliptic point. We denote by $A$ the image of $A_Y$. In this way, we get a stable minimal model $(S_{\mathrm{ell}},0),A$ with $K_{S_{\mathrm{ell}}}\sim F$, $A\cdot F=d$ and $A^2=-3d$. Hence
\begin{center}
$\vol(K_{S_{\mathrm{ell}}}+tA) = (K_{S_{\mathrm{ell}}}+tA)^2 = 2dt-3dt^2.$
\end{center}
\end{exam}
\begin{exam}\label{exruled}
Let $d\geq1$. We take the smooth rational surfaces of Construction \ref{constr4}. They have $K_Y+B_Y\sim F$. The polarization $A$ is a union of $d$ disjoint rational $(-3)$-sections. So we have $A\cdot F=d$ and $A^2=-3d$. This gives the polynomial $2dt-3dt^2$. Here is how the $(-3)$-sections come up. The exceptional $(-1)$-sections become $(-2)$-sections after the blow-ups tangent to the boundary. Then we blow up their intersections with one more pencil fibre. After this, they are $(-3)$-sections. Finally, we contract only the vertical curves of polarization degree zero. The horizontal boundary components are not contracted.
\end{exam}
\begin{exam}\label{exsharp}
Let $\Gamma$ be a smooth elliptic curve. Let $D$ be a line bundle on $\Gamma$ of degree one. We put
\begin{center}
$X_{\min}:=\mathbb P_\Gamma\bigl(\mathcal O_\Gamma\oplus\mathcal O_\Gamma(-D)\bigr) \ \text{and} \
\xi:X_{\min}\longrightarrow\Gamma .$
\end{center}
Let $C_0$ and $C_\infty$ be the negative and the positive sections. We have
\begin{center}
$C_0^2=-1,\ C_\infty^2=1,\ C_0\cap C_\infty=\emptyset \ \text{and} \ -K_{X_{\min}}\sim C_0+C_\infty .$
\end{center}
We take the tuple from Construction \ref{constr1} with $(a=4,\ e=1, \ q= 1) .$

We put $\mathcal A:=\mathcal O_{C_\infty}(C_\infty)$. This is a line bundle of degree one on the elliptic curve $C_\infty$. Let $p_0$ be the unique effective divisor in $|\mathcal A|$.

We have
\begin{center}
$\xi_*\mathcal O_{X_{\min}}(4C_\infty)\;\simeq\; \bigoplus_{i=0}^{4}\mathcal A^{i}.$
\end{center}
The restriction to $C_\infty$ is the projection onto the last summand. So it is surjective. We also have
\begin{center}
$L(-B)=\mathcal O_{X_{\min}}\bigl(4C_\infty-C_\infty\bigr) =\mathcal O_{X_{\min}}(3C_\infty).$
\end{center}
Its restriction to $C_\infty$ is $\mathcal A^{3}$, which has degree three. Since $e=1$, the line bundle $\mathcal A=\mathcal O_{C_\infty}(C_\infty)$ has degree one. So $|\mathcal A|$ consists of the single point $p_0$. Notice that a line bundle of degree at least two on an elliptic curve is base-point-free. Hence the normal first-order condition is independent at every point of $C_\infty$, also at $p_0$. Moreover, $N(-2p_0)=\mathcal A^{2}$ has degree two. Hence we can choose a section which does not vanish at $p_0$. So Proposition \ref{propdblpt} applies. For $a=3$, the condition $R_3\in|\mathcal A|$ would force $R_3=T$. For $a=4$, we take
\begin{center}
$T:=p_0\in|\mathcal A|,\ R_4\in|\mathcal A^{2}|\
\text{general}$ and $Z_4:=2T+R_4\in|\mathcal A^{4}| .$
\end{center}
Since $\deg\mathcal A^{2}=2$, a general member $R_4$ of $|\mathcal A^{2}|$ is reduced and disjoint from $T$. Consider the system $\Lambda_{Z_4}\subset|4C_\infty|$. By Proposition \ref{propdblpt}, it has the two base centres $(p_0,1)$, $(p_0',1)$ with $p_0'\succ p_0$and two other base points of multiplicity one. Hence we get
\begin{center}
$\sum_\alpha m_\alpha=2+2=4=ae .$
\end{center}
By resolving the cluster, we get $\varphi:X\longrightarrow X_{\min}$ with
\begin{center}
$C^2=a^2e-ae=12.$
\end{center}
The resulting surface $X$ contains one $(-2)$-curve $R$, which is the strict transform of the first exceptional divisor induced from blowing-up $p_0$. Since $e=1$, we have $h^0(X_{\min},\mathcal O(C_\infty))=1+h^0(\mathcal A)=2$. Hence $|C_\infty|$ is a pencil. Its base scheme is the point $p_0$. We take $\mathcal P:=|C_\infty|$ in Construction \ref{constr2}. We denote by $S$ the surface with CHS which we get. The base point of $\mathcal P$ is the support of the length-two cluster. This gives $R\cdot F=1$. Thus $R$ is a section of the induced elliptic fibration.

Now we use Construction \ref{constr3} with $R$ in place of the CHS divisor $C$. We choose a smooth fibre $F_0$. We denote by $\mu:Y\to S$ the blow-up of the point $R\cap F_0$. Let $R_Y$ and $\widetilde F_0$ be the strict transforms. We have $R_Y^2=-3$, $\widetilde F_0^2=-1$ and $K_Y+\widetilde F_0\sim\mu^*F_0$. Let $m\geq1$ be an integer. We put $A_Y=mR_Y$. Notice that the pair $(Y,\widetilde F_0+tA_Y)$ is log canonical for small $t>0$. Its log divisor has nonnegative intersection numbers with vertical curves. Its degree on $R_Y$ is $1-3mt>0$. Every other horizontal curve has positive fibre degree. So the log divisor is nef. Its square is $2mt-3m^2t^2>0$. By surface log abundance \cite[Corollary 6.6]{FujinoSurfaces12}, we can take the contraction $\pi:Y\to S_\ell$ defined by its section ring. This morphism contracts $\widetilde F_0$. It also contracts the other vertical curves of $A_Y$-degree zero. Among them are the curves over the two proper points of $R_4$. The image of $\widetilde F_0$ is a simple elliptic point of degree one. The curve $R_Y$ does not meet the contracted curves. So its image $R_\ell$ is smooth, and it does not contain the singular points. We also have $R_\ell^2=-3$ and $R_\ell\cdot F=1$. By this way, we get a stable minimal model $(S_\ell,0),A$ with $K_{S_\ell}\sim F$, $A\cdot F=m$ and $A^2=-3m^2$. Here, $A:=mR_{\ell}$.  Hence, we have
\begin{center}
$\vol(K_{S_\ell}+tA)=2m\,t-3m^2t^2 .$
\end{center}
In the notation of Theorem \ref{thrbound}, the fibre degree is $e=m$, and we have $\mathcal V_{\mathrm{fib}}=\{m\}$. Thus, we get the sharp bound $A^2=-3m^2$. For $m>1$, the polarization $A=mR_\ell$ is not reduced. The multiplicity $m$ appears squared in $A^2=-3m^2$. This is why the bound is quadratic in $m$.
\end{exam}
We put
\begin{center}
$ \mathcal M_{\alpha,\beta,e}:= \mathcal M_{slc}(2,\mathbb Z^{\geq0},\{e\},2e\beta t+\alpha t^2).$
\end{center}
\begin{theo}\label{threffect}
The moduli space $\mathcal M_{\alpha,\beta,e}$ is non-empty for every $\beta\in\mathbb N$ and every integer $\alpha\geq-3e$.
\end{theo}
\begin{proof}
We want a stable minimal model $(X,B),A$ with a fibration $f:X\longrightarrow \mathbb P^1$ and general fibre $F$ such that $A\cdot F=e$, $A^2=\alpha$ and $K_X+B\sim_{\mathbb Q} \beta F$. This is enough. Indeed, we have $(K_X+B)\cdot A=\beta e$ and $\vol(K_X+B+tA)=(K_X+B+tA)^2=2e\beta t+\alpha t^2$.\\
\\
\textbf{Case 1.} Let $e=1$, $X=\mathbb F_m$ with $m\in\{1,2,3\}$ and write $A=C_0$. Now let us take the general smooth $B\in|2C_0+(m+\beta+2)f|$. This boundary system is base-point-free. We then have $K_X+B\sim\beta f$. For small $t>0$, the divisor $\beta f+tC_0$ is ample, since $\beta-mt>0$ and the pair is log canonical. We have $A\cdot f=1$ and $A^2=-m$. By adding $k$ general fibres to $A$, we have $A^2=-m+2k$. Hence, these pairs cover every $\alpha\geq-3$. 
\\
\\
\textbf{Case 2.} Assume $e\geq2$. We first consider $-3e\leq\alpha\leq-e$. We start with the rational boundary construction of Example \ref{exruled} and Construction \ref{constr4}. We take $d$ distinct proper pencil base points on the labelled boundary. At these points, we use the length-two clusters of Subsection \ref{sec34}. As in the proof of Theorem \ref{thrm2crv}, they give disjoint sections $R_i$ with $R_i^2=-2$ and $R_i\cdot F=1$. We choose a smooth fibre $F_0\in|F|$ which meets the $R_i$ transversely. We blow up the points $q_i=R_i\cap F_0$. We take $F_0$ as a boundary component of coefficient one. Let $A_i$ be the strict transforms of the $R_i$. We have $A_i^2=-3$ and $A_i\cdot F=1$. Notice that we are free to choose the number of $(-1)$-sections from the blow-ups in Construction \ref{constr2}. The same holds for the $(-2)$-sections. So we can choose $x$ pairwise disjoint $(-1)$-sections $E_i$, $y$ pairwise disjoint $(-2)$-sections $R_j$ and $z$ pairwise disjoint $(-3)$-sections $A_\ell$ of the same fibration, with $x+y=e$ or $y+z=e$. We put
\begin{center}
$A^{(1)}=\sum_{i=1}^{x}E_i+\sum_{j=1}^{y}R_j \ \text{and} \
A^{(2)}=\sum_{j=1}^{y}R_j+\sum_{\ell=1}^{z}A_\ell .$
\end{center}
We have $A^{(i)}\cdot F=e$, $(A^{(1)})^2=-x-2y$ and $(A^{(2)})^2=-2y-3z$. We take $A^{(1)}$ for every integer $-2e\leq \alpha \leq-e$. We take $A^{(2)}$ for every integer $-3e\leq \alpha \leq -2e$. In this way, we get integral divisors $A_-$ with $A_-\cdot F=e$ and $A_-^2=-m$ for every integer $e\leq m\leq 3e$.\\
\\
\textbf{Case 3.} We now consider $\alpha>-e$. We write $\alpha=\alpha_0+2ek$ with integers $-3e\leq \alpha_0 \leq -e$ and $k\geq1$. By Case 2, we can take an integral divisor $A_0$ with $A_0\cdot F=e$ and $A_0^2=\alpha_0$. We choose general fibres $G_1,\ldots,G_k\in |F|$. We put $A=A_0+G_1+\cdots+G_k$. Then we have $A\cdot F=e$. Since the $G_i$ are fibres, we get $A^2=A_0^2+2k(A_0\cdot F)=\alpha$.

Finally, we choose an anticanonical boundary $B_0$ from the construction. Hence we have $K_X+B_0\sim0$. We put $B=B_0+H_1+\cdots+H_{\beta}$ with general $H_i\in |F|$. When we have chosen $(-3)$-sections, we take $H_1=F_0$ instead. Here $F_0$ is the fibre chosen above. Then all coefficients of $B$ lie in $\mathbb Z^{\geq 0}$. We have $K_X+B\sim_{\mathbb Q}\beta F$ and $\vol(A|_F)=A\cdot F=e$. So we get $\mathcal V_{\mathrm{fib}}=\{e\}$. This proves that the moduli space is non-empty.
\end{proof}

\begin{proof}[Proof of Theorem \ref{threxist}]
This is Theorem \ref{threffect}, since $\mathcal M_{\alpha,\beta,e} =M_{\mathrm{slc}}(2,\mathbb Z^{\geq0},\{e\},2e\beta t+\alpha t^2)$.
\end{proof}

For $c\in\mathbb Q_{>0}$, let $\Phi_c=c\mathbb Z^{\geq0}$. We consider the following coefficient sets and polynomials:
\begin{enumerate}
    \item For $c=1/k$ with $k\geq2$, we take $\mathcal V_{\mathrm{fib}}=\{e\}$ and $\sigma(t)=2e\beta t+\alpha t^2$, where $\beta\in\mathbb N$ and $\alpha\in\mathbb Z^{\geq-3e}$.
    \item For $c>1$, we take $\mathcal V_{\mathrm{fib}}=\{ce\}$ and $\sigma(t)=2ce\,t+\alpha c^2t^2$ with $\alpha\in\{-3e,-3e^2\}$, which come from Examples \ref{exell} and \ref{exsharp}.
\end{enumerate}
\begin{coro}\label{corcoeff}
For both sets of tuples, the space $M_{\mathrm{slc}}(2,\Phi_c,\mathcal V_{\mathrm{fib}},\sigma)$ is non-empty. For $c>1$, every $(2,\Phi_c)$-stable minimal model has zero boundary.
\end{coro}
\begin{proof}
For $c=1/k$, the set $\Phi_c$ contains $\mathbb Z^{\geq0}$ and Theorem \ref{threffect} applies. For $c>1$, Notice that an slc boundary has coefficients at most one. Every nonzero element of $\Phi_c$ is larger than one, hence the boundary is zero. Thus for the non-emptiness, we can take an example with zero boundary, i.e., $K_X\sim F$ and replace $A$ by $cA$. Then the intersection number of $A$ with the fibre becomes $ce$. The volume polynomial becomes $2ce\,t+c^2A^2t^2$ and we finish the proof.
\end{proof}

\subsection{Elliptic CHS components.}\label{sec63}
We let the elliptic curve $E$ vary. We choose $L\in\Pic^e(E)$. We put $Y_0=\mathbb P_E(\mathcal O_E\oplus L^{-1})$. We take the reduced base divisor $T\in|L|$ of the positive-section pencil on $Y_0$, and we resolve it. This gives $\mu:X\to Y_0$. We denote the exceptional sections by $R_i$. We put $F=\mu^*C_\infty-\sum R_i$ and $M_k=\sum R_i+kF$ for $k\geq0$. We take a general $\Theta_k\in|M_k|$ and $\beta+1$ distinct smooth pencil fibres $H_0=\widetilde C_\infty,H_1,\ldots,H_\beta$, with $\beta\geq1$. Starting with $\Delta=C_0+\sum H_j$, we blow up all intersections of $\Theta_k$ with $r-1$ of these fibres $H_i$. When $w>0$, we also blow up $w$ points on an unused fibre. Here we have that $r=1,2$ and $0\leq w<e$, or $(r,w)=(3,0)$. Finally, we contract $C_0$ and the strict transforms of the fibres whose intersections were all blown-up and get the image $\mathcal E_{e,\beta,r,w,k}$. All entries of the tuple vary, and so does the curve $E$. For $k=0$, this is the CHS construction with negative sections. For $k\geq1$, the polarization moves in the linear system which contains $A+kF$.
\begin{theo}\label{threllcmp}
For $e\geq10$, the closure of $\mathcal E_{e,\beta,r,w,k}$ is an irreducible component of $M_{\mathrm{slc}}(2,\mathbb Z_{\geq0},\{e\},2\beta et+(2ke-re-w)t^2)$. These components are unirational.
\end{theo}
\begin{proof}
First, we take a locally stable family over the spectrum $T$ of a DVR. Its closed fibre is a general member. Since $e\ge10$, the deformations of the simple elliptic points are equisingular by \cite[9.2(a)]{Pin74}. Hence we take a finite base change. Then we blow up the singular sections. This gives a simultaneous minimal resolution. Now we use the argument of Lemma \ref{lemrevMMP}. The $(-1)$-curves which meet the elliptic curves extend over $T$. So we can contract them relatively. In this way, we get the Albanese $\mathbb{P}^1$-bundle over a family $\mathcal{E}/T$ of elliptic curves. It has $\beta+2$ disjoint boundary sections. Three of these sections trivialize the bundle. Hence, it is isomorphic to $\mathcal{E}\times_T\mathbb{P}^1_T$ and its polarization lies in $|p_{\mathcal{E}}^*L\otimes p_{\mathbb{P}^1}^*\mathcal{O}(k)|$ with $L\in\Pic^e(\mathcal{E}/T)$. Here, the contracted curves determine the blown-up points. Finally, we repeat these blow-ups. Then we take the relative log canonical model of $K+B+\tau A$ for a small rational $\tau>0$ and get the recovery of the family. Hence, its general fibre lies in $\mathcal{E}_{e,\beta,r,w,k}$. Thus, the closure of $\mathcal{E}_{e,\beta,r,w,k}$ is a component.

We normalize $L\cong\mathcal{O}_E(eO)$ by translation, as in Construction \ref{constr5}. Here $E$ varies over the open $j$-line. On $E\times\mathbb{P}^1$, the divisor $\Theta_k$ varies in $|p_E^*L\otimes p_{\mathbb{P}^1}^*\mathcal{O}(k)|$. The $\beta+2$ boundary fibres vary in a rational configuration space. Let $w>0$. We blow up $w$ points on one boundary fibre $H\cong E$. These points form a divisor $D_w\le\Theta_k|_H$ of degree $w$. The choices of $D_w$ form a projective bundle over the relative $\Pic^w$. This relative $\Pic^w$ is the universal elliptic curve of Construction \ref{constr5}. This curve is rational. Indeed, on its affine chart we have $j=\frac{1728(y^2-x^3)}{y^2-x^3+3x+2}$. The divisors $\Theta_k$ with $\Theta_k|_H\ge D_w$ form a nonempty open subset of a projective bundle over it. Hence the parameter space is rational. Thus $\mathcal{E}_{e,\beta,r,w,k}$ is irreducible, and its closure is unirational.
\end{proof}

\end{document}